\documentclass[12pt,twoside]{article}

\usepackage[scale=0.75, left=1in, right=1in, asymmetric]{geometry}
\usepackage{xcolor}
\usepackage{amsthm, amsmath, amssymb, mathtools}
\usepackage{graphicx}
\usepackage{tikz}
\usetikzlibrary{calc,arrows.meta,positioning}
\usepackage{xspace} 
\usepackage{tabularx}
\usepackage[inline]{enumitem}
\usepackage[titles]{tocloft}
\usepackage{fancyhdr}
\usepackage{slantsc} 
\usepackage[style=base]{caption} 
\usepackage[
    style=alphabetic,
    maxbibnames=20,
    minbibnames=3,
    mincitenames=2,
    maxcitenames=3,
    maxalphanames=3,
    minalphanames=2,
    backend=biber]{biblatex}
\usepackage{float} 
\usepackage{hyperref} 
\definecolor{darkblue}{RGB}{0,0,160}
\hypersetup{colorlinks,
  citecolor=darkblue,
  filecolor=black,
  linkcolor=darkblue,
  urlcolor=darkblue,
  pdftitle={Topological data analysis of C. elegans locomotion and behavior},
  pdfauthor={Ashleigh Thomas, Kathleen Bates, Alex Elchesen, Iryna Hartsock, Hang Lu, and Peter Bubenik},
  pdfkeywords={persistent homology, C. elegans, behavior},
  pdfdisplaydoctitle=true,
}
\usepackage[noabbrev, nameinlink, capitalize]{cleveref} 
\crefname{subsection}{section}{sections}

\newtheorem{theorem}{Theorem}[section]
\theoremstyle{definition}
\newtheorem{remark}[theorem]{Remark}
\newtheorem{definition}[theorem]{Definition}
\newtheorem{example}[theorem]{Example}

\newcommand{\RR}{\mathbb{R}}
\newcommand{\NN}{\mathbb{N}}
\newcommand{\ZZ}{\mathbb{Z}}
\renewcommand{\tilde}[1]{\widetilde{#1}}
\newcommand{\defsty}[1]{\emph{#1}}

\newcommand*{\ie}{{i.e.}\@\xspace}

\newcommand{\wl}{l} 
\newcommand{\swe}[2]{#1^{#2}} 

\def\quarterwidth{0.215\textwidth}
\def\fifthwidth{0.17\textwidth}
\def\halfwidth{0.47\columnwidth}
\newcommand{\trimmedgraphic}[2][]{%
  \includegraphics[trim = 165 205 50 70,clip,#1]
  {#2}%
}
\newcommand{\trimmedwithaxesgraphic}[2][]{
  \includegraphics[trim = 50 80 50 55,clip,#1]%
  {#2}%
}
\newcommand{\hightrimmovespace}[2][]{%
  \includegraphics[trim=400 105 300 40,clip,#1]
  {#2}%
}

\makeatletter
\let\c@table\c@figure 
\let\ftype@table\ftype@figure 
\makeatother

\title{\vspace{-2em}Topological data analysis of \emph{C.~elegans}\\ locomotion and behavior}
\author{Ashleigh~Thomas, Kathleen~Bates, Alex~Elchesen$^\dagger$, Iryna~Hartsock$^\dagger$, \\Hang~Lu, and Peter~Bubenik}
\date{\vspace{-1em}\footnotesize{$^\dagger$ These authors contributed equally to this work.}}

\begin{document}

\maketitle
\vspace{-2em}
\begin{abstract}
  \noindent
  We apply topological data analysis to the behavior of \emph{C.\ elegans}, a widely-studied model organism in biology.
  In particular, we use topology to produce a quantitative summary of complex behavior which may be applied to high-throughput data.
  Our methods allow us to distinguish and classify videos from various environmental conditions and we analyze the trade-off between accuracy and interpretability.
  Furthermore, we present a novel technique for visualizing the outputs of our analysis in terms of the input.
  Specifically, we use representative cycles of persistent homology to produce synthetic videos of stereotypical behaviors.
  
  
\end{abstract}

\vspace{-1.5em}

\setlength\cftbeforesecskip{0pt}
{\parskip=0em \tableofcontents}

\section{Introduction}\label{sec:introduction}
Model organisms are indispensable in understanding basic principles of biology. Studies of model organisms have played a major role in discoveries of disease mechanisms, disease treatment, and neuroscience principles. The behavior of these model organisms can illuminate responses and phenotypes important for understanding the effects of experimental conditions on subjects. Behavior can be affected by neuron activity, external stimuli, and past experiences (learning), so being able to adequately measure and compare behaviors is a useful evaluation tool for a wide range of experiments.

We propose persistent homology as a new tool for assessing behavior of \emph{Caenorhabditis elegans}, worms that are a widely-used model organism. Persistence has been successfully used to study high-dimensional time series, especially those that exhibit some quasi-periodic behavior like the undulation of \emph{C. elegans} \cite{Tralie2016, Tralie2018}. But to the authors' knowledge, persistent homology has not been previously used to analyze \emph{C. elegans} behavior, though it and similar techniques have been used to study \emph{C. elegans} neural data \cite{Petri2013, Sizemore2019, Helm2020, Lutgehetmann2020, Backholm2015}.

In this paper we use persistent homology to study the locomotion of \emph{C. elegans} in two settings.
In our initial study (\Cref{sec:case-study}), we follow one worm as it moves on the surface of an agar plate.
Under these conditions there are no barriers to movement and the locomotion is both smooth and complex. We show that persistent homology is able to detect and differentiate between various characteristic behaviors such as forward crawling, backward crawling, and transitioning between the two. 
We also show a unique contribution of persistence: 
the synthesis of  skeleton data of \emph{C. elegans} performing stereotyped, periodic behaviors. This translates into videos of, for example, forward crawling that are smooth when looped (see the supplemental materials \Cref{sec:supp-material} for an example). 
Furthermore, 
this mapping from persistence features to behavior gives a concrete and biologically relevant interpretation of results: 
features of interest --- such as a feature that is detected in one sample and not another --- can be expressed as videos of synthetic behavior.

We also analyze data from a large-scale experiment of the effect of environment on \emph{C. elegans} (\Cref{sec:experimental-results}). In this setting a more controlled environment is required, so the organisms are submerged in a solution and confined to wells in microfluidic devices. Our main results study \emph{C. elegans}’ locomotion in  solutions that have various levels of viscosity. We show that we are able to use persistent homology --- and average persistence landscapes in particular --- to summarize \emph{C. elegans} locomotion in a way that allows the classification of the viscosity of an animal's environment with a high level of accuracy, and in fact a much higher level of accuracy than simpler methods based on speed and variety of postures. 
Our results indicate that persistent homology is a promising tool for quantifying the impact of changes to genotype and environment on \emph{C. elegans} locomotion.

\subsection{Related work}

\emph{Caenorhabditis elegans} is a free-living soil nematode that has been a workhorse genetic model system. The nematode’s transparent tissue, simple anatomy, and fast reproduction contribute to both ease in culture and a literal window into the internal workings of a living organism. Its completely sequenced genome contains many genes that are homologous to human genes, and importantly the ability to manipulate genes with relative ease makes it an extremely attractive model system. For neuroscience in particular, \emph{C. elegans} presents a unique opportunity with its simple nervous system (just 302 neurons) that is complex enough to exhibit many sensory modalities, including mechanosensation, chemosensation, and response to heat, osmolarity, and smell. 

Behavior characterization in \emph{C. elegans} was historically qualitative, mainly relying on experimentalists
specifying end-point assessment (e.g. whether the worm chemotaxes to a particular source of odor
within a certain amount of time), or experimentalists using heuristics to assess behavior (e.g. naming
worms genes ``\emph{unc}'' for uncoordinated). In the last decade, machine vision tools first replaced human
identifications of worms from images and videos, which allows much larger dynamic datasets to be
annotated and analyzed. In recent years, further development in quantitative behavior characterization tools such as 
tracking \cite{Yemini2013, Stirman2011, Swierczek2011, Husson2012, Porto2019}, 
eigenworms \cite{Stephens2008}, 
behavior ‘dictionaries’ \cite{Brown2013}, 
and t-SNE \cite{Berman2016, Liu2018} 
have moved the field away from merely describing the outcome to understanding the types of behavior the
brain of this simple system can generate. While many of these techniques do well in quantitatively
describing behavior and distinguishing differences in behavior, behavioral dynamics are rich and
opportunities abound in exploring behavioral dynamics using other mathematical tools.

Persistent homology has been used to analyze time series data in many different settings. Some earlier work was theoretical and studied the interaction between persistence and sliding window embeddings --- which we used in this research --- as well as proposed possible applications \cite{Perea2015b, Perea2016, Khasawneh2016}. 
Research into gene expression has used persistent homology to detect patterns or classify whether a signal is periodic \cite{Perea2015, Dequeant2008}. Frequently, persistence has been used to study neural data \cite{Stolz2017, Petri2013, Sizemore2019, Helm2020, Lutgehetmann2020, Backholm2015}, and in many cases neural data from \emph{C. elegans}, but the analysis tends to rely on clique complexes as the topological space of interest instead of sliding window embeddings.

\subsection{Acknowledgments}

The authors would like to thank referees Jose Perea and Samir Chowdhury, whose many comments considerably improved our manuscript. 
AT would also like to thank Kim Le for helpful conversations. 

This research was partially supported by the NSF-Simons Southeast Center for Mathematics and Biology (SCMB) through the grants National Science Foundation DMS1764406 and Simons Foundation/SFARI 594594.
This material is based upon work supported by, or in part by, the Army Research Laboratory and the Army Research Office under contract/grant number W911NF-18-1-0307. 
Research reported in this publication was supported in part by the National Institutes of Health under award numbers R01NS096581, R01NS115484, and R01AG056436.
The collection of the data used in this research was partially supported by the National Institutes of Health's Ruth L. Kirschstein NRSA award 1F31GM123662.

\section{Materials and methods}\label{sec:materialsandmethods}

In this section we describe the collection and preprocessing of experimental data (\Cref{sec:description-of-data}), mathematical background (\Cref{sec:swe,sec:persistent-homology}), and pipeline for using topological data analysis on \emph{C. elegans} behavior data (\Cref{sec:pipeline}).

\subsection{Description of data} \label{sec:description-of-data}
\emph{C. elegans} (N2 strain) were cultured at 20$^{\circ}$C under standard conditions on agar plates seeded with OP50
\emph{E. Coli}. Animals were age-synchronized via hatch-off and cultured on plate until they reached day 1 of
adulthood. For behavior experiments on agar, animals were prepared, imaged, and tracked as
previously described \cite{Porto2019}. For behavior experiments in 
methylcellulose media, synchronized populations were then washed off of culture plates with M9
buffer. Unless otherwise noted, video data was collected on a dissecting microscope (Leica MZ16) using
a CMOS camera (Thorlabs DCC3240M), with a frame rate of 30 frames per second and a magnification of 1.2x.

Behavior data was collected with animals confined with microfluidic devices. In these devices the cavities in
which worms are loaded have only slightly greater depth than the width of an adult
worm, which restricts worms to the focal plane of the microscope and to almost entirely 2-dimensional
behavior. Microfluidic devices were fabricated as described previously
\cite{Chung2011}. Methylcellulose solutions were prepared at
concentrations of $0.5\%$, $1\%$, $2\%$, and $3\%$ weight in volume of M9 buffer. To ensure that single animals
could be isolated in single chambers of the unbonded microchamber microfluidic device, we first picked
animals onto a room-temperature, unseeded plate. To ensure that animals were fully immersed in
methylcellulose mixture, we used a glass pipet to aspirate a small amount of methylcellulose solution,
and then aspirated animals from the unseeded plate one at a time into the methylcellulose solution.
Then, single animals surrounded by methylcellulose mixture were pipetted into individual chambers of
an unbonded PDMS chamber device. The device could then be flipped over onto a sterile 10cm Petri
dish and gently pressed down until the individual chamber walls came into contact with the Petri dish,
preventing animals from leaving their chambers. Animals were then imaged in devices for about 5
minutes at $30$ frames per second, resulting in time series data with $10665$ points.

To extract midline data from videos, we first found masks for each frame to isolate the worm from
the background using a combination of Otsu thresholding \cite{Otsu1979}, image smoothing using a Gaussian kernel, and size filtration. Otsu thresholding is a thresholding algorithm based on the gray-level histogram of an image. The threshold is identified by the grayscale pixel value that minimizes the intra-class variance of background and foreground pixels. We then broadly followed the method used in Stephens et al. \cite{Stephens2008} to represent the worm’s posture in ``worm-centric'' coordinates. Briefly, we found the midline of the worm in each frame by thinning the mask to a single line and interpolating between pixels of this line such that the midline was represented by 101 evenly-spaced points. We calculated the tangent angle between each pair of adjacent points along the midline so that the animal’s posture could be represented as a vector of angles, and then transformed those vectors with PCA so users could balance accuracy requirements and resource limitations via truncation of the data. We replaced frames in which animals were self-occluded with the data from the most recent non-self-occluded frame. We used untruncated PCA data for most computations because it has the same persistence output as the raw angle data. We used truncated PCA data (the first $5$ principal components) for the cycle representative computations in \Cref{sec:case-study}. 

The videos for this study were selected from a much larger set of data based on how well  they could be segmented and skeletonized. Some videos have subsequences that are difficult to automatically skeletonize because the animals self-occlude, \ie bend in such a way as to cross over themselves. Thus, this dataset is likely biased towards less complex behaviors like thrashing and there are some cases where there are multiple videos of the same animal. The resulting data set has $40$ samples of $10665$ points each with $10$ samples for each viscosity condition.

\subsection{Sliding window embeddings} 
\label{sec:swe}
Sliding window embeddings turn time series data into point cloud data in a way that does not forget the temporal information of the time series. There are some additional benefits to sliding window embeddings, including that they ``separate" points that intersect each other in a time series, such as in \Cref{eg:figure-eight,eg:metronome}.

\begin{definition} \label{def:time-series}
  A \defsty{time series} is a sequence of vectors $(x_t)_{t\in T}=(x_t,x_{t+1},x_{t+2},\ldots)$ where each $x_t$ is in the same finite-dimensional vector space $V$ and $T$ is a totally ordered set. 
\end{definition}

\begin{remark}
  The totally ordered set $T$, which indexes the time series, can be $\ZZ$, $\NN$, or a finite set like $[N]=\{1,2,\ldots,N\}$. For many applications including the ones in this paper, the indexing set is finite and will be omitted in notation for brevity, as in $(x_t)_t$.
\end{remark}

Given any time series we can construct a new time series called a sliding window embedding, which is also known as a time delay embedding with a lag or delay time of $1$. 

\begin{definition} \label{def:swe}
  Given a time series $\tau=(x_t)_{t}$ with vectors $x_t\in V$, a \defsty{sliding window embedding} of \defsty{window length} $\wl$ of $\tau$ is a new time series, $\swe{\tau}{\wl}=(\tilde{x_t})_{t}$, with 
  \[\tilde{x_t} = [x_t \quad x_{t+1} \quad\ldots\quad x_{t+\wl-1}] \in V^\wl  \]
  where $[ \cdot ]$ is concatenation of vectors. 
\end{definition}

That is, the $t^{th}$ vector in the new time series is the concatenation of $\wl$ consecutive vectors in the original time series and has dimension equal to $\wl \cdot \dim(V)$. 

\begin{remark} \label{rem:length-of-swe}
  If the original time series has $N$ points, then the sliding window embedding of window length $\wl$ has $N-\wl+1$ points, as one can see in \Cref{eg:size-of-swe}. 
\end{remark}
\begin{example} \label{eg:size-of-swe}
  Consider the time series $\tau=([1, 2], [3, 4], [5, 6], [7,8], [9,10])$ in $\RR^2$. The sliding window embedding of $\tau$ of window length $\wl=3$ is 
  \[\swe{\tau}{3}=([1,2,3,4,5,6], [3,4,5,6,7,8], [5,6,7,8,9,10]) \subseteq \RR^6\]
  which has $5-3+1=3$ points. 
\end{example}

We applied persistent homology (\Cref{sec:persistent-homology})
to sliding window embeddings of \emph{C. elegans} video data in order to quantify behavior 
Degree $1$ persistent homology detected cycles in these sliding window embeddings which we show correspond to particular behaviors. 

The cycles that persistent homology detects may consist of collections of points that trace out a closed curve. A cycle is ``large'' or highly persistent if it encloses an area that could fit a large ball; a cycle that is tall and skinny has small persistence.

Below we see two examples where change{a} time series exhibits a single periodic behavior but persistent homology will detect either two or zero non-trivial cycles. In contrast, the persistent homology of a sliding window embedding detects exactly one non-trivial cycle in both examples. 

\def\egwidth{0.29\textwidth}
\begin{example} \label{eg:figure-eight}
  \Cref{fig:figure-eight} (A) displays one period of a periodic time series in $\RR^2$ with the property that if successive points are connected by line segments then the path of the time series self-intersects. 
  To discover this figure-eight-shaped loop, one might try to use persistent homology (\Cref{sec:persistent-homology}).
  However, persistent homology would detect two distinct loops, each comprising half of the period. 
  See \cref{fig:fig8-simpl-cx-persistence} (D) for an illustration of these loops. 
  
  \Cref{fig:figure-eight} (B) and (C) show two-dimensional PCA projections of sliding window embeddings of the figure-eight for $\wl=10$ and $\wl=20$, respectively, using the first and third principal components. In these point clouds the time series draw out simple closed curves, and in fact in each of these cases persistence detects a single loop. 
  
  Notice that as the window length increases, the ``size'' of the loop increases. This increase in the loop's persistence makes it easier for persistent homology to robustly detect it.
   
  \begin{figure}[H] 
    \begin{center}
      \begin{tabular}{c c c}
        \multicolumn{1}{l}{(A)} & \multicolumn{1}{l}{(B)} & \multicolumn{1}{l}{(C)}\\
        \includegraphics[trim=0 150 0 125,clip, width=\egwidth]{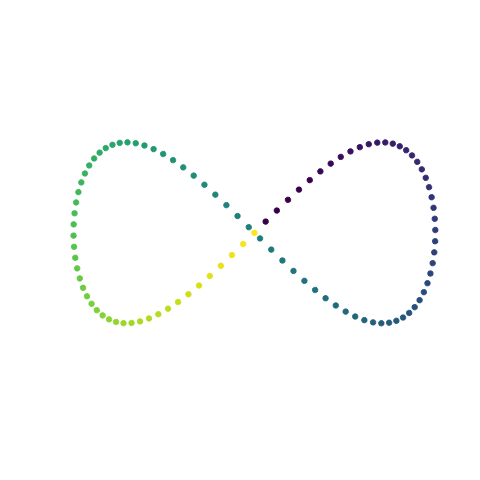}
        & \includegraphics[trim=0 150 0 125,clip, width=\egwidth]{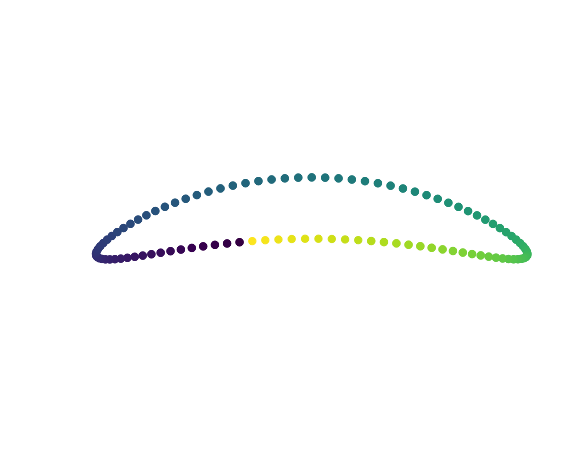}
        & \includegraphics[trim=0 150 0 125,clip, width=\egwidth]{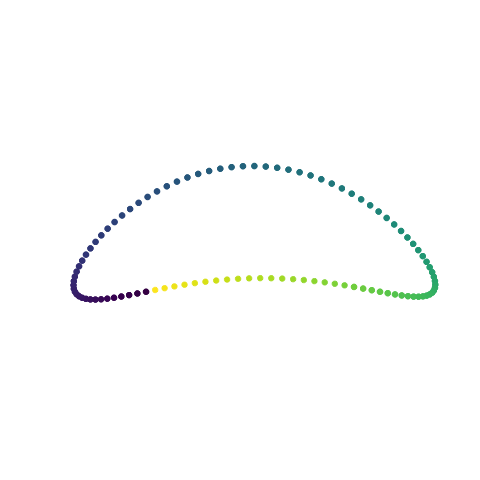} \\
        time series $\tau$
        & projection of $\swe{\tau}{10}$ 
        & projection of $\swe{\tau}{20}$ 
      \end{tabular}
    \end{center}
    \caption[Figure eight example of sliding window embeddings]{(A) A time series in $\RR^2$ determines a self-intersecting curve. (B) The sliding window embedding of window length $10$ separates the previously intersecting segments of the curve. (C) A sliding window embedding with a higher window length separates the intersecting segments even further. With too small of a window length $\wl$, the resulting loop will be relatively flat and long and will therefore have small persistence and be difficult to differentiate from noise. }
    \label{fig:figure-eight}
  \end{figure}
\end{example}

\begin{example} \label{eg:metronome}
  \Cref{fig:metronome} (A) shows a $1$-dimensional time series that is a discretization of a sine wave. This periodic behavior creates no loops --- in fact, because the points take values in $\RR$, the time series cannot produce degree $1$ homology. However, a sliding window embedding, in this case of window length $4$, creates a loop that is detected by persistent homology. That loop in $\RR^4$ is projected down to two  dimensions in \Cref{fig:metronome} (B). 
  \begin{figure}[H]
    \begin{center}
      \begin{tabular}{c c}
        \multicolumn{1}{l}{(A)} & \multicolumn{1}{l}{(B)}\\
        \includegraphics[trim=0 150 0 125,clip, width=\egwidth]{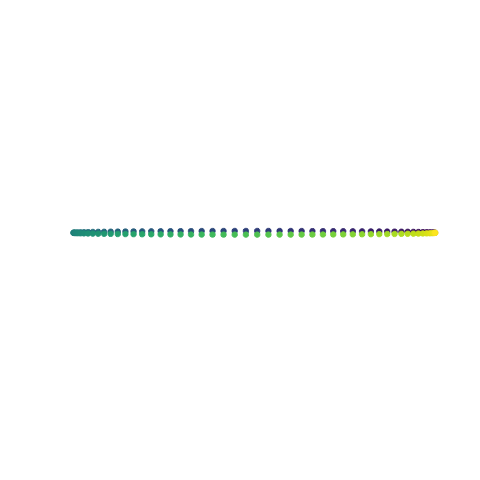}
        & \includegraphics[trim=0 150 0 125,clip, width=\egwidth]{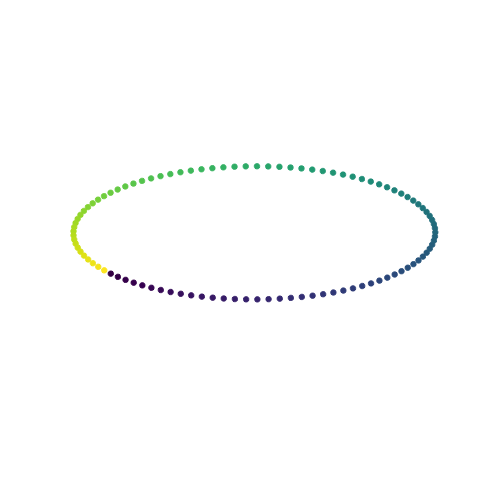} \\
        time series $\tau$
        & projection of $\tau^{4}$
      \end{tabular} 
    \end{center}
    \caption[Metronome examples of sliding window embeddings]{Sliding window embeddings encode periodic behavior in the form of loops. (A) A periodic time series that has trivial first homology. (B) A sliding window embedding with $\wl=4$ of the original time series that has nontrivial first homology. }
    \label{fig:metronome}
  \end{figure}
\end{example}

\newpage
\subsection{Persistent homology} 
\label{sec:persistent-homology}

In this section we provide an overview of persistent homology and how it may be used to produce quantitative summaries of the shape of a collection of points such as the sliding window embedding discussed above.

\begin{definition} \label{def:simpl-cx}
  A \defsty{simplicial complex} on a set of vertices $V$ is a collection $K$ of non-empty subsets of $V$ such that if $\tau \in K$ and $\tau' \subset \tau$, then $\tau' \in K$. An element $\tau\in K$ is called a \defsty{simplex}. An \defsty{$n$-simplex} or simplex of dimension $n$ is a simplex $\tau\in K$ with $|\tau| = n+1$.  
  The \defsty{1-skeleton} of a simplicial complex $K$ is the set of simplices with dimension at most one.
  A \emph{filtered simplicial complex} or \emph{filtration} is a collection $\{K_r\}_{r \in S}$ of simplicial complexes $K_r$ where $S \subseteq \mathbb{R}$ such that $K_r \subseteq K_s$ for all $r,s \in S$ with $r \leq s$.
\end{definition}

\begin{definition}\label{def:VR-cx-scale-r}
  Let $X\subset \mathbb{R}^d$ be a finite set and let $r\geq 0$. The \defsty{Vietoris-Rips complex of $X$ at scale $r$}, denoted $\mathcal{R}_r(X)$, is the simplicial complex with vertex set $X$ and whose simplices are given as follows. A subset $\{x_0,\dots,x_n\}\subset X$ is an $n$-simplex in $\mathcal{R}_r(X)$ if and only if $|x_i-x_j|\leq r$ for all $i,j\in \{0,\dots,n\}$.
\end{definition}

\begin{definition}\label{def:VR-cx}
  The \defsty{Vietoris-Rips filtration} of a finite set $X\subset \mathbb{R}^d$ is the collection $\mathcal{R}(X) := \{\mathcal{R}_r(X)\}_{r\geq 0}$.
\end{definition}

Note that while the Vietoris-Rips complex of $X$ is parameterized by the non-negative reals, the finiteness of $X$ guarantees that $\mathcal{R}(X)$ consists of only finitely many distinct simplicial complexes. 

\begin{figure}[H]
  \begin{center}
    \begin{tabular}{cccc}
      \multicolumn{1}{l}{(A)}\\
      radius=$0.5$ & radius=$0.9$ & radius=$5.0$ & radius=$8.9$ \\
      \includegraphics[trim = 50 150 50 70, clip, width=\quarterwidth] {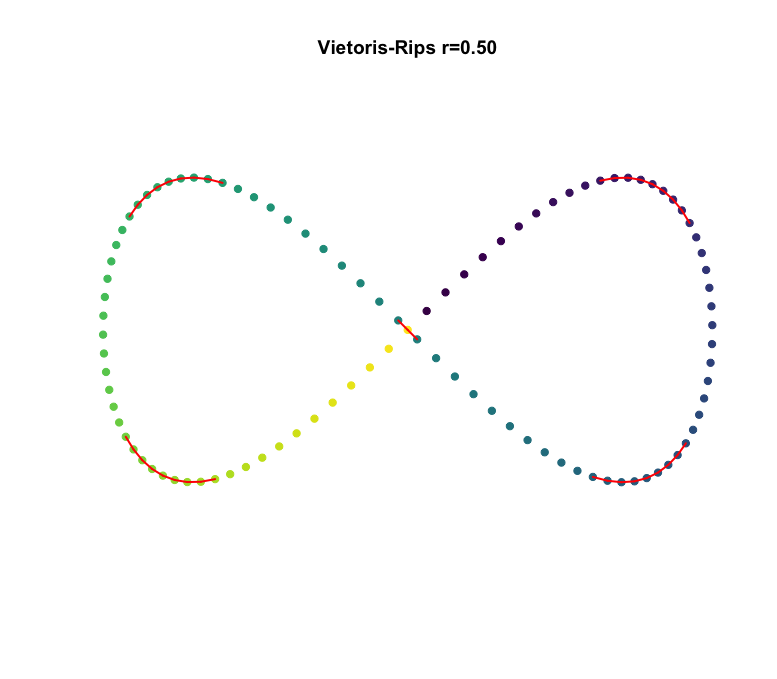}
      &\includegraphics[trim = 50 150 50 70, clip, width=\quarterwidth] {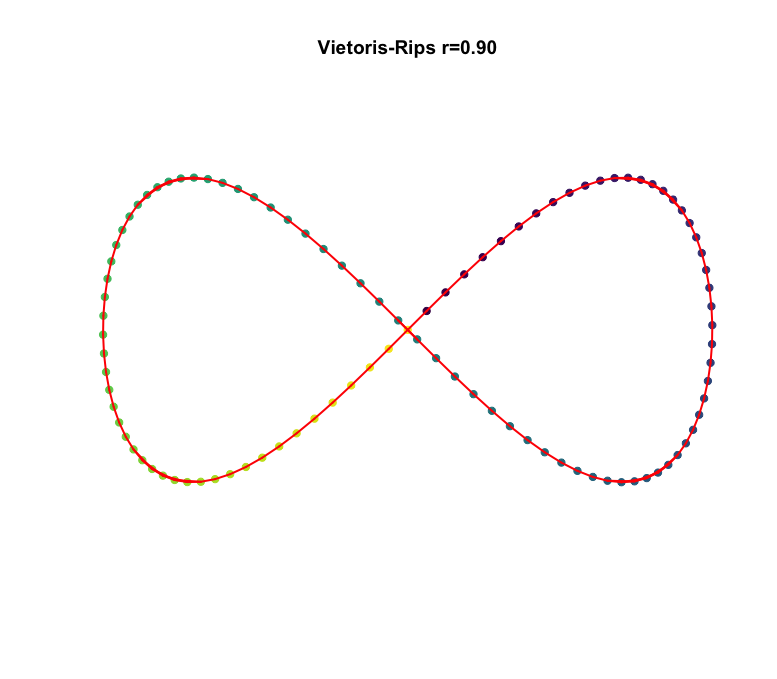}
      &\includegraphics[trim = 50 150 50 70, clip, width=\quarterwidth] {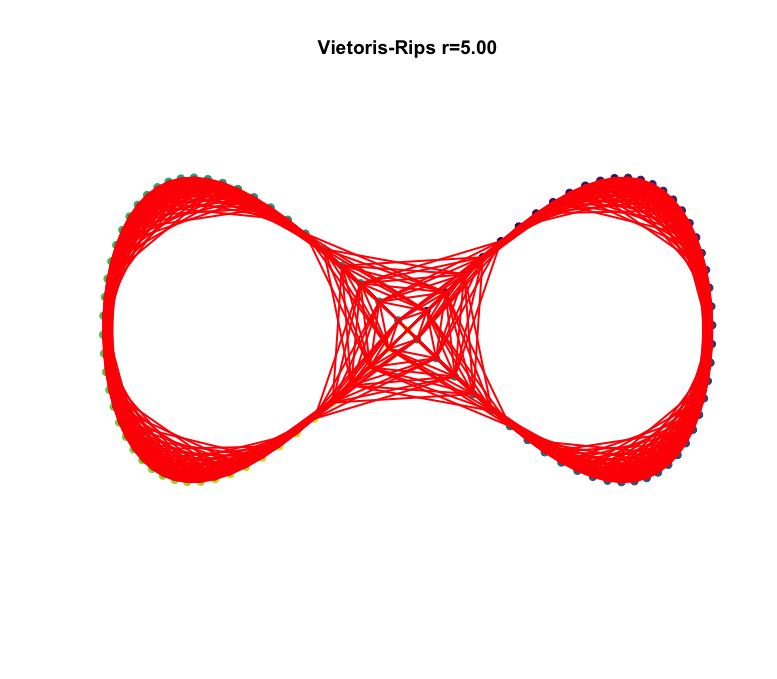}
      &\includegraphics[trim = 50 150 50 70, clip, width=\quarterwidth] {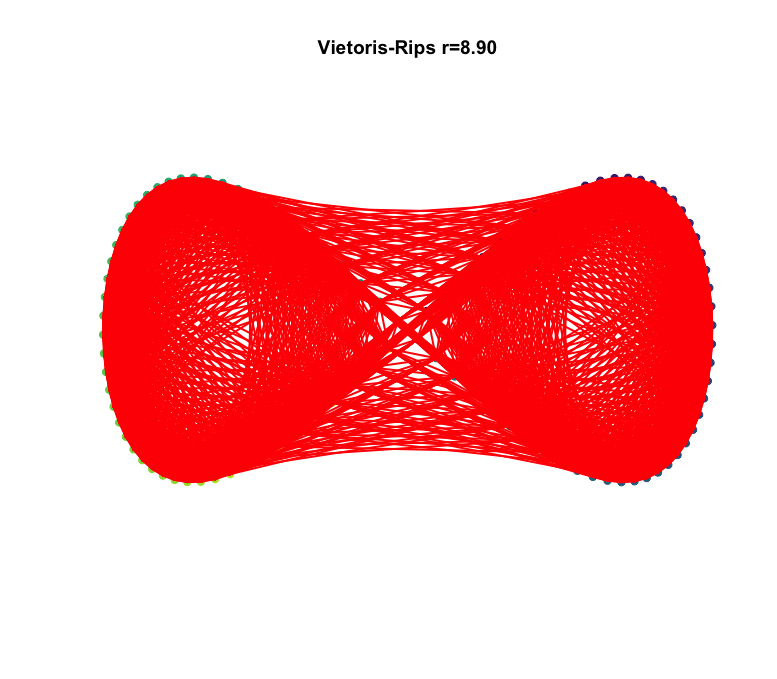}
      \\
      $0$ loops & $2$ loops & $2$ loops & $0$ loops 
    \end{tabular}
    \\[2em]
    \begin{tabular}{cc}
      \multicolumn{1}{l}{(B)} & \multicolumn{1}{l}{(C)}\\[1.5em]
      \includegraphics[trim = 0 0 0 50, clip, height=0.3\columnwidth]{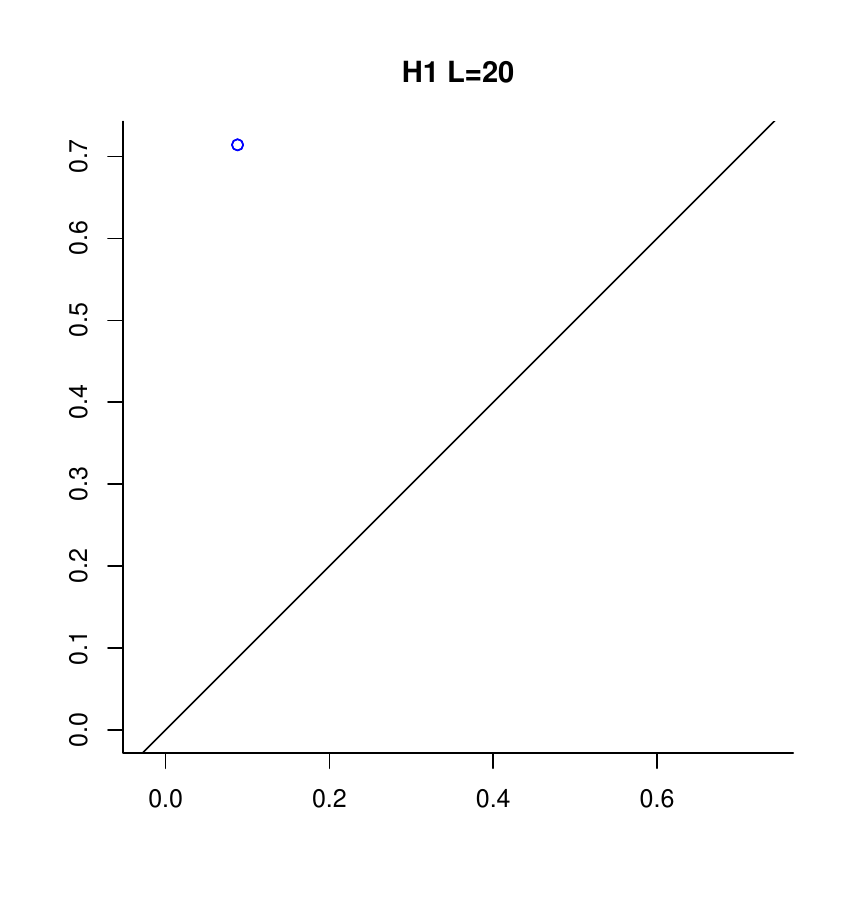}
      &\includegraphics[trim = 0 0 0 50, clip, height=0.3\columnwidth]{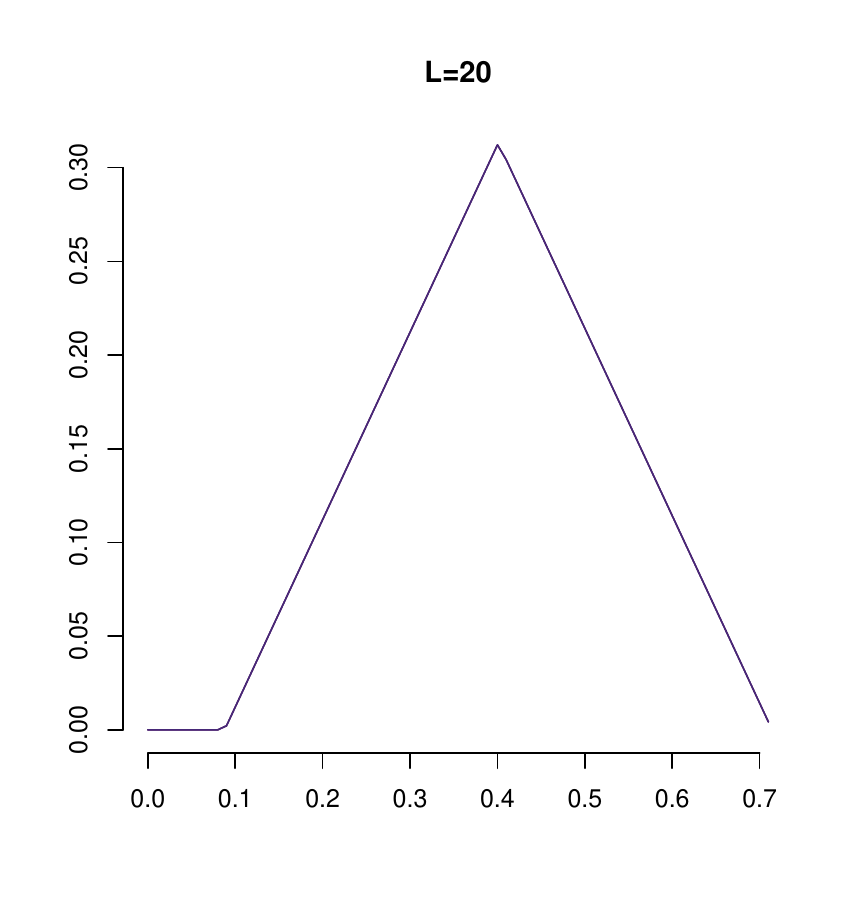}
      \\[1.5em]
      \multicolumn{1}{l}{(D)} \\
      \trimmedwithaxesgraphic[width=\halfwidth]{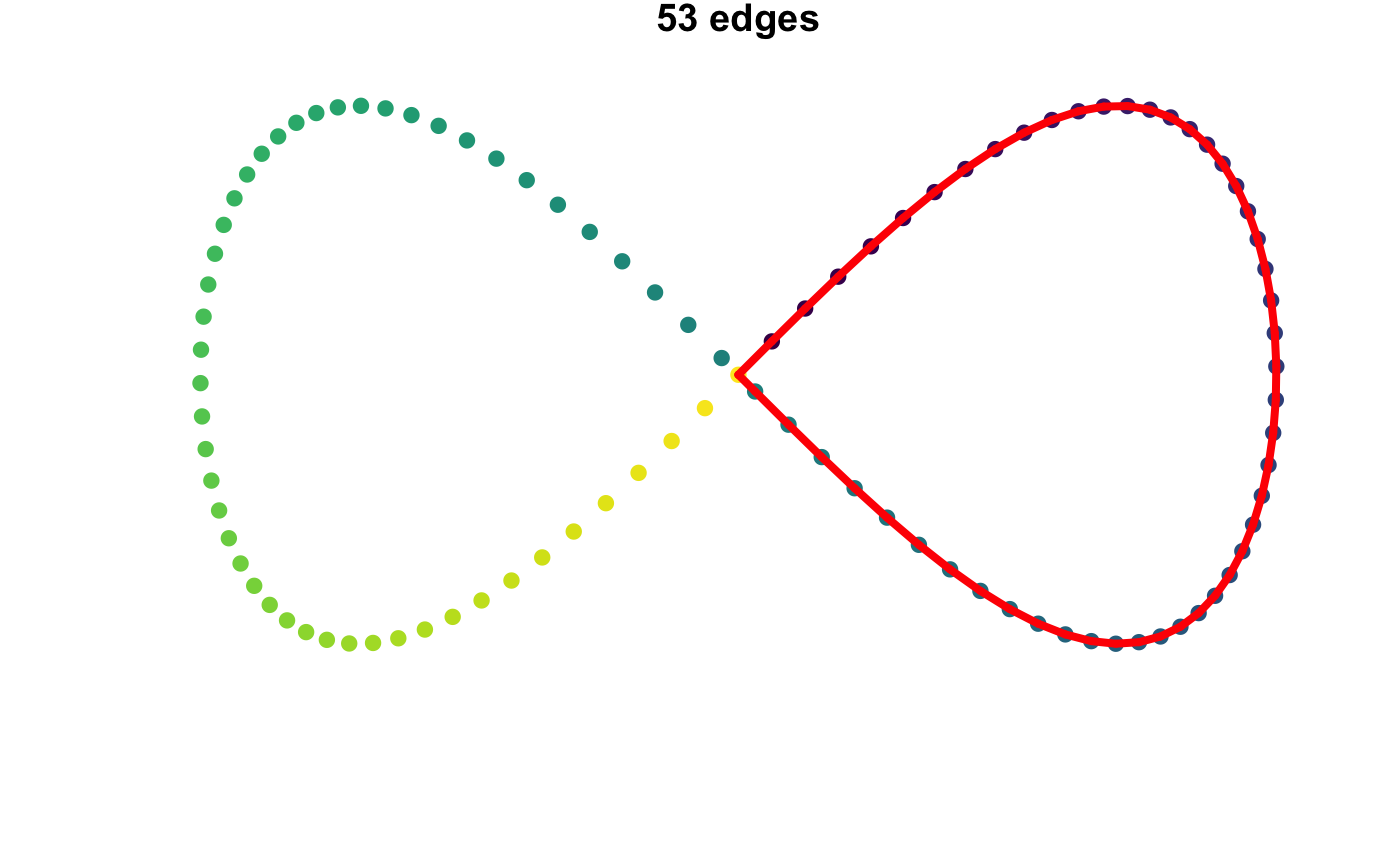}
      & \trimmedwithaxesgraphic[width=\halfwidth]{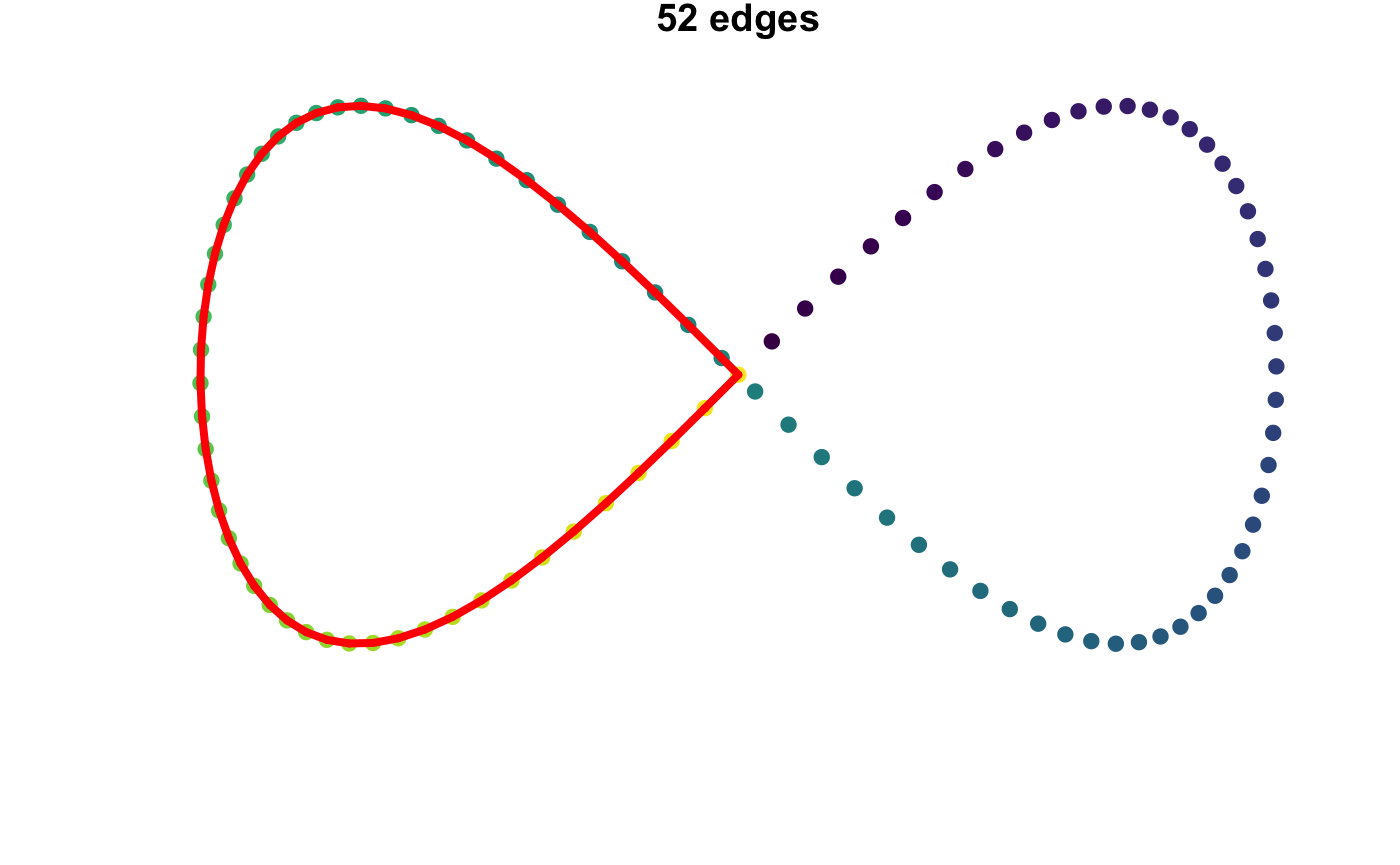}
    \end{tabular}
  \end{center}
  \caption[Vietoris-Rips and persistence of figure-eight]{
    (A) The $1$-skeleton of the Vietoris-Rips filtered simplicial complex of a figure-eight-shaped point cloud at four scales. 
    (B) The degree $1$ persistence diagram of the figure eight in $2$ dimensions. Note that the point has multiplicity $2$. 
    (C) The corresponding degree $1$ persistence landscape. The first and second landscapes are nonzero and identical and all other landscapes are trivial. 
    (D) The two loops that generate the homology of the Vietoris-Rips complex on the figure eight. }
  \label{fig:fig8-simpl-cx-persistence} 
\end{figure}

\begin{example} \label{eg:fig8-filt-simpl-cx}
  \Cref{fig:fig8-simpl-cx-persistence} (A) shows the 1-skeleton of the Vietoris-Rips complex  of a pointcloud in $\RR^2$ at four scales. Notice that each simplicial complex includes into the next. 
\end{example}

The persistent homology of a Vietoris-Rips filtration can be represented by a multiset in $\mathbb{R}^2$ called a \emph{persistence diagram} in which each point gives the  scale of the appearance and disappearance of a topological feature (such as a loop) in the filtration.

\begin{example} \label{eg:fig8-persistence}
  \Cref{fig:fig8-simpl-cx-persistence} (B, C) show the persistent homology in degree 1 of \Cref{eg:fig8-filt-simpl-cx}. Notice that both (B) the persistence diagram and (C) the persistence landscape (see \Cref{def:landscape}) show two cycles, but because they are born and die at exactly the same radius parameters they are plotted in the same place. The two cycles are shown in (D). 
\end{example}

It is difficult to apply standard tools of statistics and machine learning directly to persistence diagrams, which, for example, need not have unique averages \cite{Mileyko2011}. 
A solution is to map persistence diagrams into a vector space or Hilbert space. One such mapping is the persistence landscape. See \cite{persistence-landscapes} for the following definitions and results.
\begin{definition} \label{def:landscape}
  For $a<b$ let $f_{a,b}: \mathbb{R} \to \mathbb{R}$ be the piecewise-linear function given by 
  \[ f_{a,b}(t)=\begin{cases}
    t-a, \text{ if } a \leq t \leq \frac{a+b}{2} \\
    b-t, \text{ if } \frac{a+b}{2} \leq t \leq b \\
    0, \text{ otherwise.}
  \end{cases} \] 
  Given a persistence diagram $\textup{Dgm}_p(\mathcal{K})$, the corresponding \defsty{$k$th persistence landscape} is the function $\lambda_k: \mathbb{R} \to \mathbb{R}$ given by defining $\lambda_k(t)$ to be the $k$th largest value of $f_{a,b}(t)$ over all points $(a,b)\in \textup{Dgm}_p(\mathcal{K})$. The \defsty{persistence landscape} is the sequence $(\lambda_k)_k$. The parameter $k$ is called the \defsty{depth} of the persistence landscape. For a point cloud $X$, we will denote by $\textup{PL}(X)$ the persistence landscape obtained by applying degree $1$ persistent homology to the Vietoris-Rips filtration of $X$.
\end{definition}
Persistence landscapes have unique averages, satisfy the law of large numbers and central limit theorems, and can be discretized for computations. Because the sequence of functions that make up a landscape are nested, they can all be graphed on the same plot as in the right column of	\Cref{fig:case-study-persistence}. 

While the persistence landscape is defined to be an object in a space of continuous functions, it can be discretized and turned into a finite-dimensional vector. Through discretization, each depth of the landscape transforms from a continuous function on $\RR$ to a vector where the $i^{th}$ entry in the vector corresponds to the function value at the $i^{th}$ discrete parameter value.  The vectors for each depth of the landscape are concatenated together to produce a single high-dimensional vector.
These discrete landscapes can be computed directly, which we did for the computations outlined in \Cref{sec:pipeline}.

This vector space (in fact, Hilbert space) setting lets us use linear algebra-based statistical and machine learning techniques such as principal component analysis (PCA). The principal components from PCA on discretized landscapes can be converted into a format much like a persistence landscape --- a sequence of continuous functions on $\RR$ --- but the principal components are not themselves persistence landscapes because the functions fail to be nonnegative. 

\newpage
\subsection{Pipeline} \label{sec:pipeline}

In this section we give details for our analysis of \emph{C. elegans} data. 
The input consist of piecewise linear midlines of \emph{C. elegans} from video recordings as described in \Cref{sec:description-of-data}. These midlines were parameterized by the $100$ angles between adjacent segments and then were transformed using PCA, so each sample input to our system was a time series $\tau$ of $100$-dimensional vectors measured in radians. See Figure~1 in \cite{Stephens2008} and the accompanied description for more details on this parameterization of the \emph{C. elegans} midlines or see our short summary of the procedure in \Cref{sec:description-of-data}.

The time domain of this time series was divided into overlapping patches of a given size called the \defsty{patch length}, resulting in a collection $\{\tau_i\}_i$ of smaller time series. For our experiments, a patch length of $300$ was chosen, with adjacent patches overlapping by half of the patch length. The sliding window embeddings of the $\tau_i$ were then computed with window length parameter $\wl=20$, resulting in a new collection $\{(\tau_i)^\wl\}_i$ of time series of length $300-\wl+1=281$.
This analysis is not particularly sensitive to the choices of the hyperparameters patch length and window length; only extreme changes in either parameter lead to significant changes in results. The hyperparameter choices were motivated by the timescales at which \emph{C. elegans} complete meaningful behaviors: for patch length, $150$ frames of $30$ fps video is $5$ seconds of behavior; for window lengths, $20$ frames is $0.67$ seconds and corresponds to roughly one period of forward crawling in adult \emph{C. elegans} submerged in the $0.5\%$ methylcellulose environment. Strategies for choosing appropriate window lengths are described in \cite{Perea2015b}. The method for cross validation of the choice of window length is described in \Cref{sec:null-model}.  

Persistence diagrams were computed for each of the patches $(\tau_i)^\wl$. This step accounts for the vast majority of the computational resources of the pipeline, 
and the computational costs are made worse by the concatenation of vectors in a sliding window embedding. This is where we greatly benefit from the preprocessing that turns video data, which is extremely high-dimensional (see \cite{Tralie2018} Section 3.1), into a $100$-dimensional time series. 
On a 2017 $15$-inch MacBook Pro with a $2.8$ GHz Intel Core i7 processor and $16$ GB of RAM, this step took $22588.632$ seconds, or about $6$ hours and $15$ minutes. 

\begin{figure}[t!]
  \begin{center}
    \begin{tikzpicture}
      \tikzstyle{box} = [rectangle, rounded corners, minimum width=1.5cm, minimum height=1cm, align=center, text centered, draw=black, font=\small]
      \tikzstyle{arrow} = [thick,->,>=stealth]
      
      \node (p1) [box] {Patch 1}; 
      \node (sw1) [box, right of=p1, xshift=1.2cm] {Sliding \\ window \\ embedding}; 
      \node (pd1) [box, right of=sw1, xshift=1.5cm] {Persistence \\ diagram}; 
      \node (pl1) [box, right of=pd1, xshift=1.5cm] {Persistence \\ landscape}; 
      
      \node (inp) [box, left of=p1, yshift=-2.1cm, xshift=-0.5cm] {Input}; 
      
      \node (apl) [box, right of=pl1, yshift=-2cm, xshift=1.5cm] {Average \\  persistence \\ landscape};
      \node (pn) [box, below of=p1, yshift=-3cm] {Patch n}; 
      \node (swn) [box, below of=sw1, yshift=-3cm] {Sliding \\ window \\ embedding}; 
      \node (pdn) [box, below of=pd1, yshift=-3cm] {Persistence \\ diagram}; 
      \node (pln) [box, below of=pl1, yshift=-3cm] {Persistence \\ landscape};
      \node  at ($(p1)!.5!(pn)$) {\Large\vdots};
      \node at ($(sw1)!.5!(swn)$) {\Large\vdots};
      \node at ($(pd1)!.5!(pdn)$) {\Large\vdots};
      \node at ($(pl1)!.5!(pln)$) {\Large\vdots};
      \draw [arrow] (inp) -- (p1);
      \draw [arrow] (inp) -- (pn);
      \draw [arrow] (p1) -- (sw1);
      \draw [arrow] (sw1) -- (pd1);
      \draw [arrow] (pd1) -- (pl1);
      \draw [arrow] (pn) -- (swn);
      \draw [arrow] (swn) -- (pdn);
      \draw [arrow] (pdn) -- (pln);
      \draw [arrow] (pl1) -- (apl);
      \draw [arrow] (pln) -- (apl);
    \end{tikzpicture}
  \end{center}
  \caption[Pipeline for a single worm] {Pipeline for a single worm.}
  \label{fig:pipeline}
\end{figure}
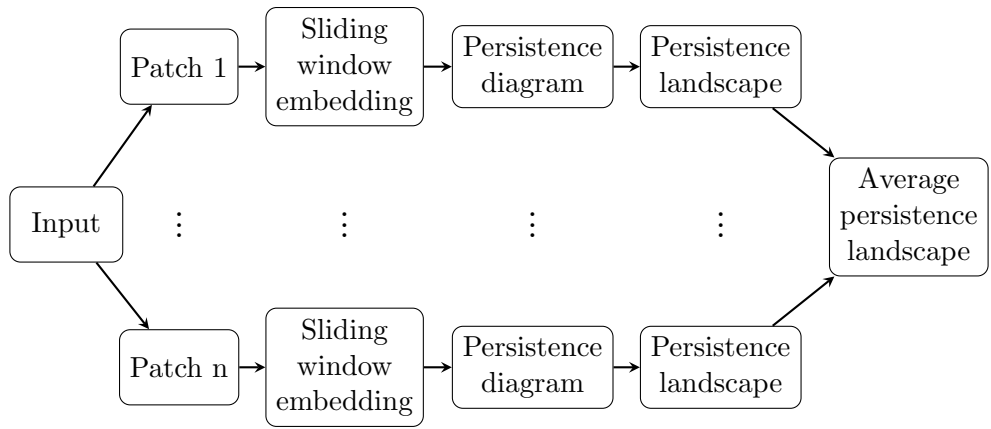

For each $(\tau_i)^\wl$, a (discretized) persistence landscape $\textup{PL}((\tau_i)^\wl)$ was computed from the persistence diagram of the Vietoris-Rips complex $\mathcal{R}((\tau_i)^\wl)$.
The grid of parameter values on which the persistence landscape was evaluated to produce its discretization was chosen to include all of the bars of the barcode and to be sufficiently fine to produce nice visualizations. Since the persistence landscapes are piecewise linear with slope bounded by $\pm 1$, the step size of this discretization bounds the error and there is eventually little to be gained from a finer discretization. The step size of the discretization we used was $0.1$.  
The maximum depth was chosen to include all nonzero depths of the persistence landscapes.
The collection $\{\textup{PL}((\tau_i)^\wl)\}_i$ of persistence landscapes was then assembled into a single summary for a given video 
by averaging the persistence landscapes across patches to result in a single average persistence landscape associated to each video. 
These steps are summarized in \Cref{fig:pipeline}.
For each environment viscosity, the average persistence landscapes for each video were averaged to give an average persistence landscape of the class.

The persistence landscapes for each class were used for analysis of the sliding window embedding as follows. Distances between each class' average persistence landscapes were computed using the usual Euclidean distance. The pairwise distances were visualized using multidimensional scaling to give a $2$-dimensional visualization of the similarities between classes; see \Cref{fig:multid-scaling}(A).

Principal component analysis (PCA) was applied to the set of average persistence landscapes for each video and the first two principal components were plotted as sequences of functions.
We plotted the PCA projection of these video average persistence landscapes together with the average of each class; see \Cref{fig:full-pca}.
These plots visualize some of the similarity between classes and variation within classes.

Next, we further studied the variation within classes in two ways.
First, the standard deviations of each coordinate were computed for the average persistence landscapes of the videos in each class.
These were visualized as sequences of functions in \Cref{fig:class-avg-PL-variation}(B) to show
the variation in different parts of the average persistence landscapes.
Second, we applied PCA to the average persistence landscapes of the videos in each class, along with the cumulative variances explained by the first $n$ principal components for $n=1\ldots,10$ ($10$ is the number of samples in each class). The first three principal components were also computed. See \Cref{fig:pca-by-class}. 

We conducted a permutation test on the pairwise Euclidean distances between the average persistence landscapes of each class.
We used $10,000$ permutations for each permutation test.
The approximate $p$-value equals the percentage of cases in which the distance is at least as large as the observed distance. Results can be found in \Cref{table:permutationtest-SVMconfusion} (A). 

We applied multiclass support vector machines (SVM) to classify samples according to viscosity of their environments. We use the ksvm function from the kernlab package in R on the average persistence landscapes of the videos. Accuracy was estimated using 10-fold cross validation with cost set to 10. Cross validation was repeated 20 times and the results were averaged. A confusion matrix for one instance of SVM with 10-fold cross validation was also computed and is shown in \Cref{table:permutationtest-SVMconfusion} (B). 

We used support vector regression (SVR) to approximate viscosity of the worm environments given the worm's behavior data. The goal was to assess how predictive our techniques are. Accuracy was estimated by averaging 10 repetitions of  10-fold cross validation. Results are plotted in \Cref{fig:SVR}. 

As a final step we applied cross validation to the hyperparameter window length. We cross-validated the choice of window length by running the above pipeline for window lengths of $1$, $10$, $20$, and $30$ on a subset of the data. To reduce total computation time,  we restricted to the first minute of each video, which corresponds to $1800$ frames. We then compared the permutation test and multiclass SVM results to see which hyperparameter choice gave the best  results. The results from cross validation of window length can be found in \Cref{sec:wl-cross-validation}. 

\subsection{Validation} \label{sec:null-model}

To help validate our pipeline, we compared our results to those obtained by applying the same computational procedure to a null model given by randomly permuting the frames in each video.

We also studied the effects of using a preprocessing step different from sliding window embeddings: moving average filters. A moving average filter of window length $\wl$ of time series data creates a new time series. This time series has the same length as the corresponding sliding window embedding and each point in the time series is constructed using the same set of points from the original time series. The moving average filter, however, takes the average of the vectors in its window, instead of concatenating them. 

Furthermore, we compared our results to the ones obtained from two simpler techniques. 
For the first technique we attempted to characterize \emph{C. elegans} behavior using the speed of the worm. To do this, we computed $2$-norms of the differences between two consecutive frames of angle data and then averaged all of those discrete derivatives, which resulted in a single value per sample. 
For the second technique, which could be described as measuring the variance of the worm's pose over a video, we computed standard deviations for the angle data coordinate-wise, which resulted in a vector of length $100$ per sample. In each case we performed a permutation test and multiclass SVM on the resulting feature vectors. The results of these experiment appears in \Cref{sec:null-model-results}.

\section{Results}\label{sec:results}

We present the results of a case study of a single sample of behavior data and an experiment on the effects of viscosity of the surrounding environment on \emph{C. elegans} locomotion and behavior. The case study assesses a significantly smaller data set and directly links topological features to specific behaviors. The experimental results in \Cref{sec:experimental-results} are more difficult to directly interpret in terms of specific behaviors but nonetheless we show the effectiveness of persistent homology in distinguishing variations in behaviors here. We leave more explicit interpretation of persistence output in terms of specific behaviors for future work.

\subsection{An illustrative case study} \label{sec:case-study}

The following results were obtained by carefully analyzing a sample of \emph{C. elegans} behavior data from a video of a worm crawling on agar. The sample consists of $400$ frames of a $30$ frames per second video, so roughly $13$ seconds of behavior. Having a solid surface to provide friction forces slower but more complicated behavior than we would see in an aqueous environment, and we take advantage of the resulting clarity of the data. 

The data we analyzed consists of a $20$ second video where the subject exhibits the following behaviors in chronological order:
\begin{enumerate}[nolistsep]
  \item crawl forward,
  \item crawl backward,
  \item pause, and
  \item crawl backward again. 
\end{enumerate}

Below we analyze the time series $\tau$ from this data, the corresponding sliding window embedding $\tau^{20}$, the corresponding moving average filter, and the sliding window embedding of the corresponding null model. These comparisons illustrate various strengths of the sliding window embedding: it smooths the noise from the original time series; it retains more geometric data than the moving average filter; and it captures temporal data from the original time series that is destroyed in the null model. Then, using representative cycles we construct synthetic \emph{C. elegans} midline data that produce a stereotyped, forward crawl and explain how this process gives concrete interpretations of persistence features in terms of synthetic behavior data. See \Cref{sec:representative-cycles}.  

\begin{figure}[H] 
  \begin{center}
    \begin{tabular}{c c c c}
      PC1 $\times$ PC2 & PC2 $\times$ PC3 & diagram & landscape
      \\
      \multicolumn{1}{l}{(A) $\tau$}\\
      \includegraphics[width=\quarterwidth,trim=50 80 10 40,clip]{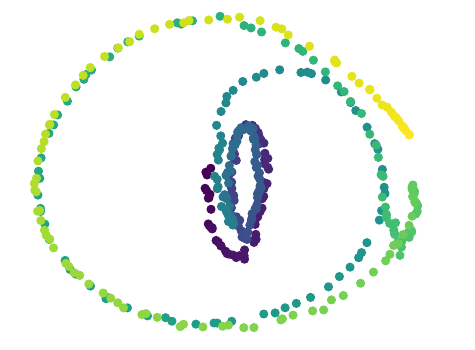} 
      &\includegraphics[width=\quarterwidth,trim=55 120 25 70,clip]{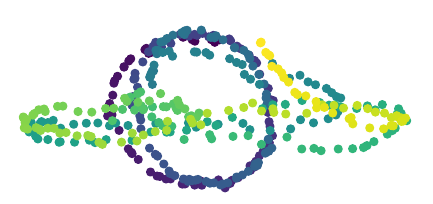}
      &\includegraphics[width=\quarterwidth,trim=0 35 0 20,clip]{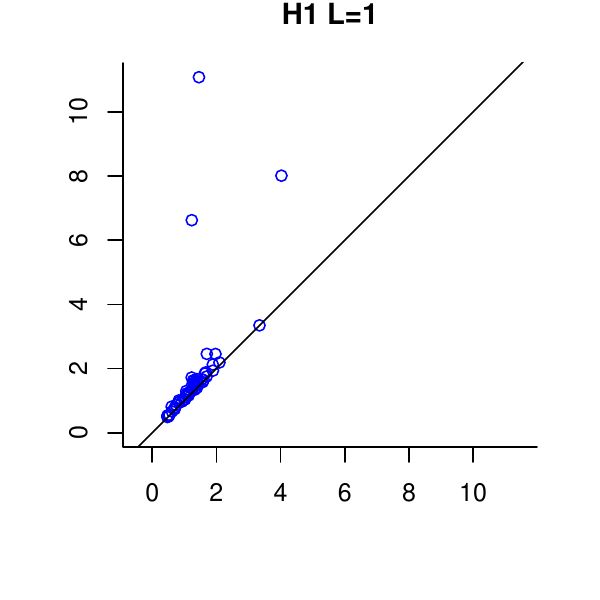}
      &\includegraphics[width=\quarterwidth,trim=0 35 0 20,clip]{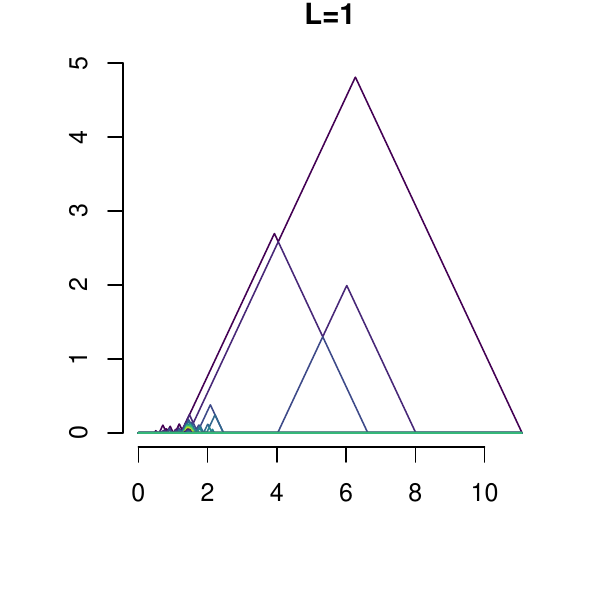}
      \\
      \multicolumn{1}{l}{(B) $\tau^{20}$}\\
      \includegraphics[width=\quarterwidth,trim=50 80 10 40,clip]{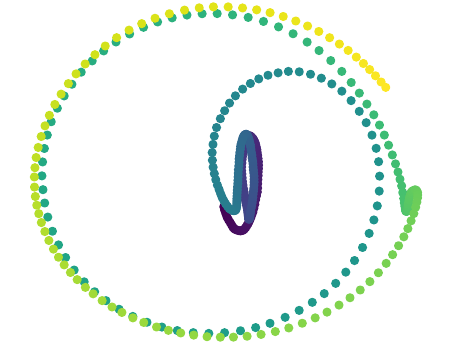} 
      &\includegraphics[width=\quarterwidth,trim=55 120 25 70,clip]{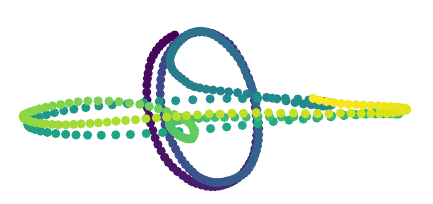}
      &\includegraphics[width=\quarterwidth,trim=0 35 0 20,clip]{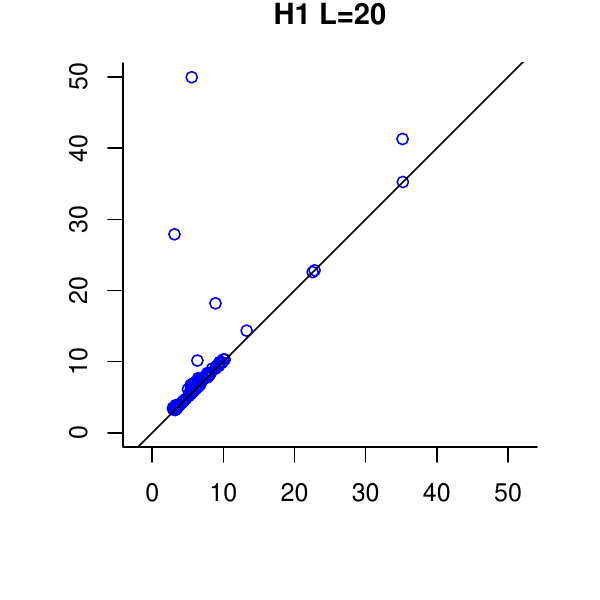}
      &\includegraphics[width=\quarterwidth,trim=0 35 0 20,clip]{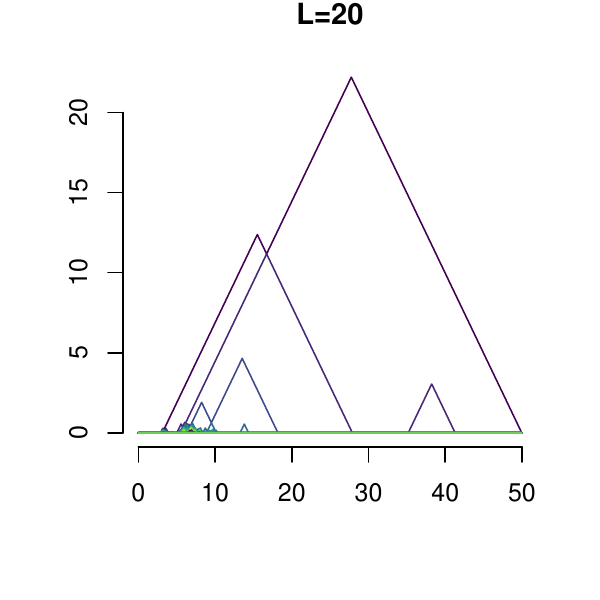}
      \\
      \multicolumn{1}{l}{(C) smoothed $\tau$}\\ 
      \includegraphics[width=\quarterwidth,trim=50 80 10 40,clip]{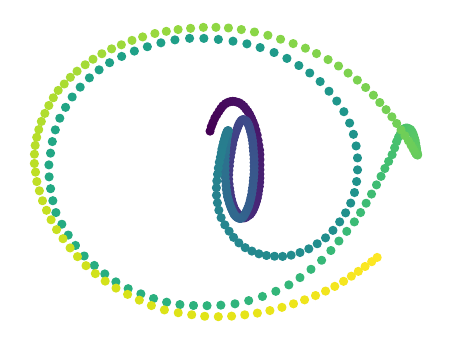} 
      &\scalebox{-1}[1]{\includegraphics[width=\quarterwidth,trim=50 90 10 50,clip]{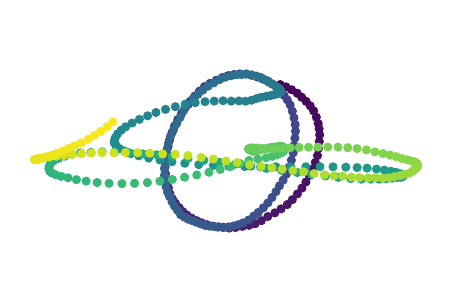}}
      &\includegraphics[width=\quarterwidth,trim=0 35 0 20,clip]{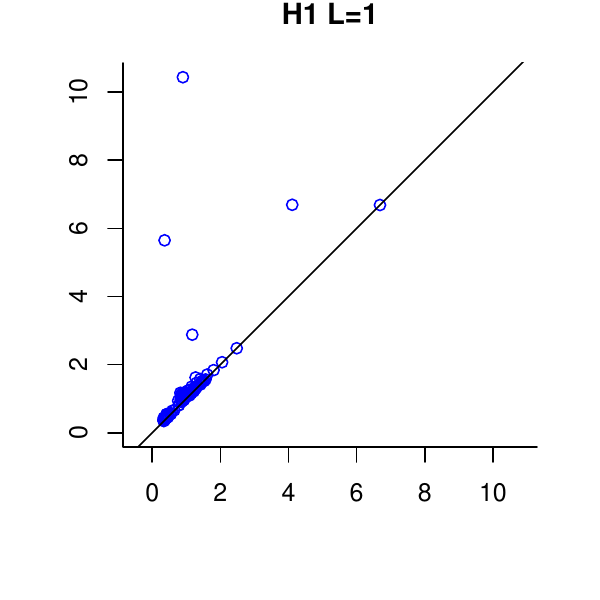}
      &\includegraphics[width=\quarterwidth,trim=0 35 0 20,clip]{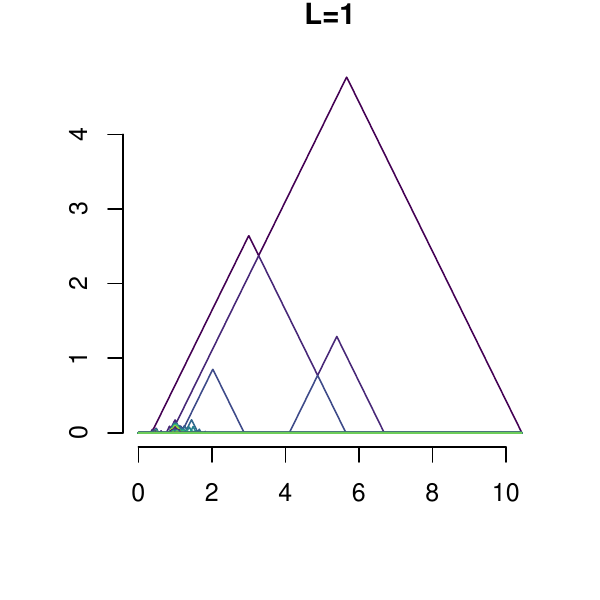}
      \\
      \multicolumn{1}{l}{(D) (scrambled $\tau$)$^{20}$}\\
      \includegraphics[width=\quarterwidth,trim=30 40 10 20,clip]{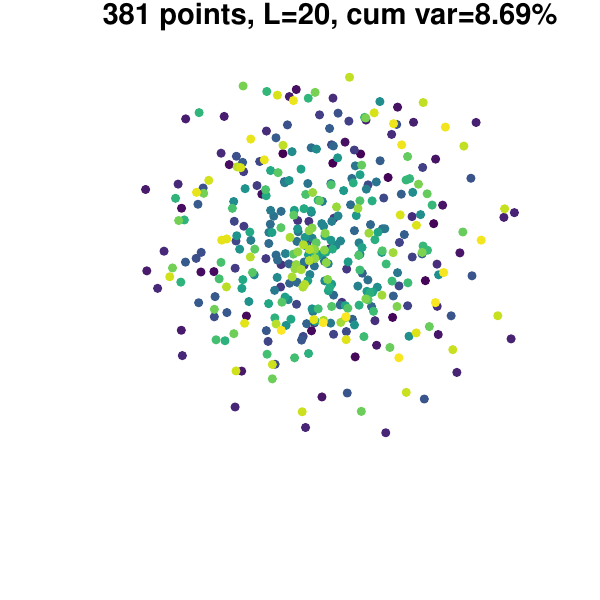} 
      &\includegraphics[width=\quarterwidth,trim=50 60 10 20,clip]{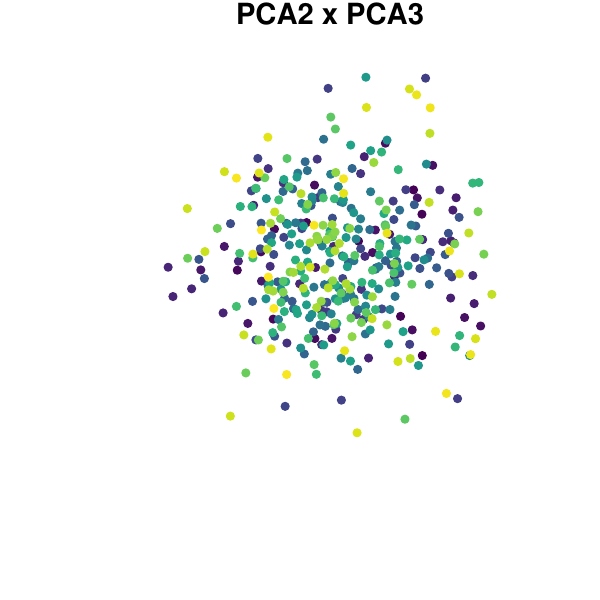}
      &\includegraphics[width=\quarterwidth,trim=0 35 0 20,clip]{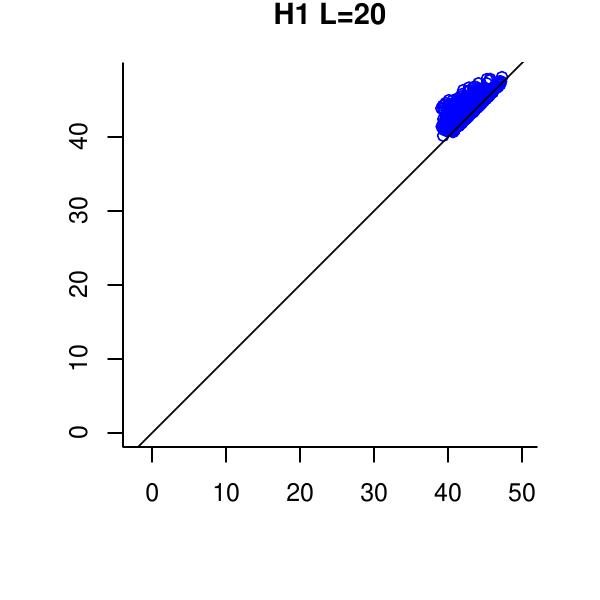}
      &\includegraphics[width=\quarterwidth,trim=0 35 0 20,clip]{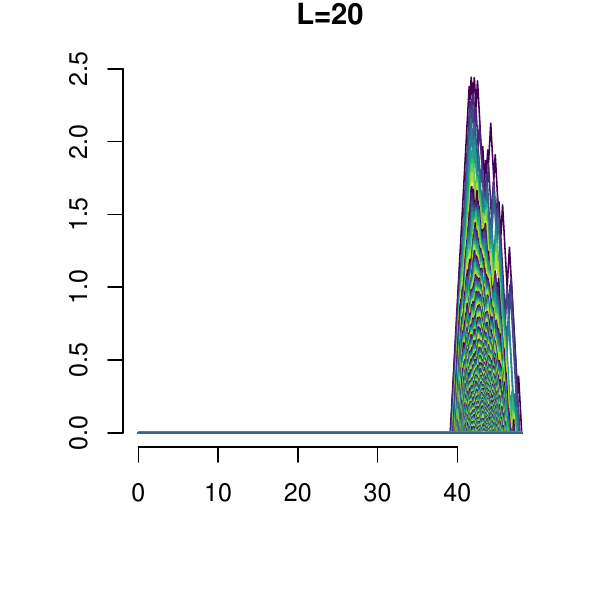}
    \end{tabular}
  \end{center}
  \caption[time series/swe with diagram and landscape]{(A) A time series $\tau$, (B) its sliding window embedding of window length $\wl=20$, (C) its smoothing by moving average filter of window length $20$, and (D) the null model,
    all projected onto $2$-dimensional axes of principal components. The corresponding persistence diagrams and persistence landscapes are to the right. }
  \label{fig:case-study-persistence}
\end{figure}

To visualize the four point clouds on which we will compute persistence, we apply PCA and project onto the first few principal components. Some of these projections are shown in the two left-most columns of \Cref{fig:case-study-persistence}. In contrast to the three other time series, the null model time series in (D) has no discernible geometric structure. It appears that the corruption of temporal information has destroyed the interpretability of the visualizations of sliding window embeddings. Meanwhile, the original time series $\tau$ in (A) has a similar shape to its sliding window embedding $\tau^{20}$ in (B), with the caveat that the sliding window embedding has the effect of smoothing the data and making it more robust to noise. The moving average filter in (C) also has this smoothed property. Though the three time series in \Cref{fig:case-study-persistence} (A-C) have similar shapes, persistence diagrams, and persistence landscapes, they vary in one important feature: the pause. 

In \Cref{fig:case-study-behaviors} (A) the points in the sliding window embedding that correspond to frames where the worm is performing a specific behaviors are highlighted. The points corresponding to the pause behavior deviate from the path of points corresponding to crawling backwards. This deviation is small compared to the noisiness of the original time series so the pause deviation does not create a large topological feature in the original time series. 
This is reflected in the persistence diagrams and landscapes of \Cref{fig:case-study-persistence} as well; the original time series diagram and landscape (A) show only $3$ significant topological features, while the diagram and landscape of the sliding window embedding (B) and moving average filter (C) show $4$. 

The sliding window embedding and moving average filter both smooth the input data and detect the pause behavior, but they are qualitatively different. One piece of geometric information that the sliding window embedding retains that the moving average filter does not is the direction of the time series. Consider plotting a time series and the corresponding reverse-chronological-order time series in the same ambient space. The points of the original time series would align exactly with its reverse, so the two corresponding point clouds are the same. This is also true of a moving average filter on that time series. The sliding window embedding, however, can have distinct point clouds. 

Retaining the direction of time in a time series is particularly important for data that has certain types of symmetry. A natural occurrence of such data is the metronome data in \Cref{eg:metronome}. Because the data in this time series follows a path and then backtracks along that path, it never produces a loop with any significant persistence. There is no loop in the moving average filter of the data, either.

\subsubsection{Interpretability: Mapping persistence features to (synthetic) behaviors}
\label{sec:representative-cycles}

\begin{figure}[H]
  \begin{center}
    \begin{tabular}{ccccc}
      \multicolumn{1}{l}{(A)}\\
      forward & transition & backward $1$ & pause & backward $2$ \\
      \hightrimmovespace[width=\fifthwidth]{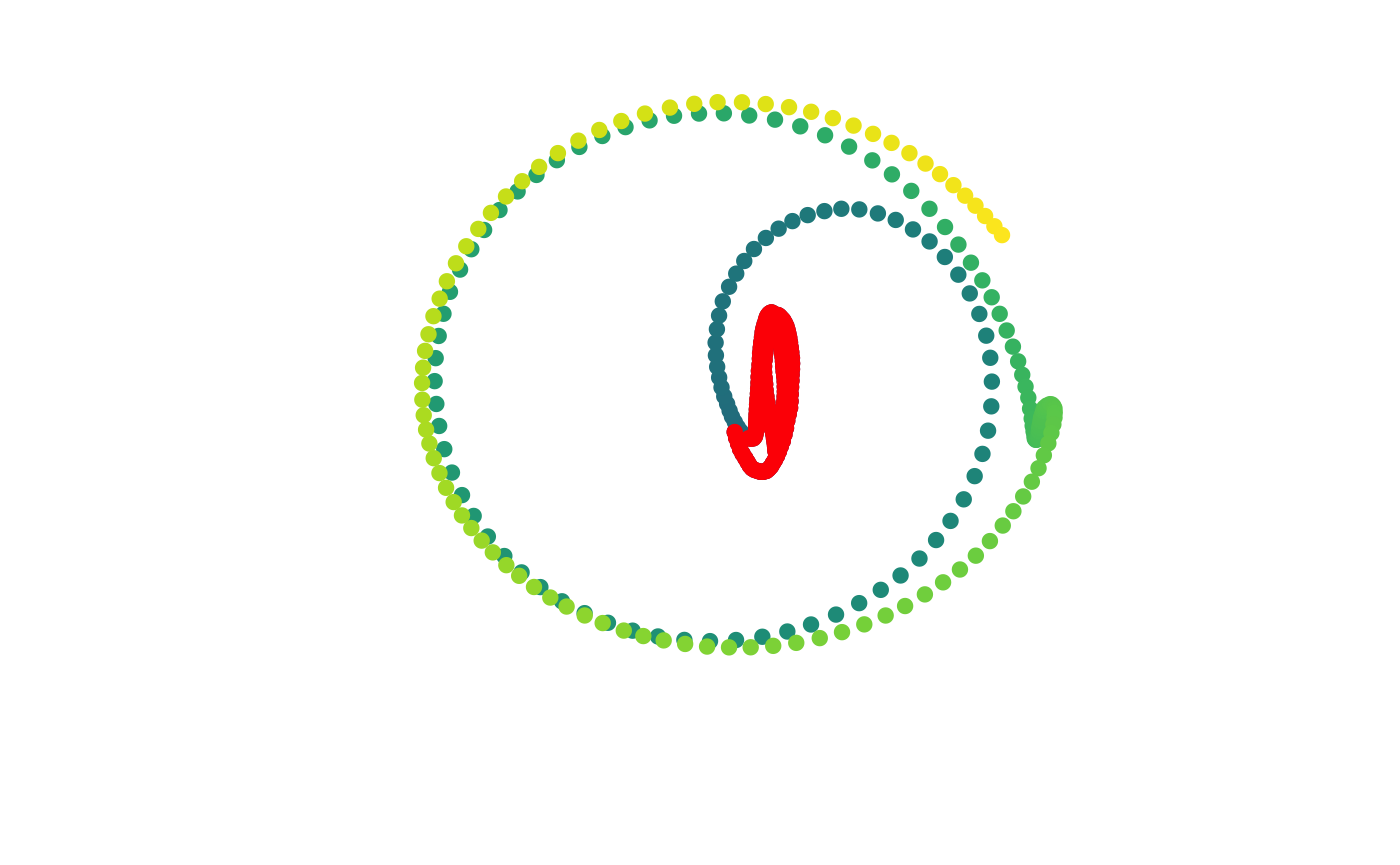}
      &\hightrimmovespace[width=\fifthwidth]{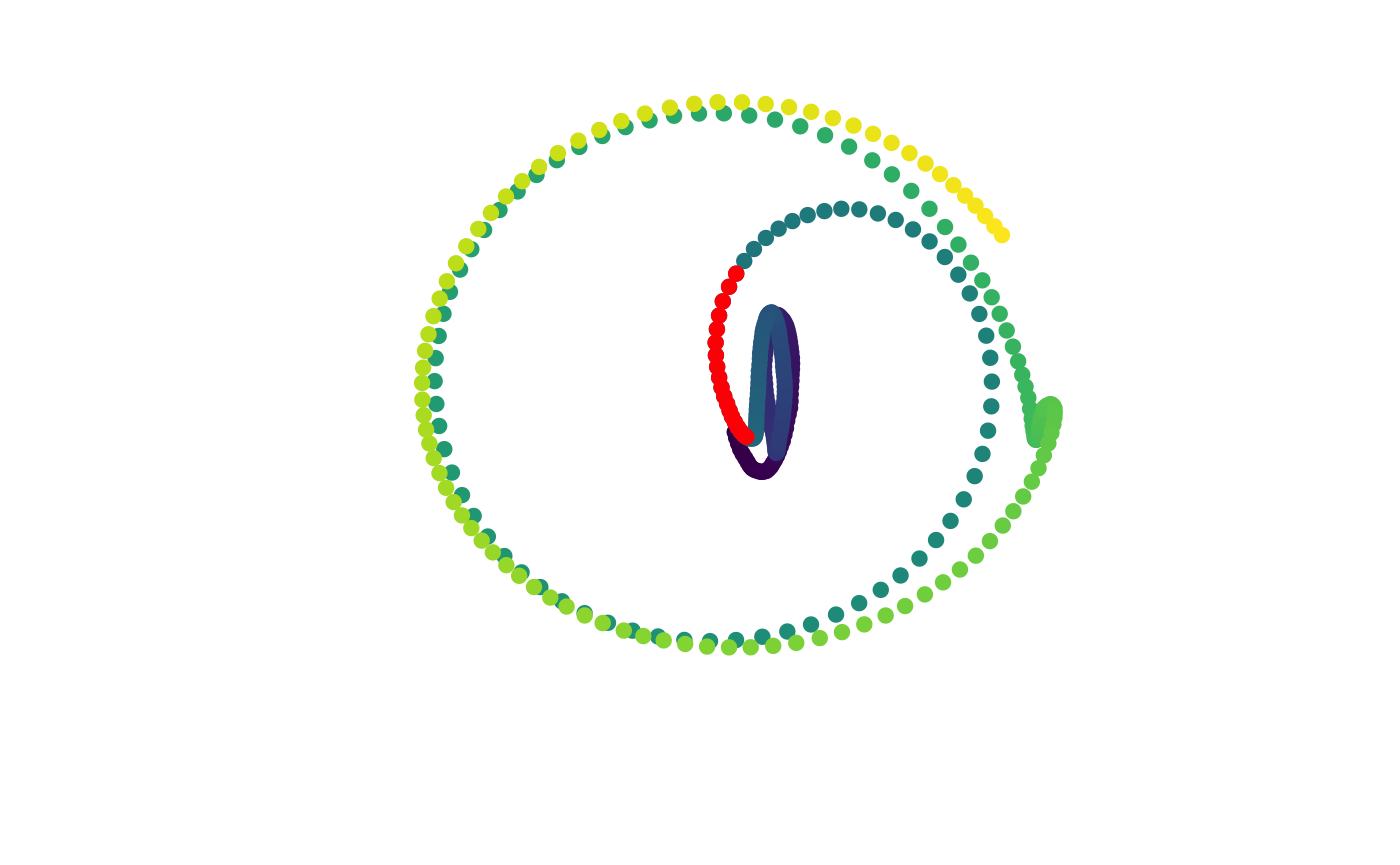}
      &\hightrimmovespace[width=\fifthwidth]{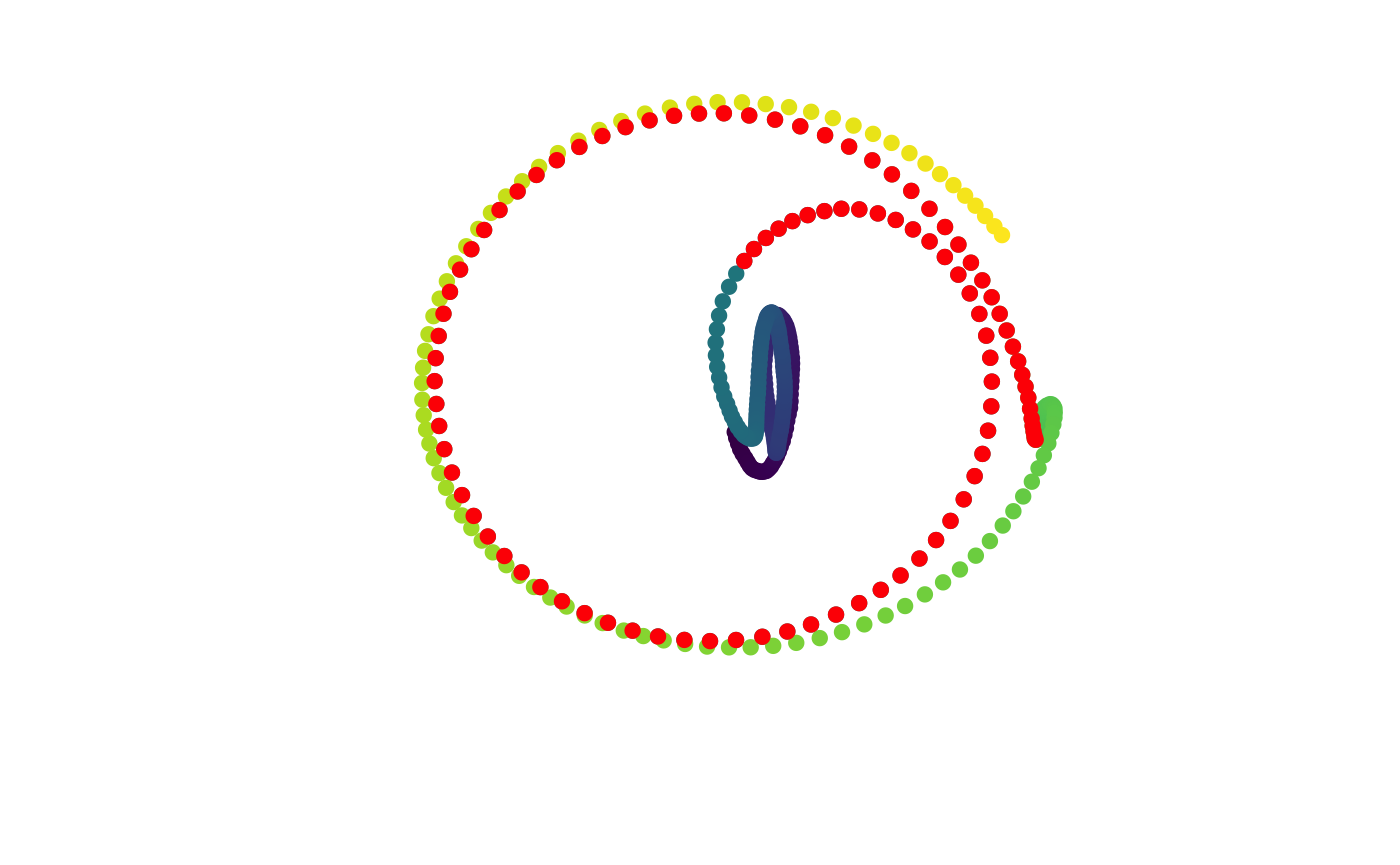}
      &\hightrimmovespace[width=\fifthwidth]{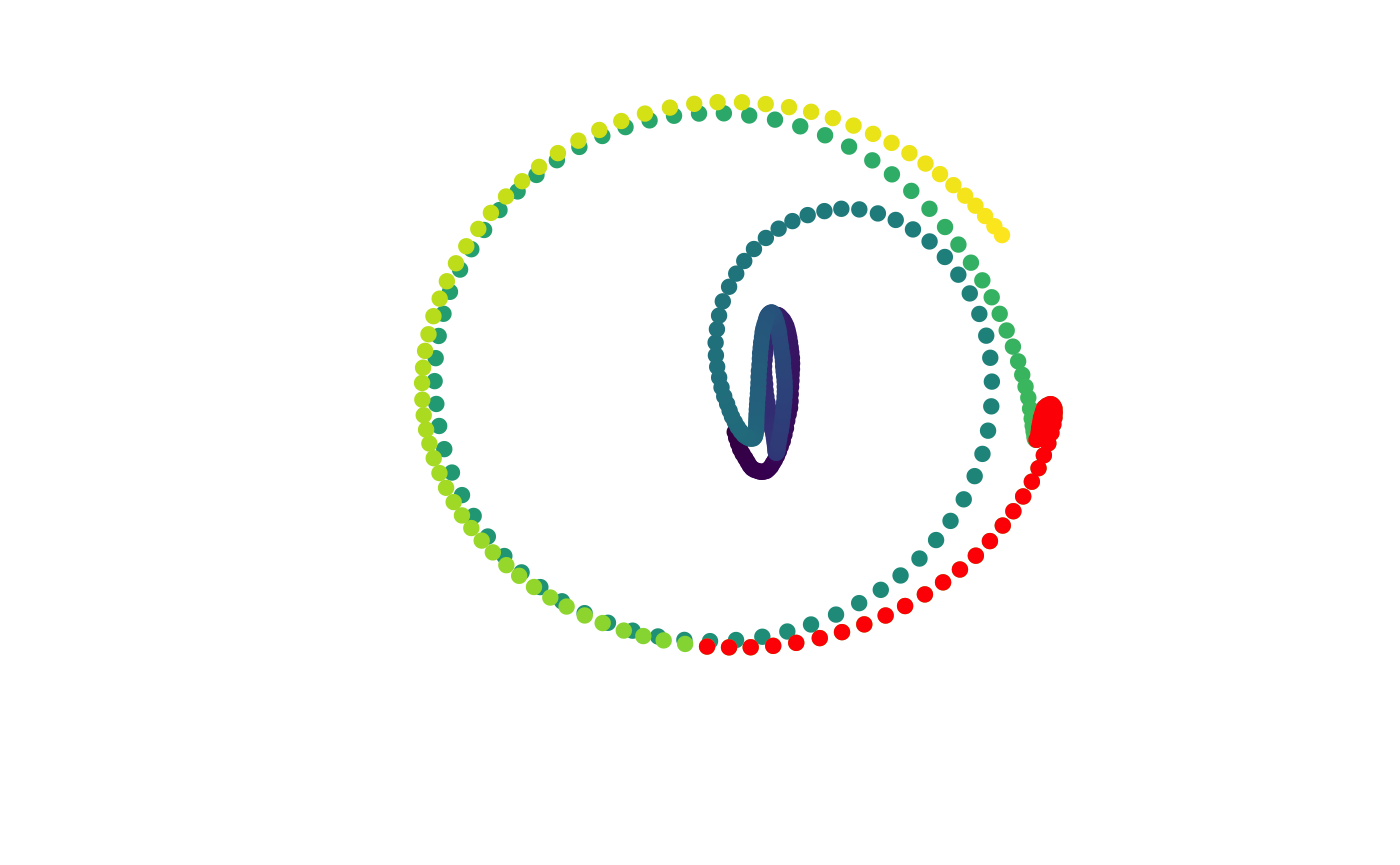}
      &\hightrimmovespace[width=\fifthwidth]{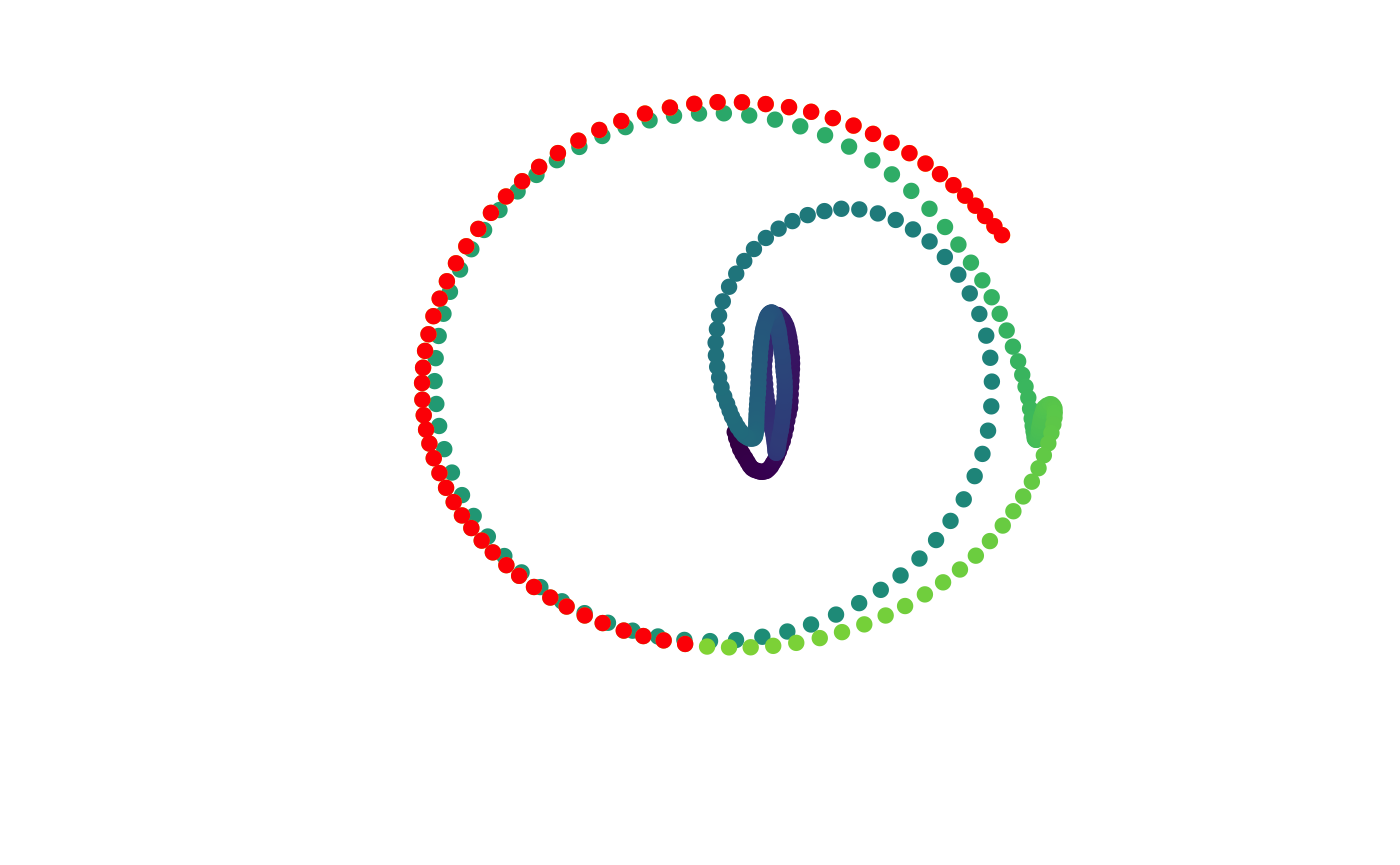}
    \end{tabular}
    \\
    \begin{tabular}{cccc}
      \multicolumn{1}{l}{(B)}\\
      forward & transition & backward & pause \\
      \hightrimmovespace[width=\quarterwidth]{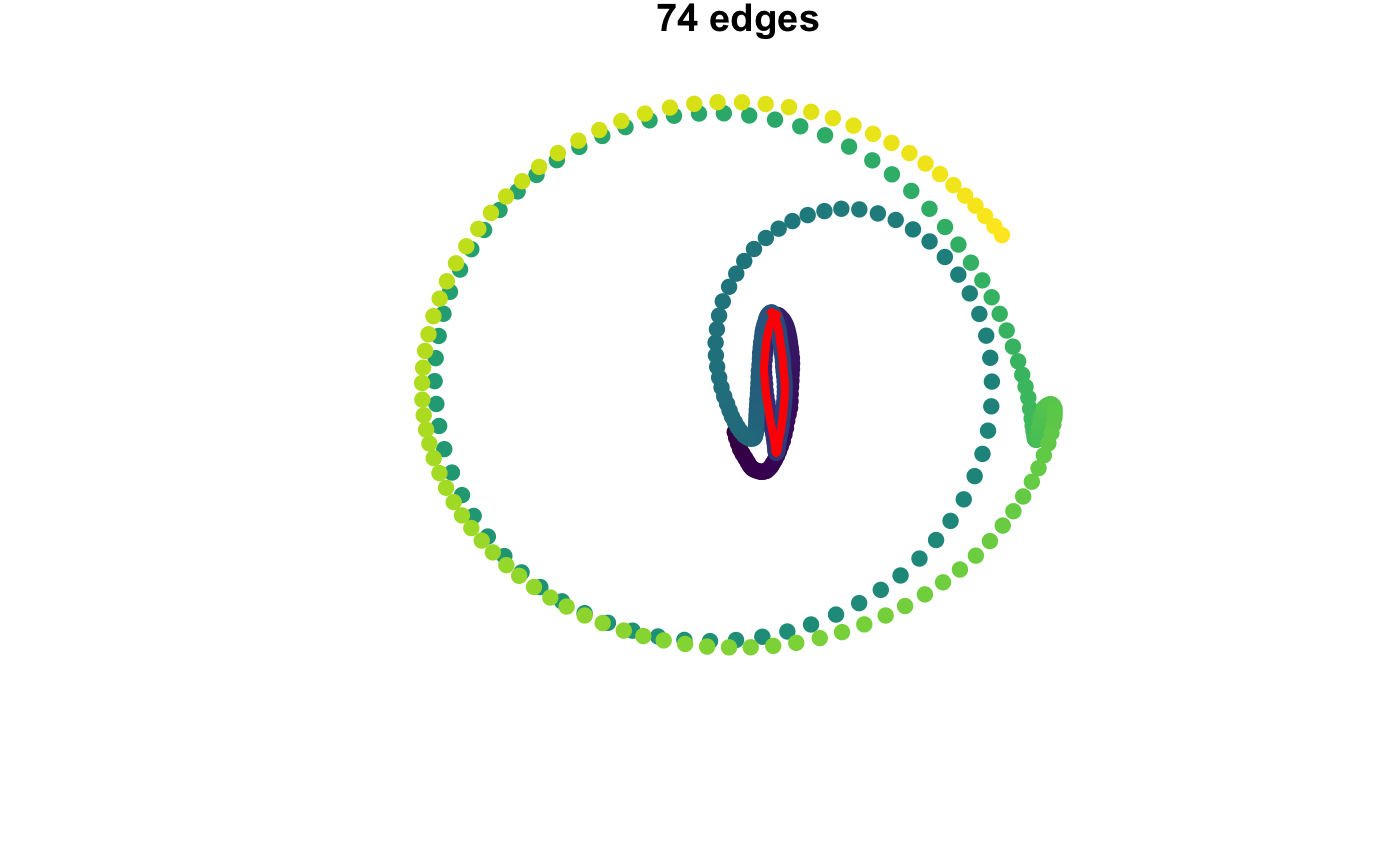}
      &\hightrimmovespace[width=\quarterwidth]{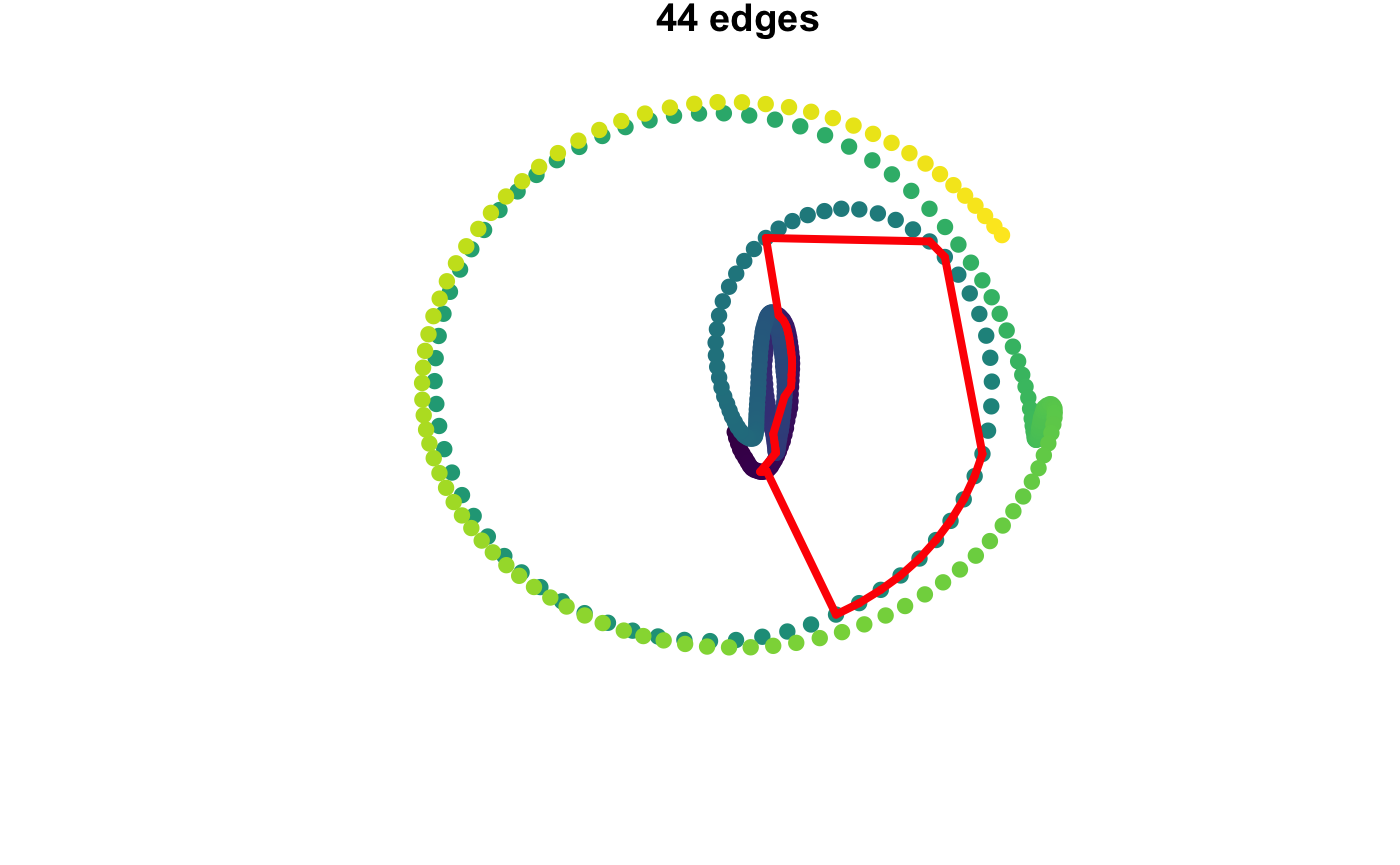}
      &\hightrimmovespace[width=\quarterwidth]{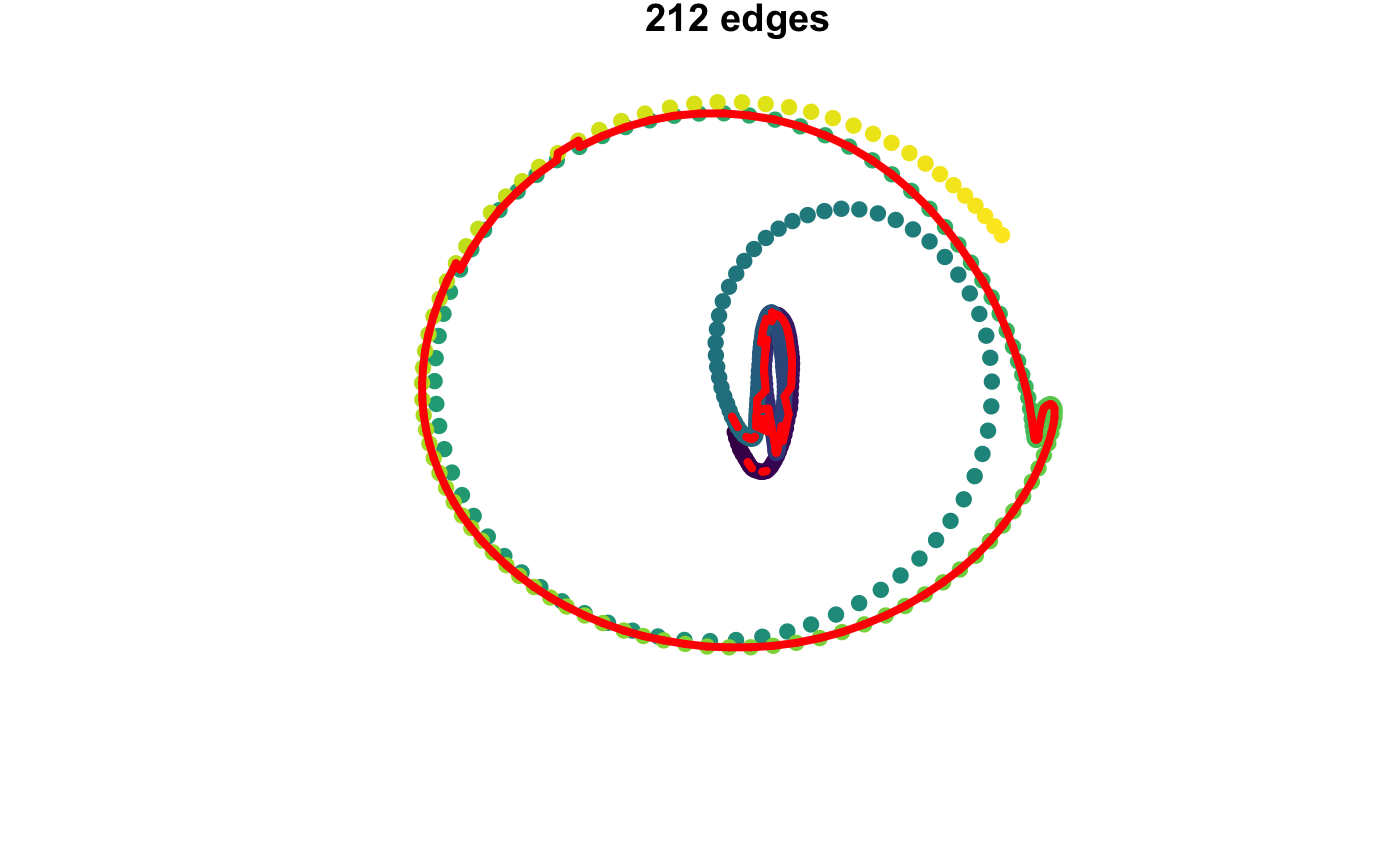}
      &\hightrimmovespace[width=\quarterwidth]{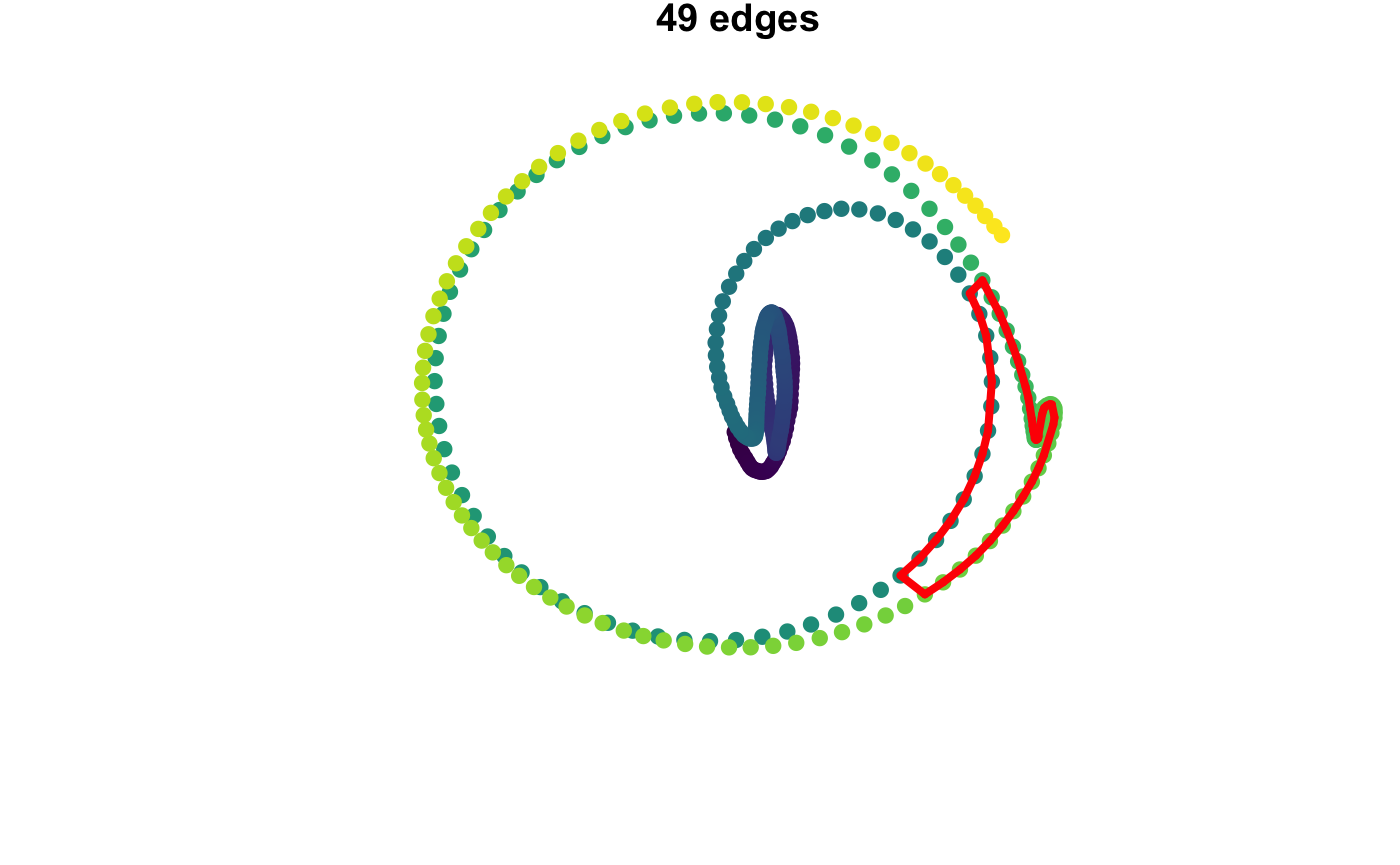}
      \\
      \multicolumn{1}{l}{(C)}\\
      \trimmedwithaxesgraphic[width=\quarterwidth]{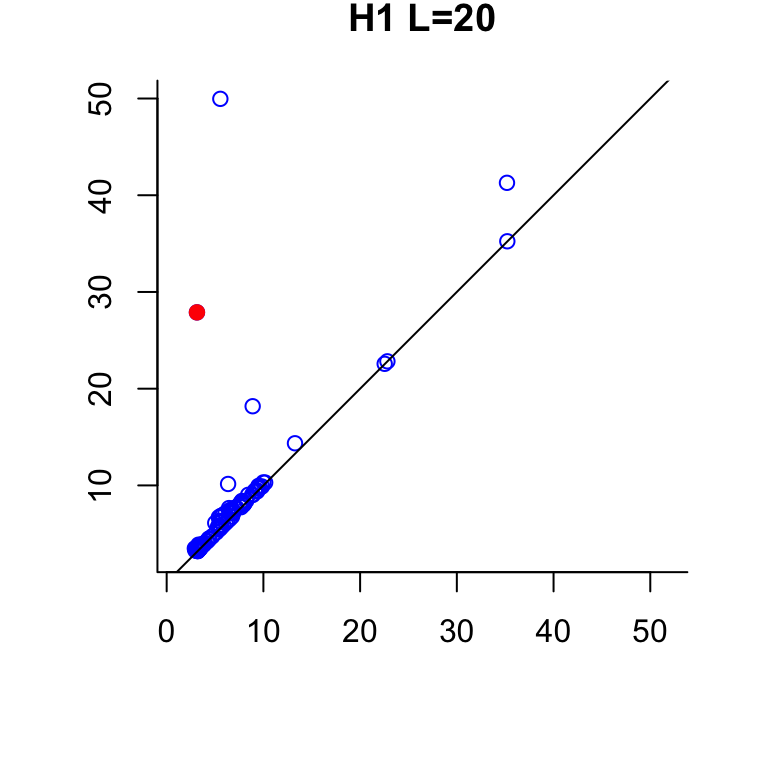}
      &\trimmedwithaxesgraphic[width=\quarterwidth]{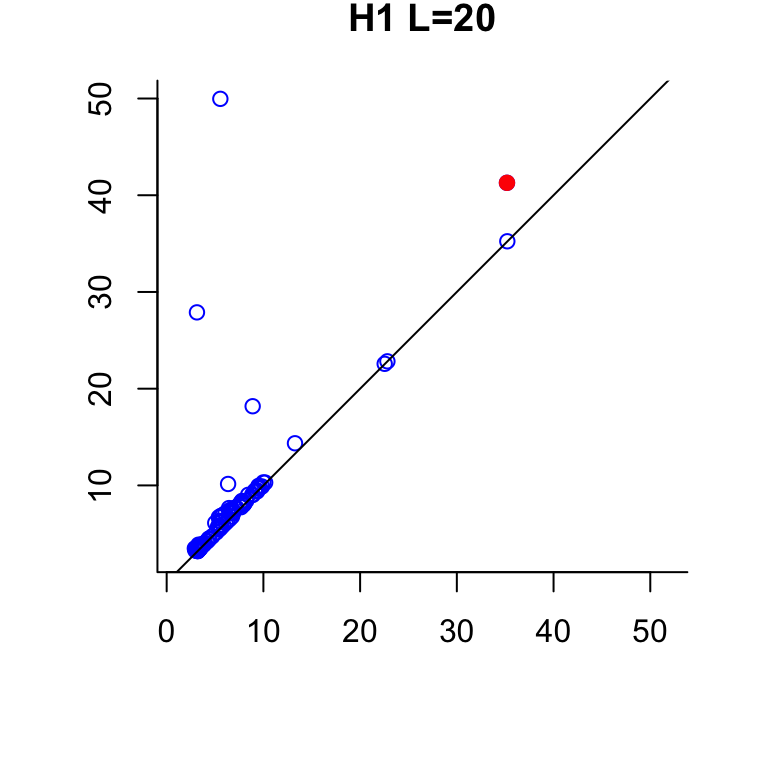}
      &\trimmedwithaxesgraphic[width=\quarterwidth]{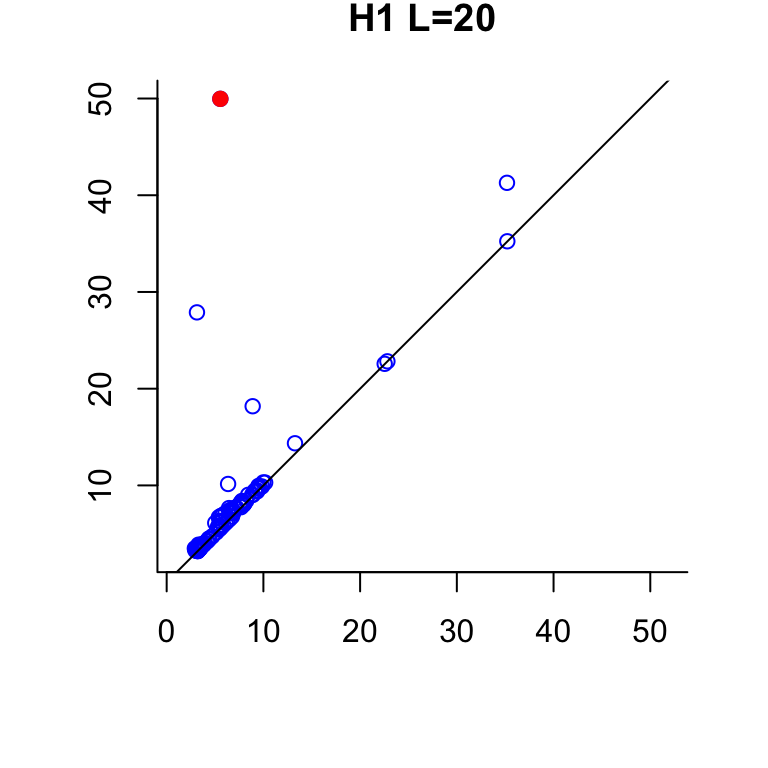}
      &\trimmedwithaxesgraphic[width=\quarterwidth]{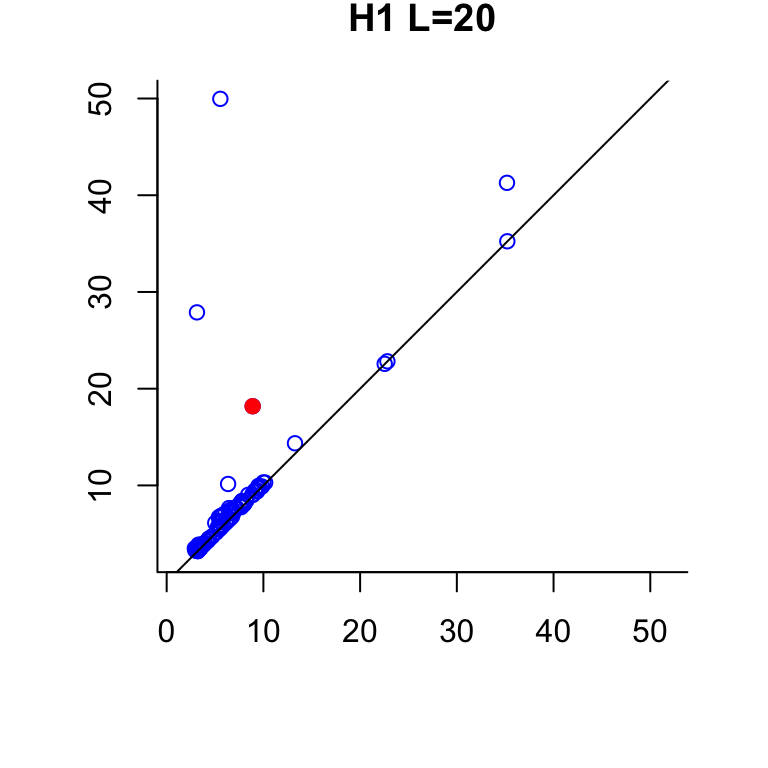}
      \\
      \multicolumn{1}{l}{(D)}\\
      \includegraphics[trim = 150 150 150 100,clip, width=\fifthwidth]{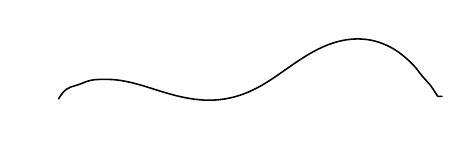}
      & \includegraphics[trim = 150 150 150 100,clip, width=\fifthwidth]{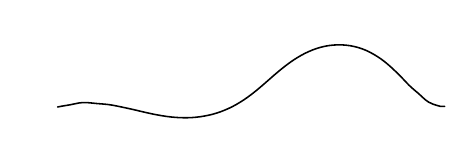}
      & \includegraphics[trim = 150 150 150 100,clip, width=\fifthwidth]{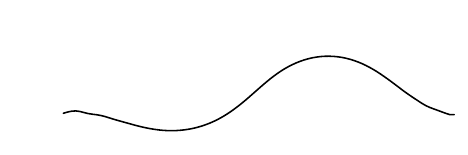}
      & \includegraphics[trim = 150 150 150 100,clip, width=\fifthwidth]{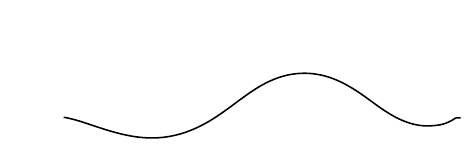}
    \end{tabular}
  \end{center}
  \caption[Points and cycles corresponding to behaviors ]{(A) Points in the sliding window embedding $\tau^{20}$ that correspond to each of the labeled behaviors are highlighted. (B) The representative cycles with longest persistence from automated persistence software correspond to specific behaviors. (C) The persistence diagram with the homology class corresponding to the above cycle representative highlighted. (D) Still frames from a looping video of forward crawling data. This synthetic data was constructed from the forward representative cycle in (B). See the supplemental materials (\Cref{sec:supp-material}) for the full video. }
  \label{fig:case-study-behaviors}
\end{figure}

We computed representative cycles for persistent homology classes for each of the longest-persisting topological features
in the sliding window embedding
using Dionysus \cite{Dionysus}.
These are shown in \Cref{fig:case-study-behaviors} (B). The homology classes that correspond to each of these representative cycles are highlighted
in \Cref{fig:case-study-behaviors} (C). 
We remark that instead of the representative cycles produced by Dionysus one may want to use (approximate) shortest cycle representatives
\cite{Dey2019,Dey2018,Obayashi2018a,Erickson:2012}.

One of the benefits of using persistence for behavior analysis is that these representative cycles give a direct translation from persistence back into \emph{C. elegans} behavior. Each point in the cycle corresponds to $\wl$ poses and these points have a defined sequence in the cycle. The cycle lacks a direction (which way is forward in time vs which way is backward) but in many cases a direction can be inferred by subsets of the sequence that correspond to contiguous sequences in the original time series. Given all this data and a way to combine $\wl$ poses into one ``average'' pose, we can construct synthetic, periodic behavior data from representative cycles. An example of frames from such a video is shown in \Cref{fig:case-study-behaviors} (D) and the corresponding video is available in the supplemental materials (\Cref{sec:supp-material}). 

\newpage
\subsection{Experimental results} \label{sec:experimental-results}

We present our analysis of an experiment where \emph{C. elegans} are submerged in solutions with varying viscosities. The viscosity of the solution is correlated with how much methylcellulose is added and experimental conditions are labeled with their methylcellulose content, usually in order from low to high methylcellulose and viscosity. 

\begin{figure}[H] 
  \begin{center}
    \begin{tabular}{cccc}
      $0.5\%$  & $1\%$  
      & $2\%$  & $3\%$ \\[0.5em]
      \multicolumn{1}{l}{(A)}\\
      \trimmedgraphic[width=\quarterwidth]{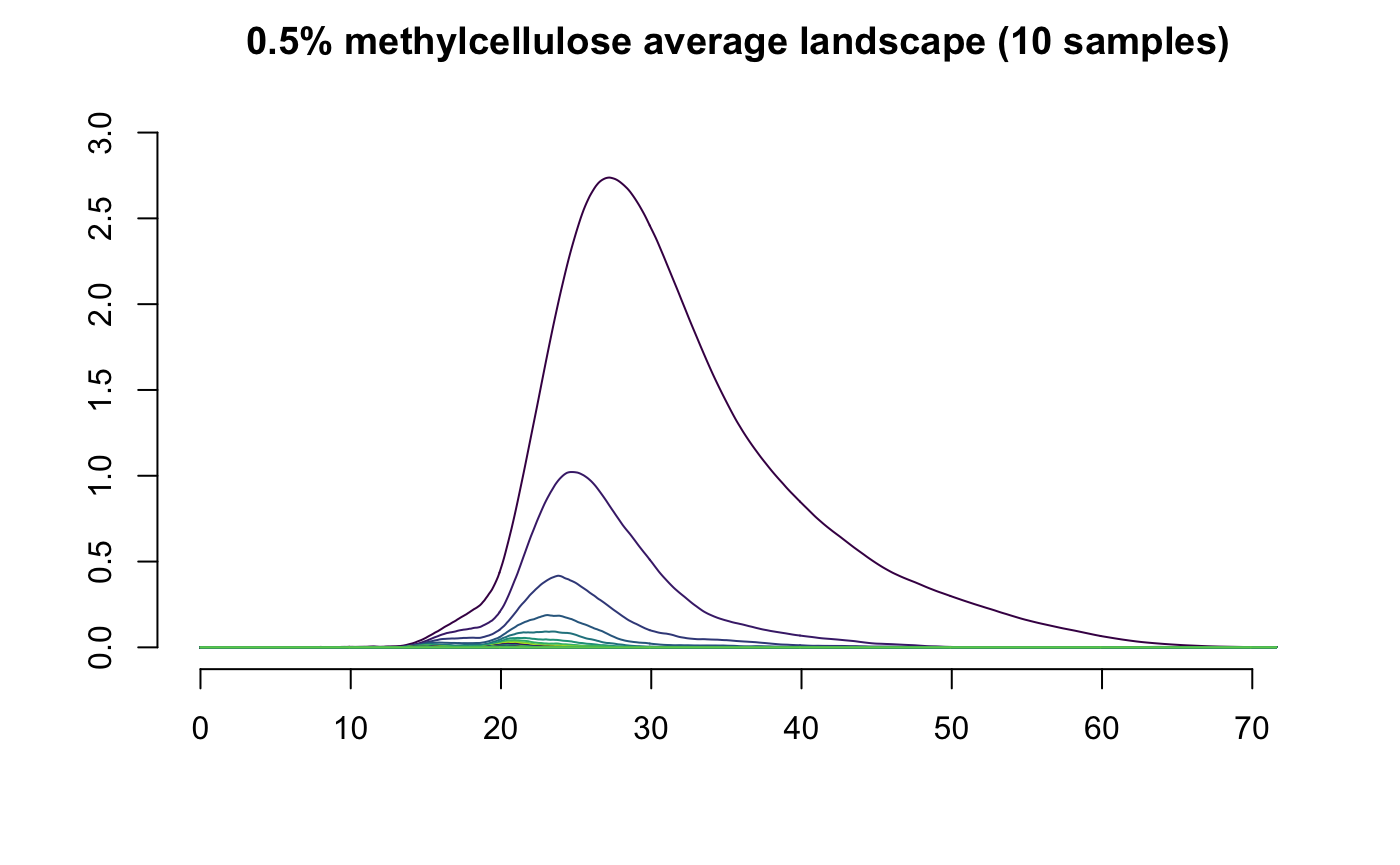}
      &\trimmedgraphic[width=\quarterwidth]{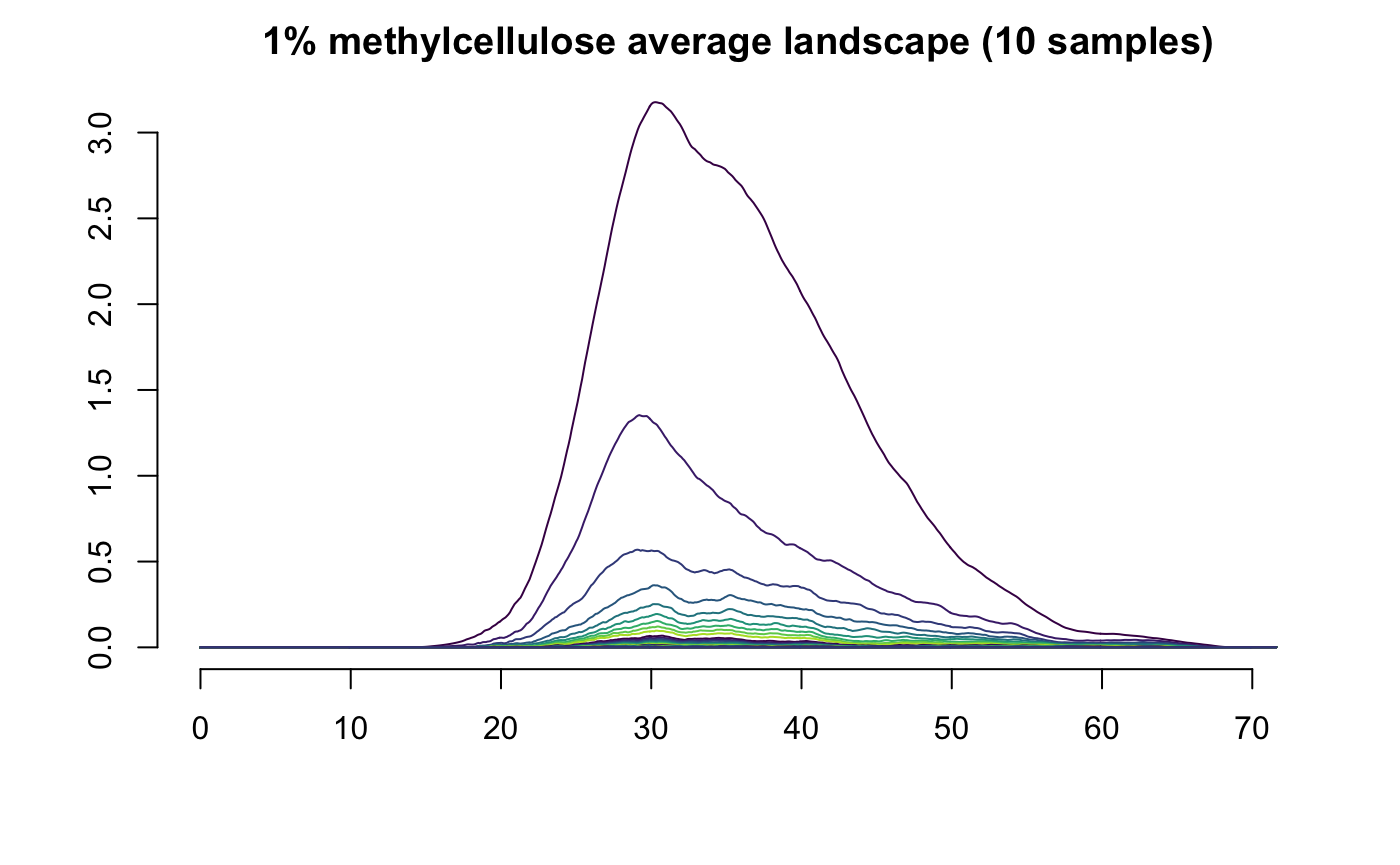}
      &\trimmedgraphic[width=\quarterwidth]{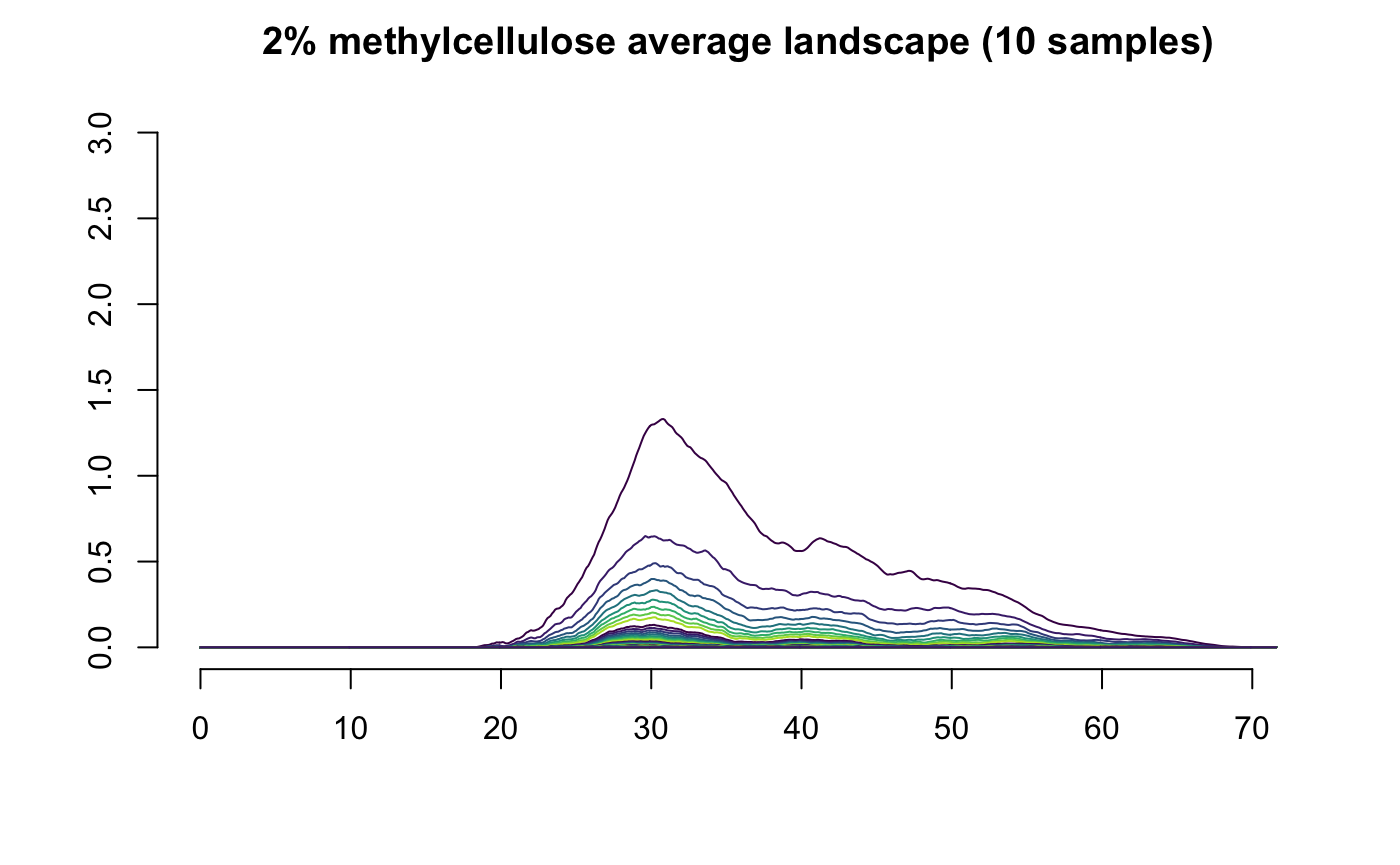}
      &\trimmedgraphic[width=\quarterwidth]{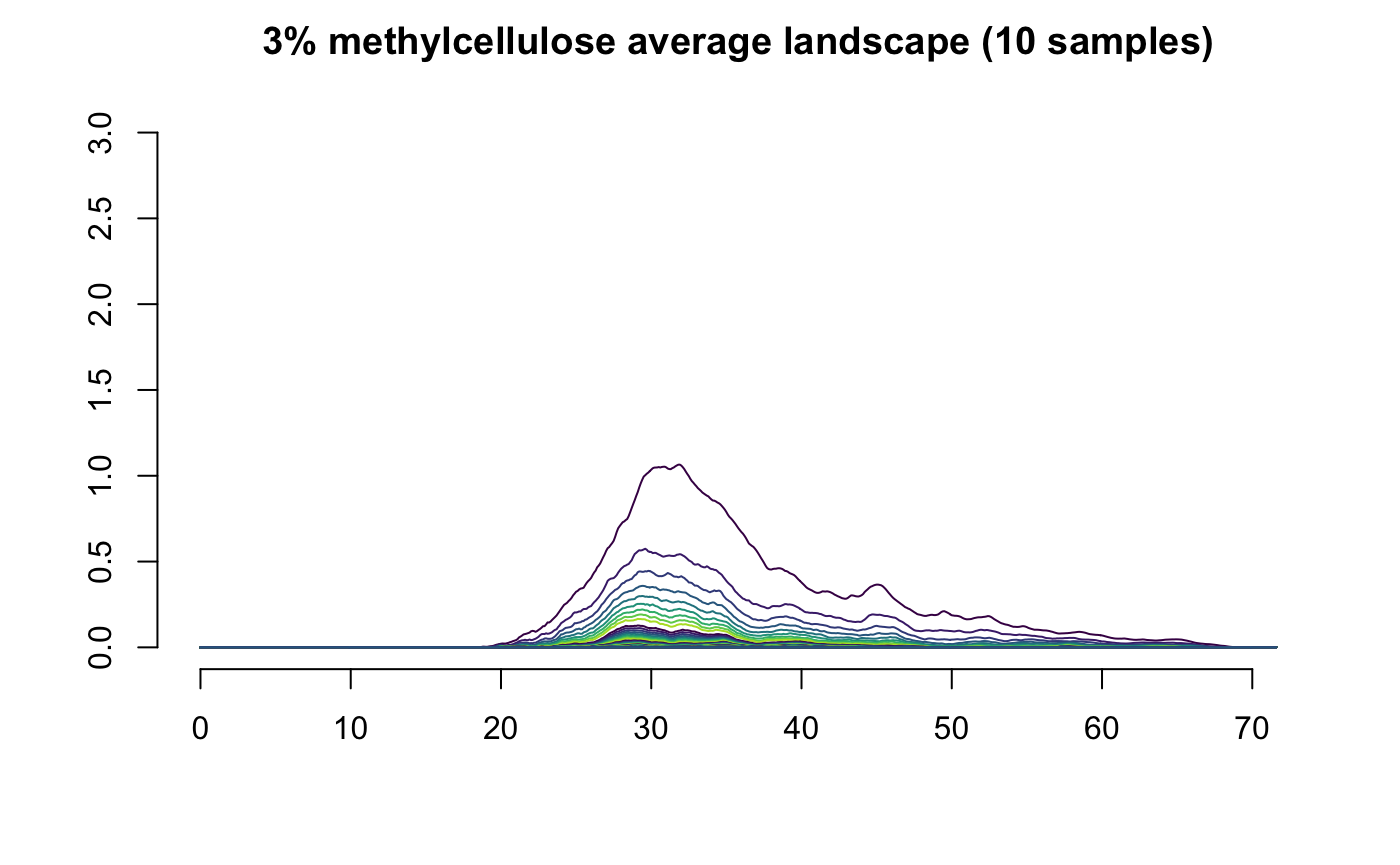}
      \\[1.5em]
      \multicolumn{1}{l}{(B)}\\
      \trimmedgraphic[width=\quarterwidth]{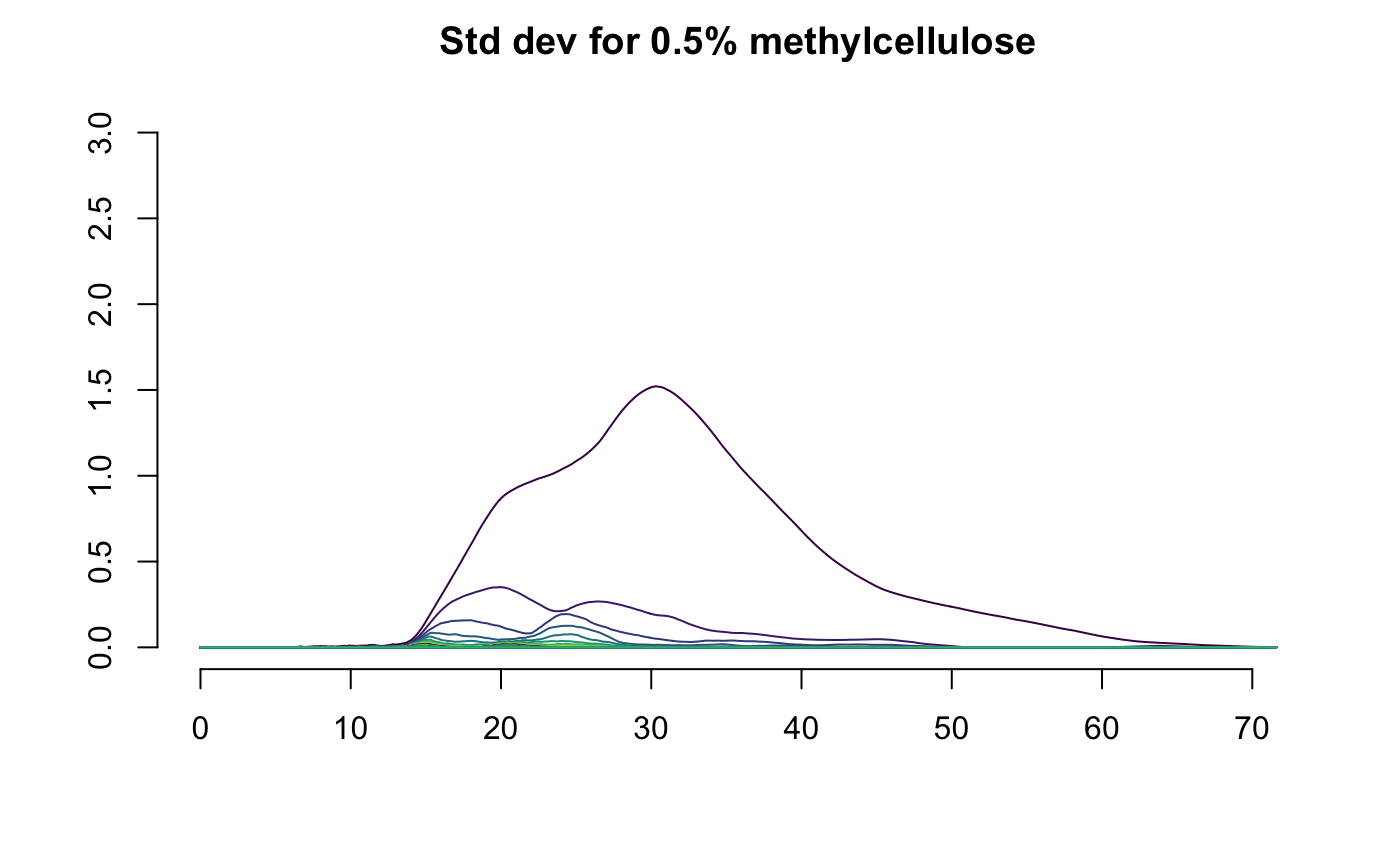}
      &\trimmedgraphic[width=\quarterwidth]{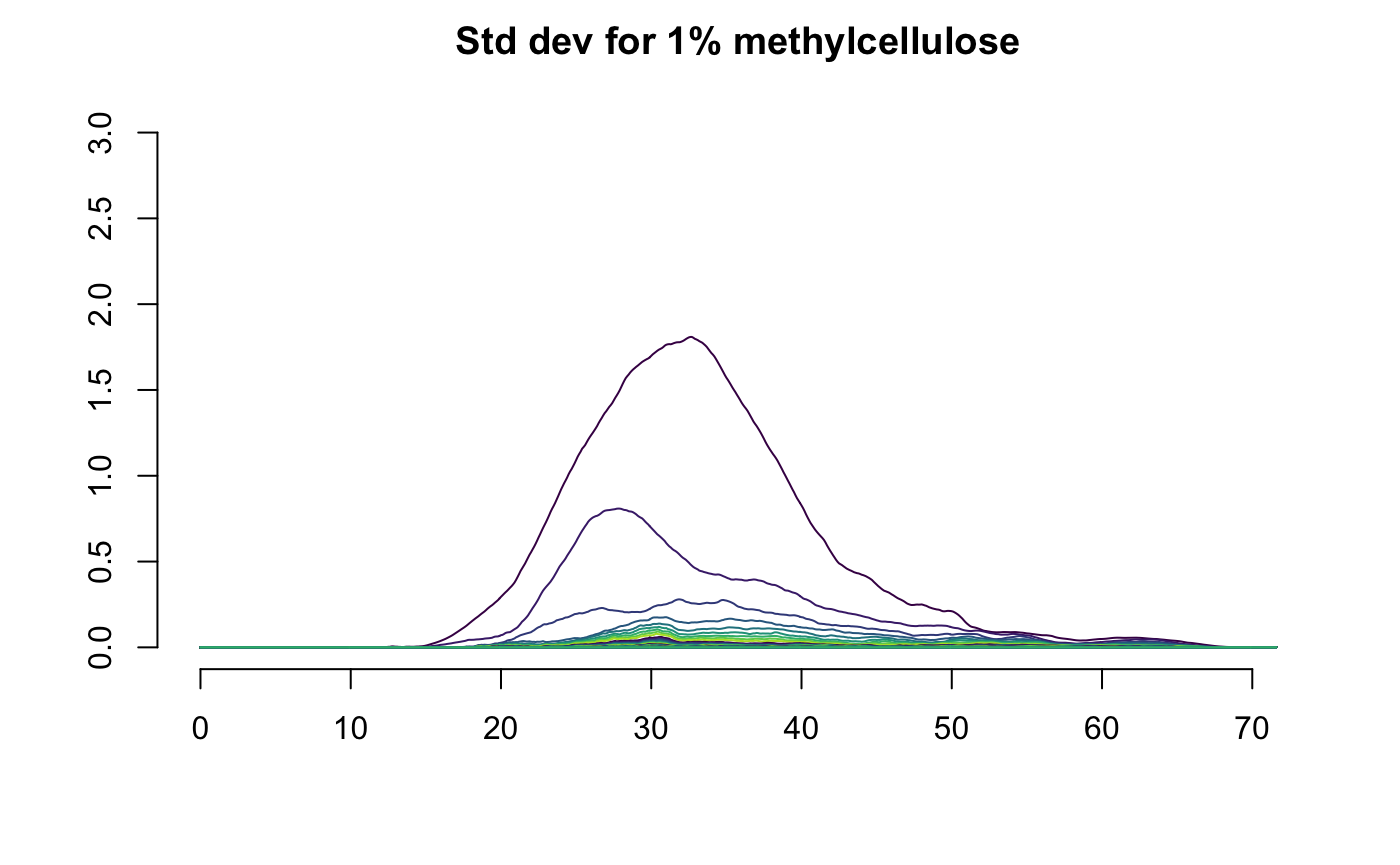}
      &\trimmedgraphic[width=\quarterwidth]{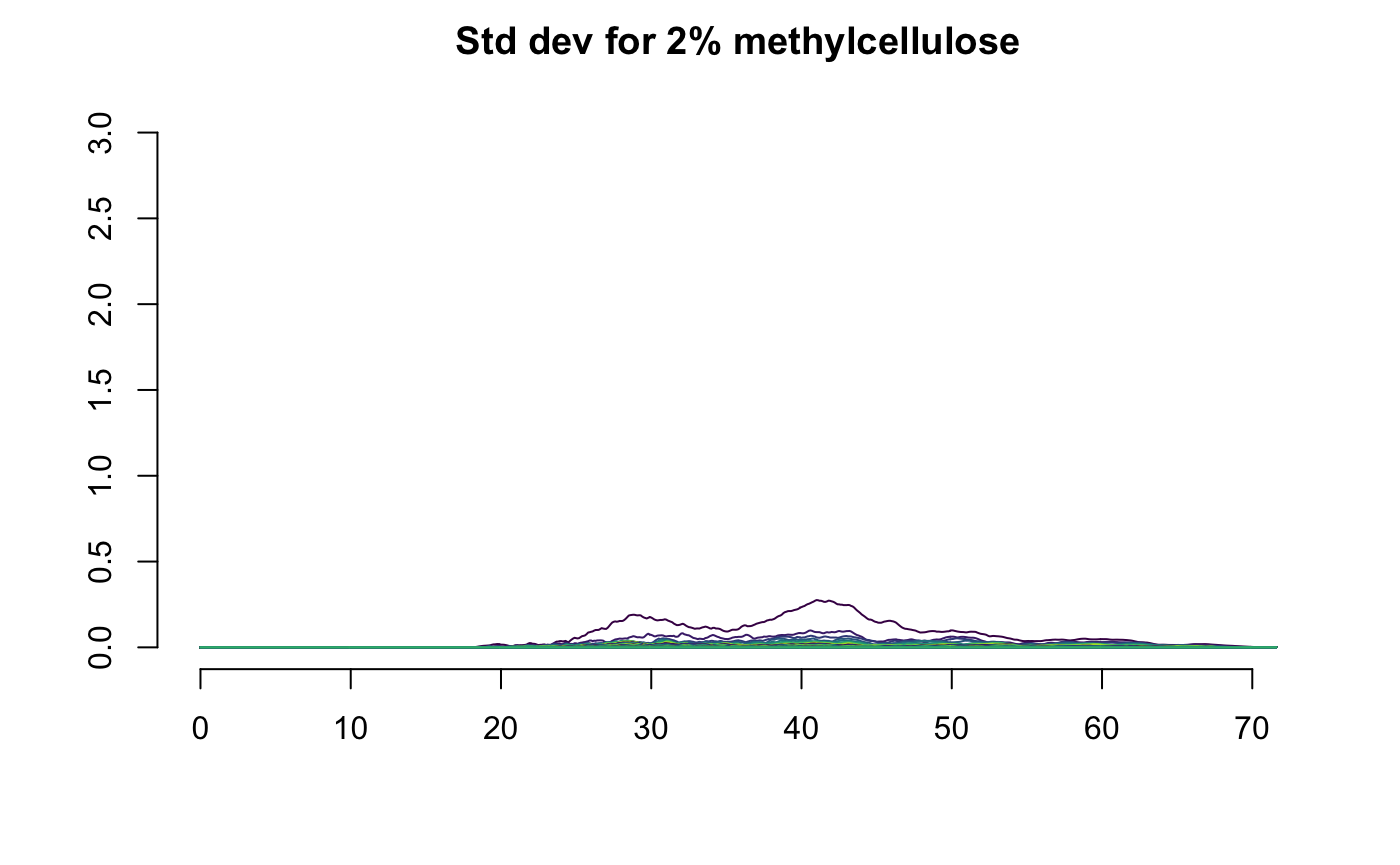}
      &\trimmedgraphic[width=\quarterwidth]{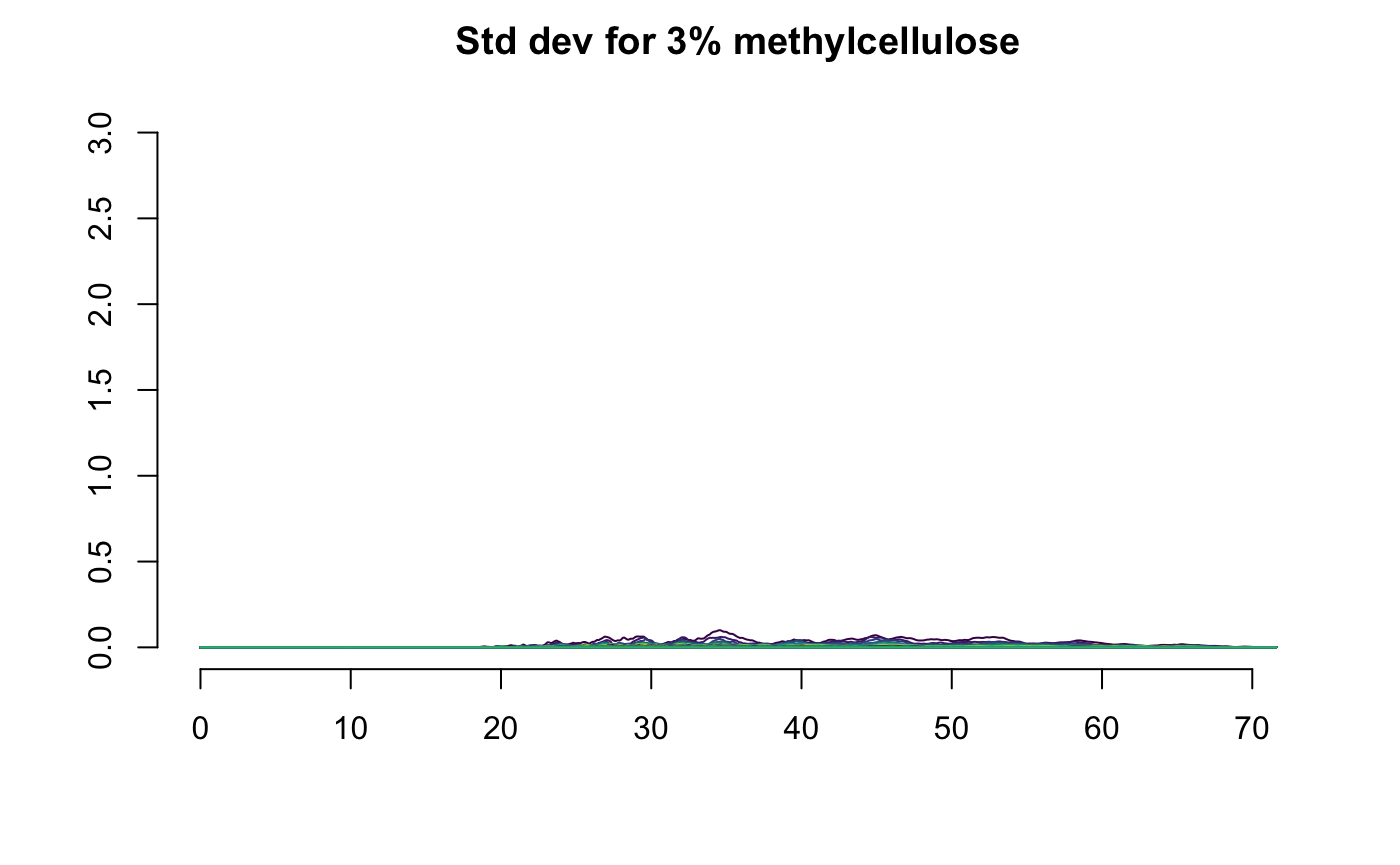}
      \\[1.5em]
      \multicolumn{1}{l}{(C)}\\
      \includegraphics[trim = 60 75 25 25,clip, width=\quarterwidth]{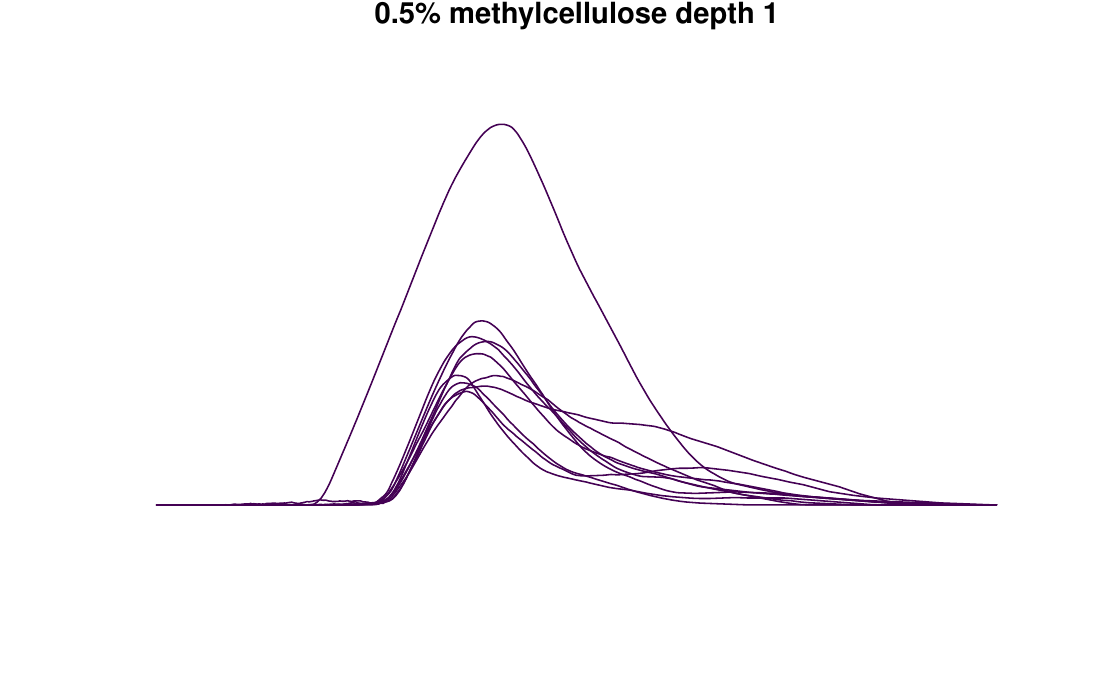}
      & \includegraphics[trim = 60 75 25 25,clip, width=\quarterwidth]{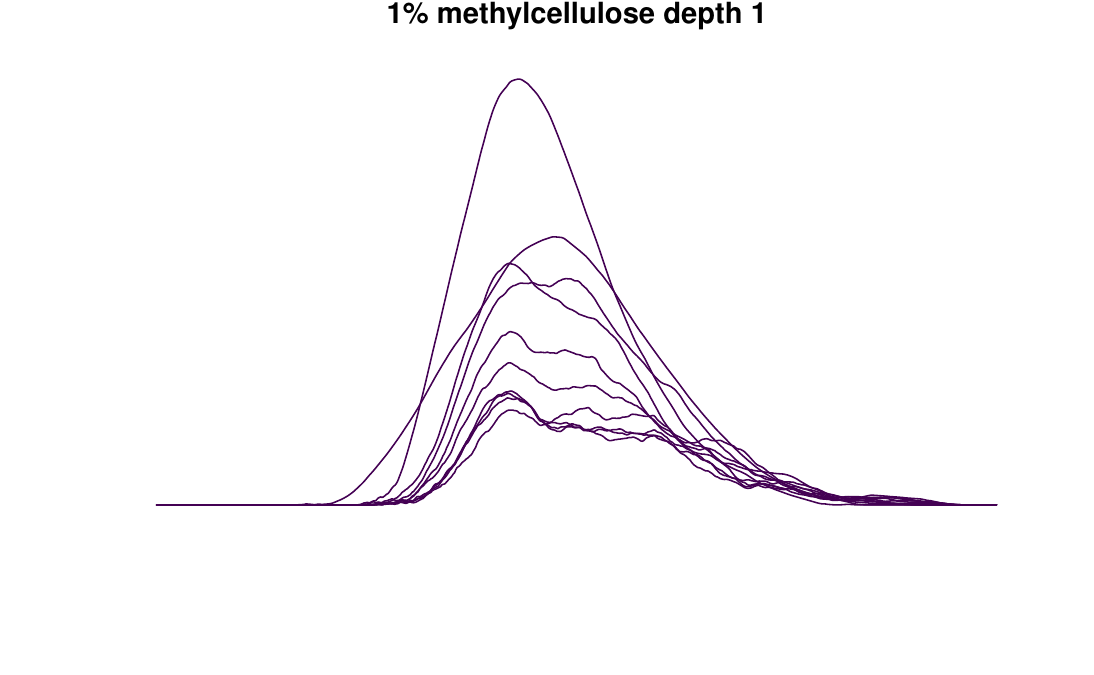}
      & \includegraphics[trim = 60 75 25 25,clip, width=\quarterwidth]{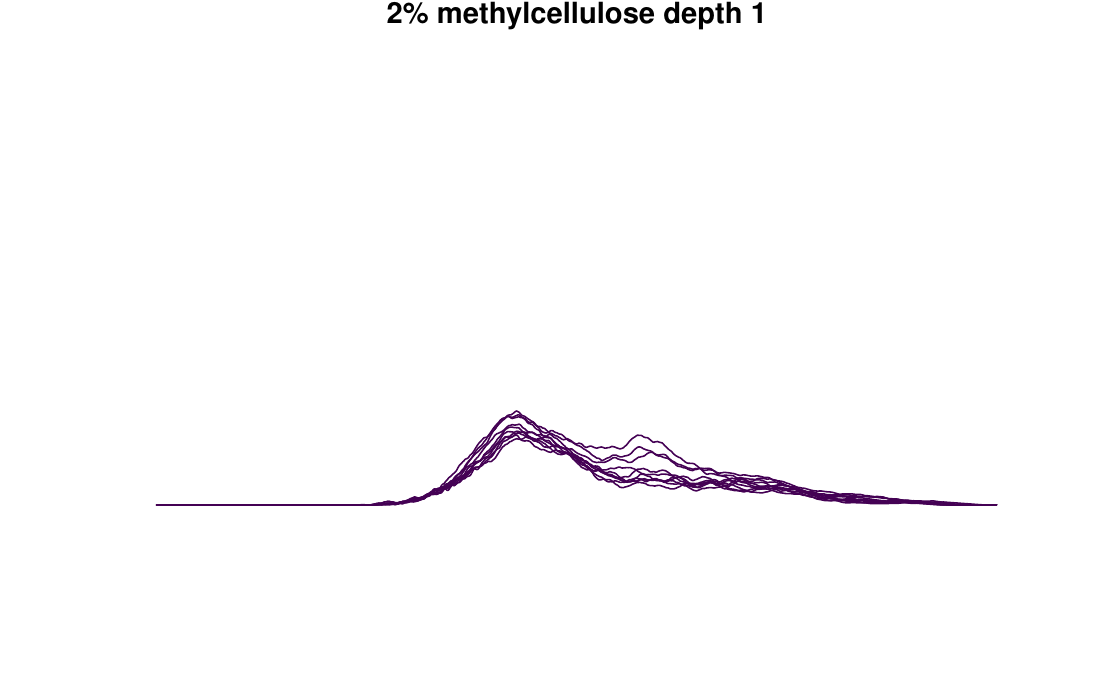}
      & \includegraphics[trim = 60 75 25 25,clip, width=\quarterwidth]{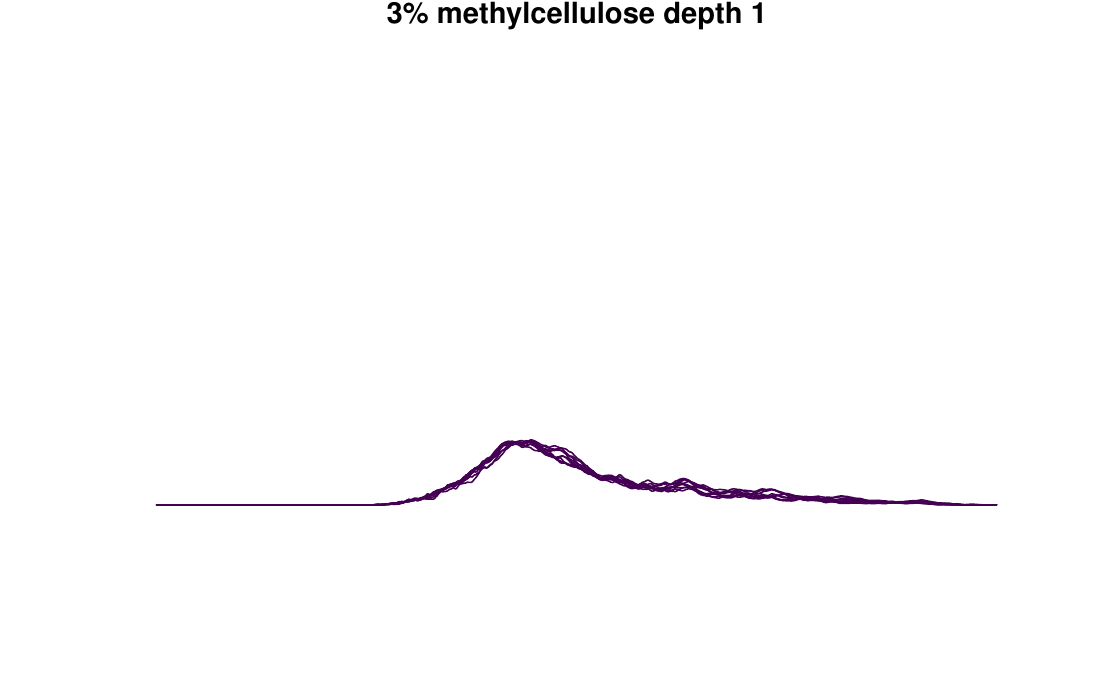}
      \\[1.5em]
      \multicolumn{1}{l}{(D)}\\
      \includegraphics[trim = 60 75 25 25,clip, width=\quarterwidth]{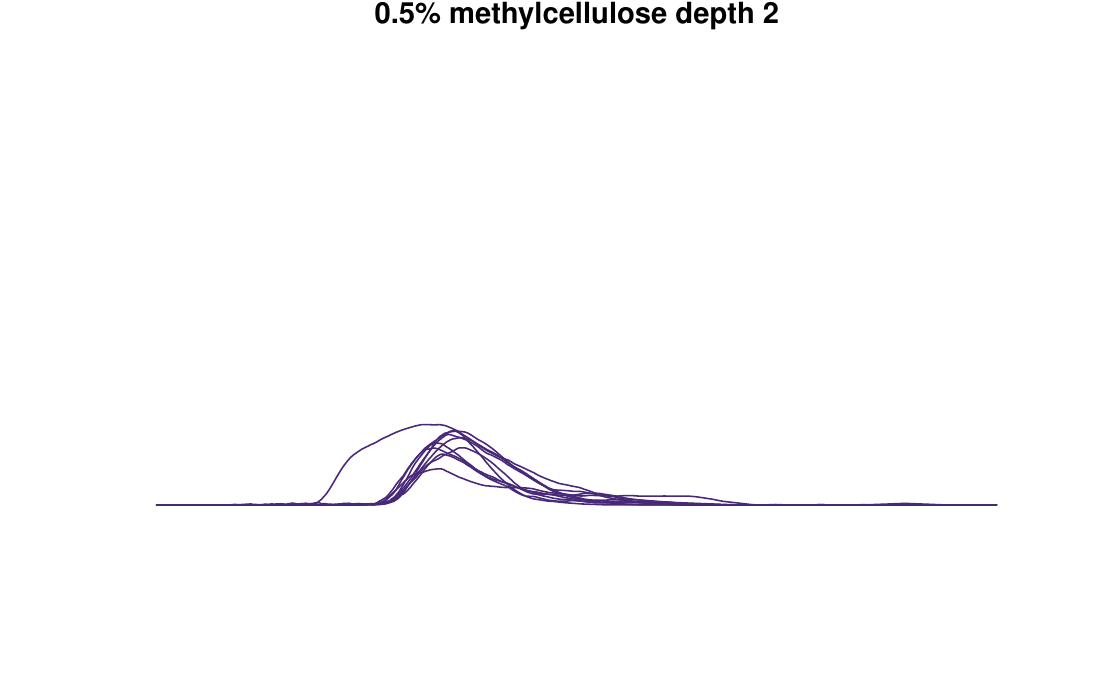}
      & \includegraphics[trim = 60 75 25 25,clip, width=\quarterwidth]{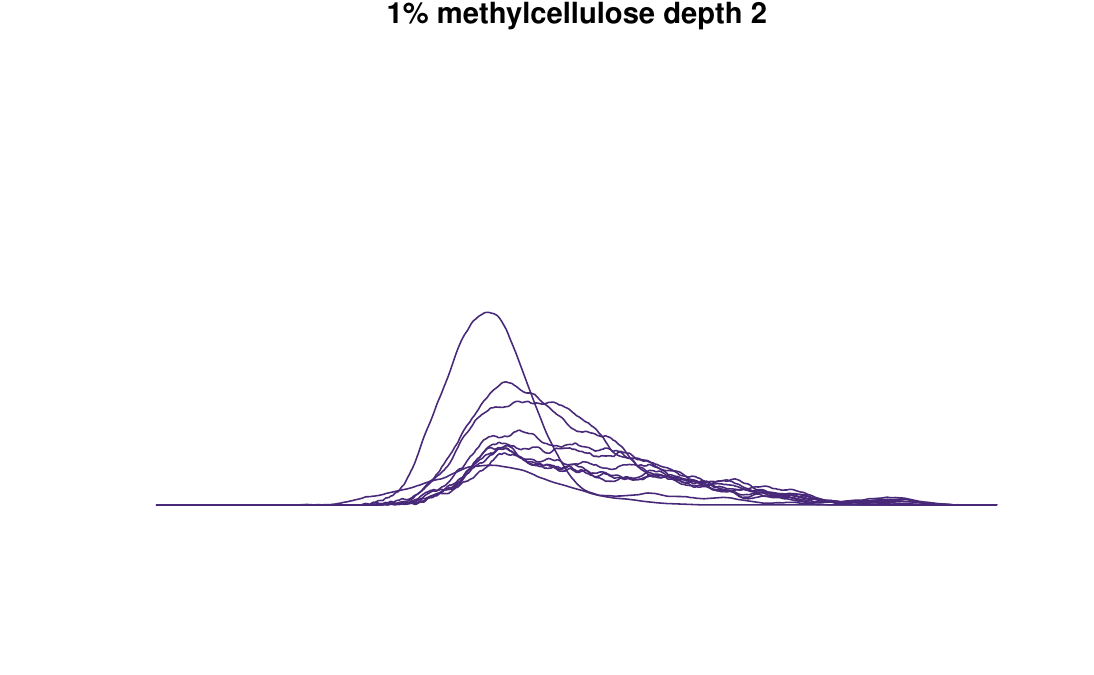}
      & \includegraphics[trim = 60 75 25 25,clip, width=\quarterwidth]{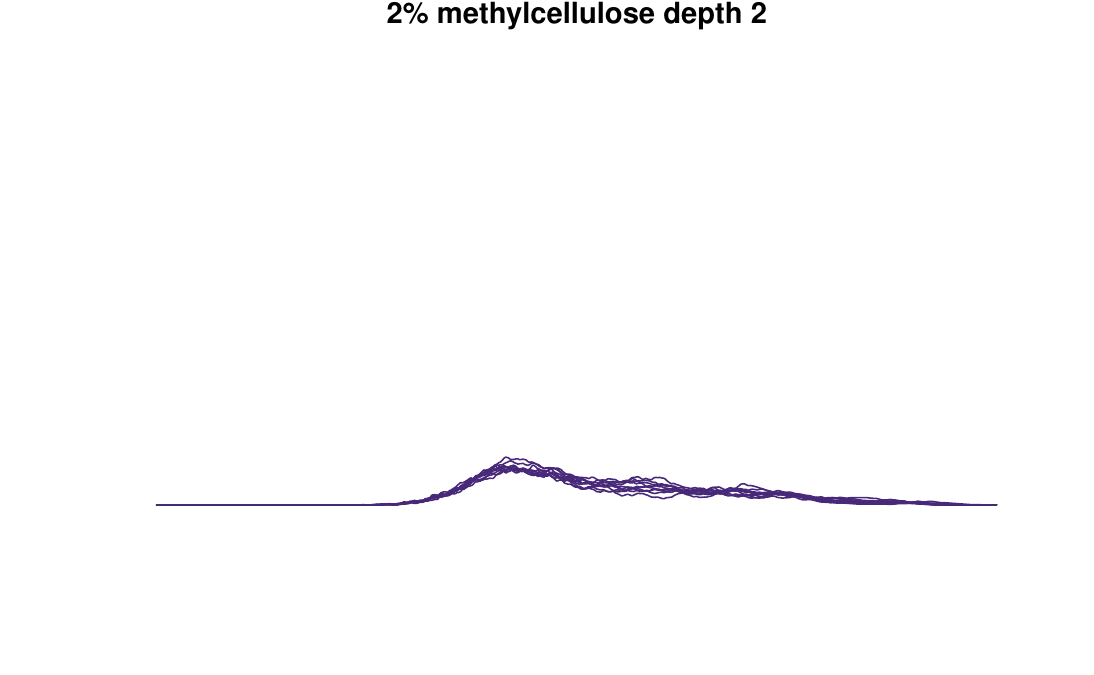}
      & \includegraphics[trim = 60 75 25 25,clip, width=\quarterwidth]{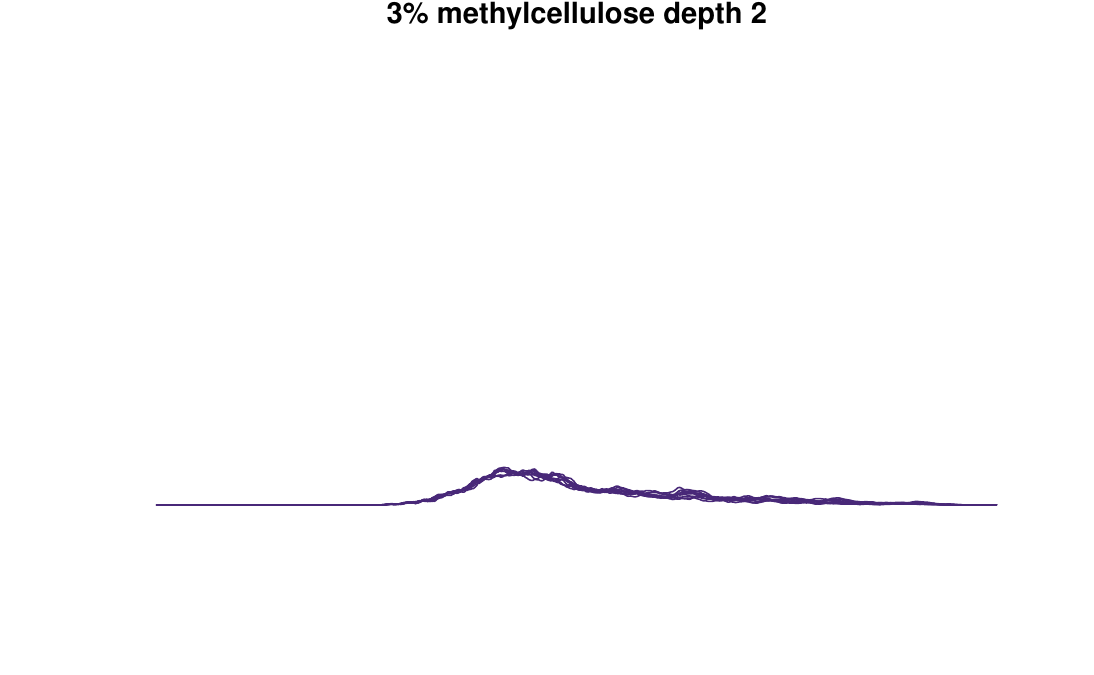}
      \\[1.5em]
      \multicolumn{1}{l}{(E)}\\
      \includegraphics[trim = 60 75 25 25,clip, width=\quarterwidth]{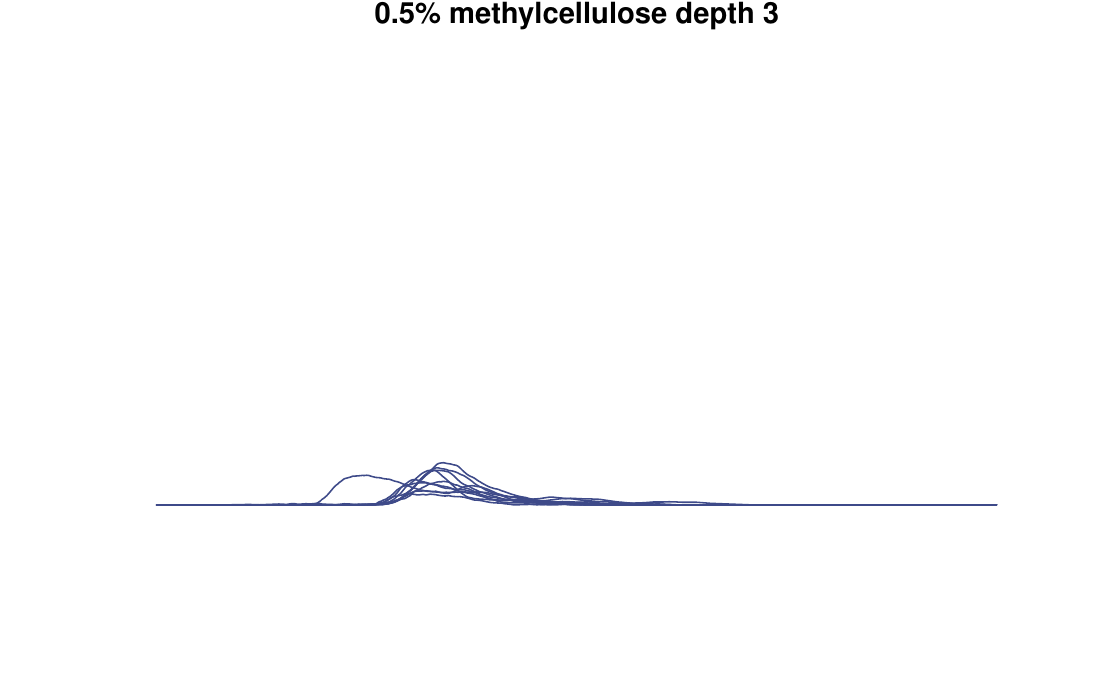}
      & \includegraphics[trim = 60 75 25 25,clip, width=\quarterwidth]{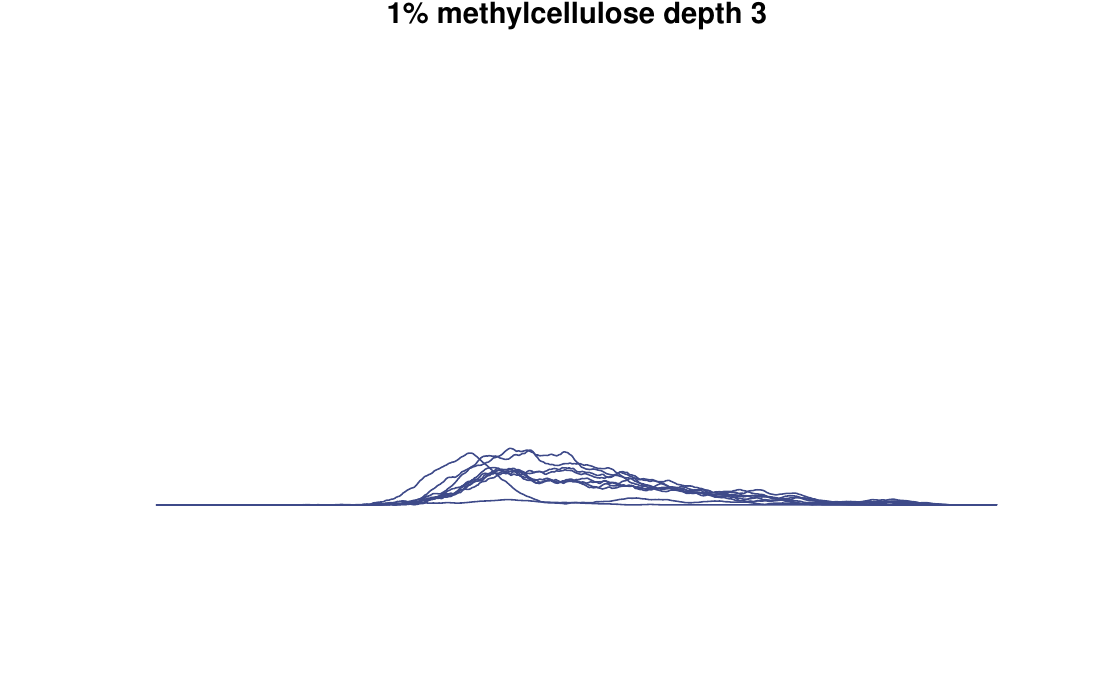}
      & \includegraphics[trim = 60 75 25 25,clip, width=\quarterwidth]{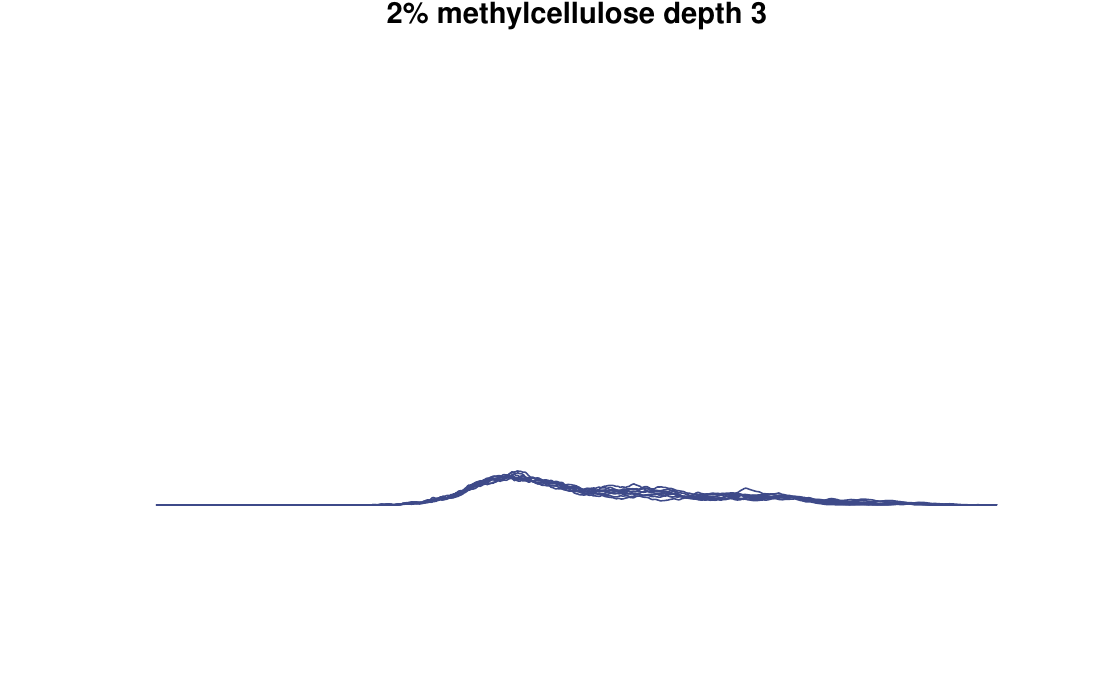}
      & \includegraphics[trim = 60 75 25 25,clip, width=\quarterwidth]{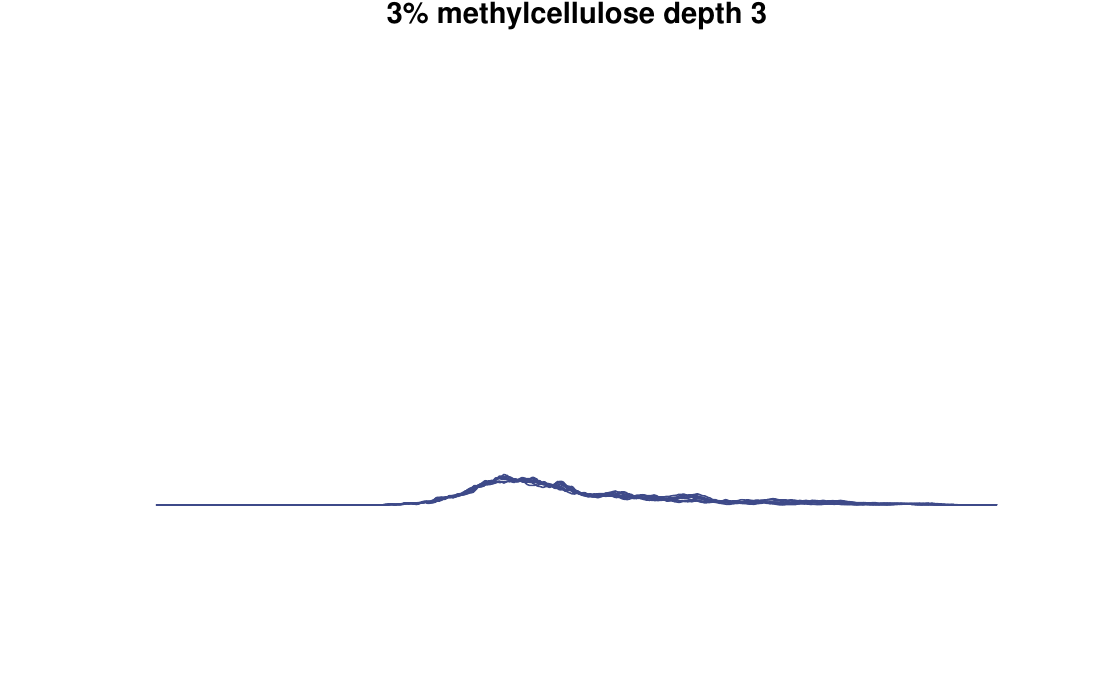}
    \end{tabular}
  \end{center}
  \caption[Average persistence landscapes for each environmental condition]{(A) The average landscapes for each environmental viscosity, distinguished by percentage of methylcellulose present. The averages here are taken over the average persistence landscape of videos (i.e. individuals) in a given viscosity class. (B) Standard deviations of each coordinate in the average persistence landscapes for each class. (C) The first landscapes of each sample landscape, organized by class. These concurrently plotted first landscapes show the variation in the samples for each class. (D-E) The second and third landscapes, respectively, for each sample according to its class. All plots share the same $x$-axis; the groups of plots in (A), (B), and (C-E) each have their own $y$-axis scale. }
  \label{fig:class-avg-PL-variation}
\end{figure}

Average persistence landscapes for each class are shown in \Cref{fig:class-avg-PL-variation} (A). The lower viscosity conditions allow for larger depth $1$ landscapes ($\lambda_1$) but have relatively few non-zero higher-depth landscapes, which means there are a smaller number of larger loops detected in the sliding window embeddings. This indicates that in lower-viscosity environments, \emph{C. elegans} exhibit behaviors of higher amplitude but either demonstrate fewer distinct behaviors or have much less variation between repetitions of behaviors. Conversely, the high-viscosity classes show many more cycles in the sliding window embeddings, with each cycle being small compared to the cycles found in the low-viscosity environments. These observations suggest that at high-viscosity, behaviors do not involve large changes in posture and are more varied. From observing the raw video data, it is apparent that in higher-viscosity environments \emph{C. elegans} can make smaller, tighter body bends, which is consistent with these results. 

We also observed that as viscosity increases, the support of the persistence landscapes stretches further to the right and cycles are born at lower radius values. 
The worms seemed to exhibit less varied behaviors in lower-viscosity environments, so perhaps in such environments they continued ``retracing their steps'' through the sliding window embedding space which resulted in more densely-sampled curves and thus homology classes formed at lower radii. 

\begin{figure}[H] 
  \begin{center}
    \includegraphics[trim = 0 40 0 0, clip, height=6.5cm]{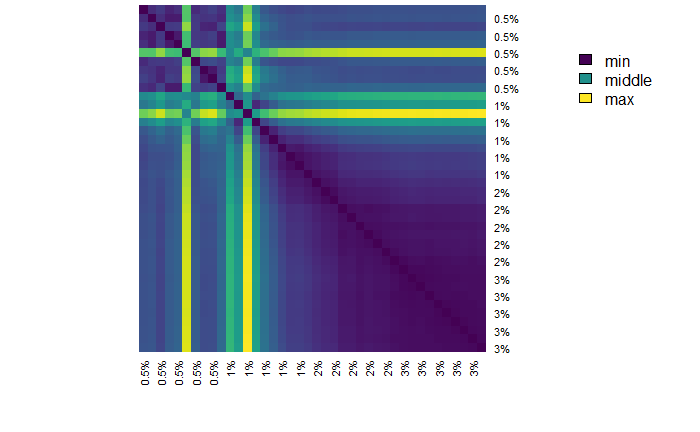}
  \end{center}
  \caption[Heatmap of sample distance]{Heatmap of the distance matrix of average persistence landscapes of samples.}
  \label{fig:heatmap}
\end{figure}

\begin{figure}[H] 
  \begin{center}
    \begin{tabular}{cc}
      \multicolumn{1}{l}{(A)} & \multicolumn{1}{l}{(B)}\\[1.0em]
      \includegraphics[trim = 20 30 20 30, clip, width=\halfwidth]{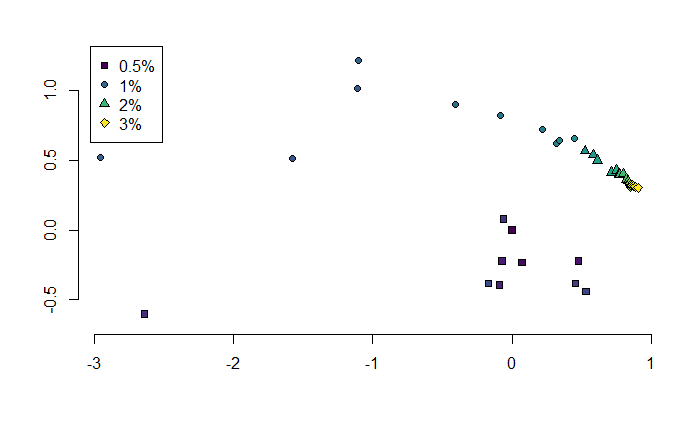}
      & \includegraphics[trim = 50 70 50 70, clip, width=\halfwidth]{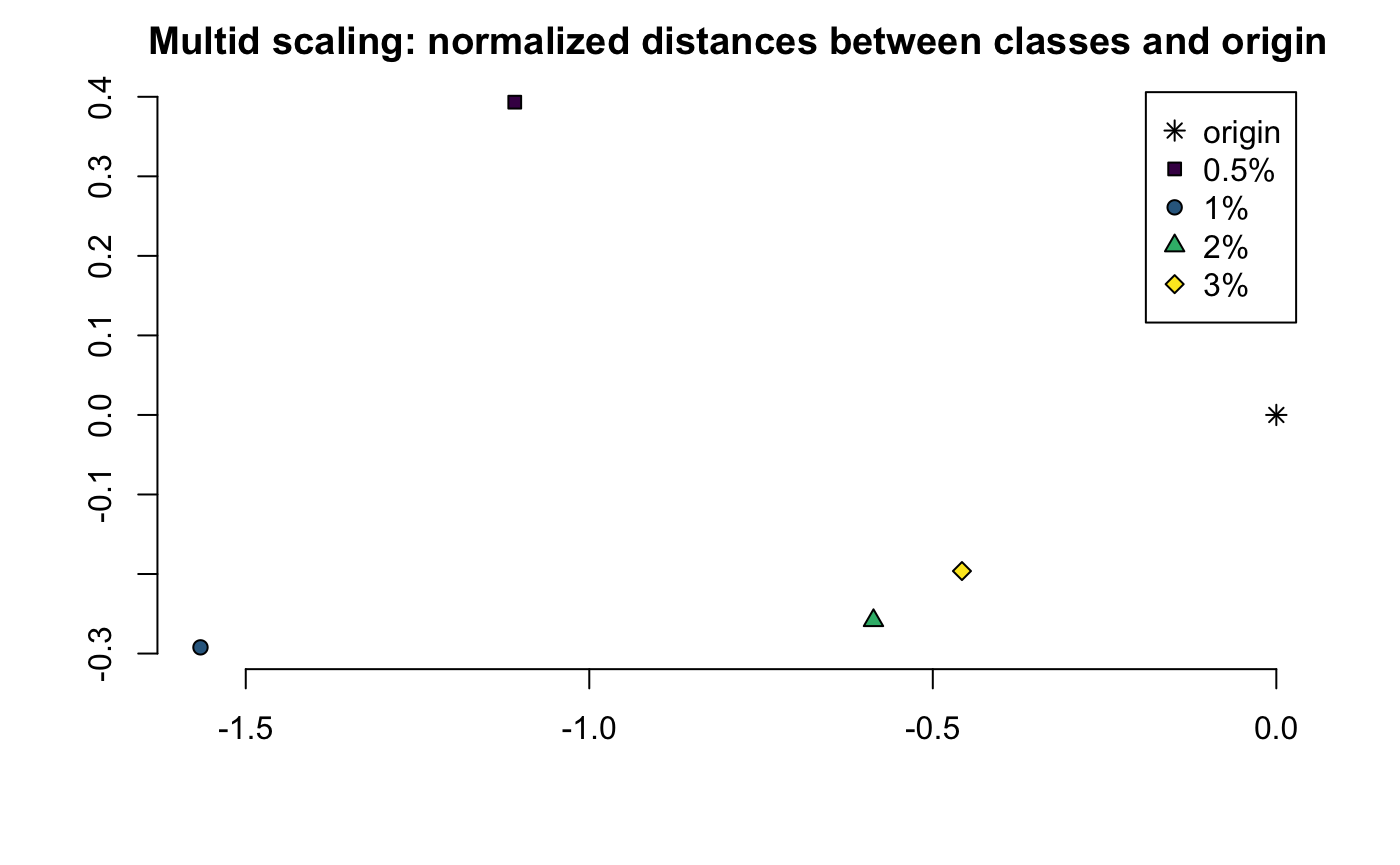}
    \end{tabular}
  \end{center}
  \caption[Multidimensional scaling of pairwise distances between samples]{(A) Multidimensional scaling of average persistence landscape of each sample. (B) Multidimensional scaling of average persistence landscapes of the classes and the origin.}
  \label{fig:multid-scaling}
\end{figure}

The pairwise distances between landscapes for each sample are visualized in \Cref{fig:heatmap}. 
The normalized pairwise distances between the average persistence landscapes for each class are shown in \Cref{table:distances}. We include the origin ---  the zero persistence landscape, \ie the $0$ vector --- in these distance computations to complete the normalization. Normalization is such that the average distances between each class and the origin is $1$. Multidimensional scaling on these distances visualizes the similarities between samples and classes, respectively, and are shown in \Cref{fig:multid-scaling}. From the raw distances and the multidimensional scaling of the distances, we can see that the high-viscosity classes ($2\%$ and $3\%$ methylcellulose) are closest together and that this pair, the $0.5\%$ class, and the $1\%$ class are roughly equidistant from one another. 

\begin{table}[H] 
  \begin{center}
    \begin{tabular}{c | c c c c}
      & $0.5\%$& $1\%$ & $2\%$ & $3\%$\\
      \hline
      origin & $1.1873106$ & $1.5934976$ & $0.6777196$ & $0.5414722$\\
      $0.5\%$ & & $0.8330355$ & $0.8380992$ & $0.8809327$\\
      $1\%$ & & & $0.9972502$ & $1.1255496$\\
      $2\%$ & & & & $0.1758501$
    \end{tabular}
  \end{center}
  \caption[Pairwise distances between classes]{Normalized pairwise distances between average persistence landscapes of each class, where distance is Euclidean distance between vectors in $\RR^{255969}$ and the normalization is such that the average distance to the origin is $1$. }
  \label{table:distances}
\end{table}

\begin{figure}[H] 
  \begin{center}
    \begin{tabular}{m{0.7\columnwidth} m{0.3\columnwidth}}
      \multicolumn{1}{l}{(A)} & \multicolumn{1}{l}{(B)}\\[1.0em]
      \includegraphics[trim = 50 70 50 64,clip, height=6cm]{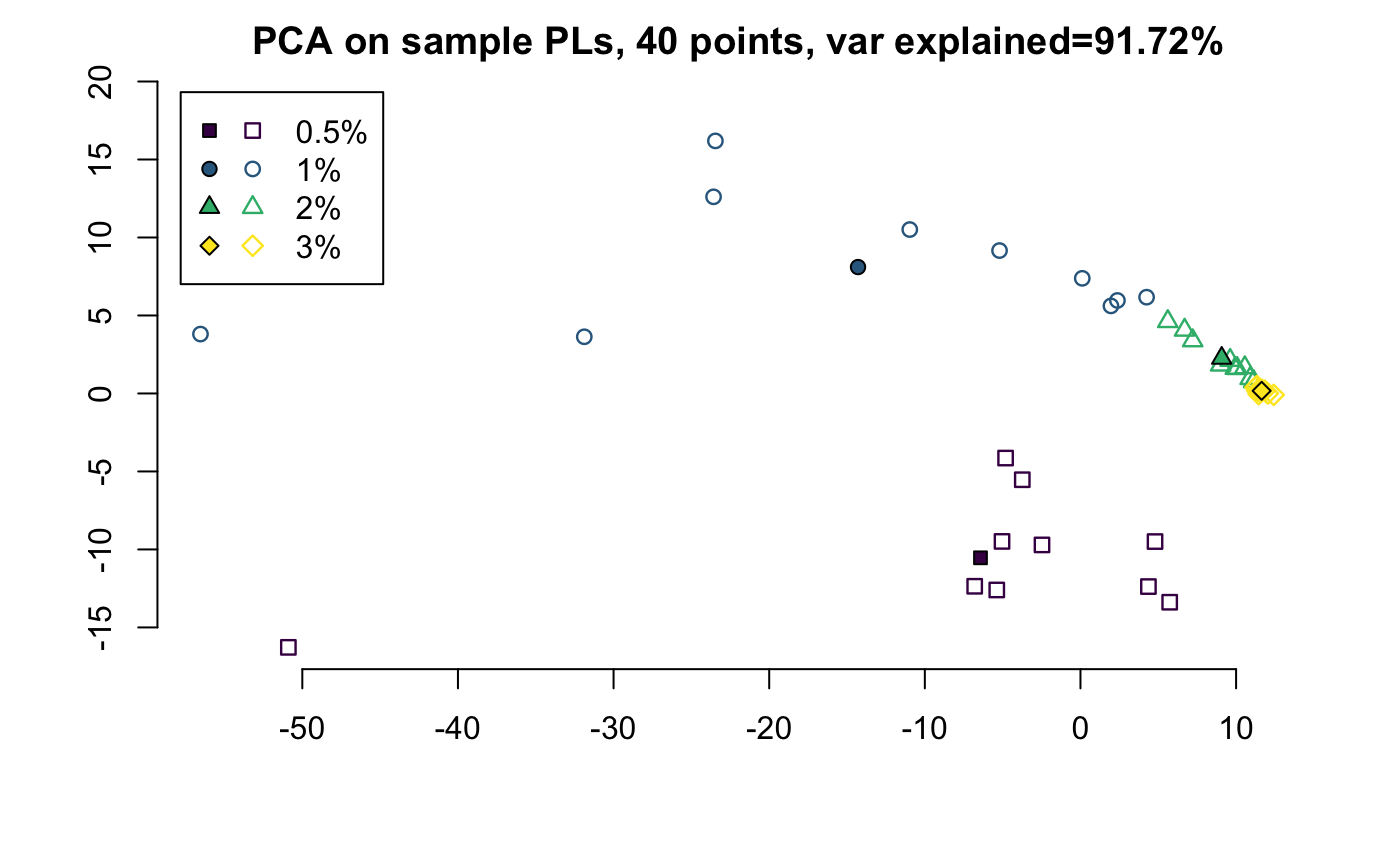}
      & \begin{tabular}{c}
        \raisebox{3em}{\includegraphics[trim = 165 205 50 300,clip, width=0.3\columnwidth]{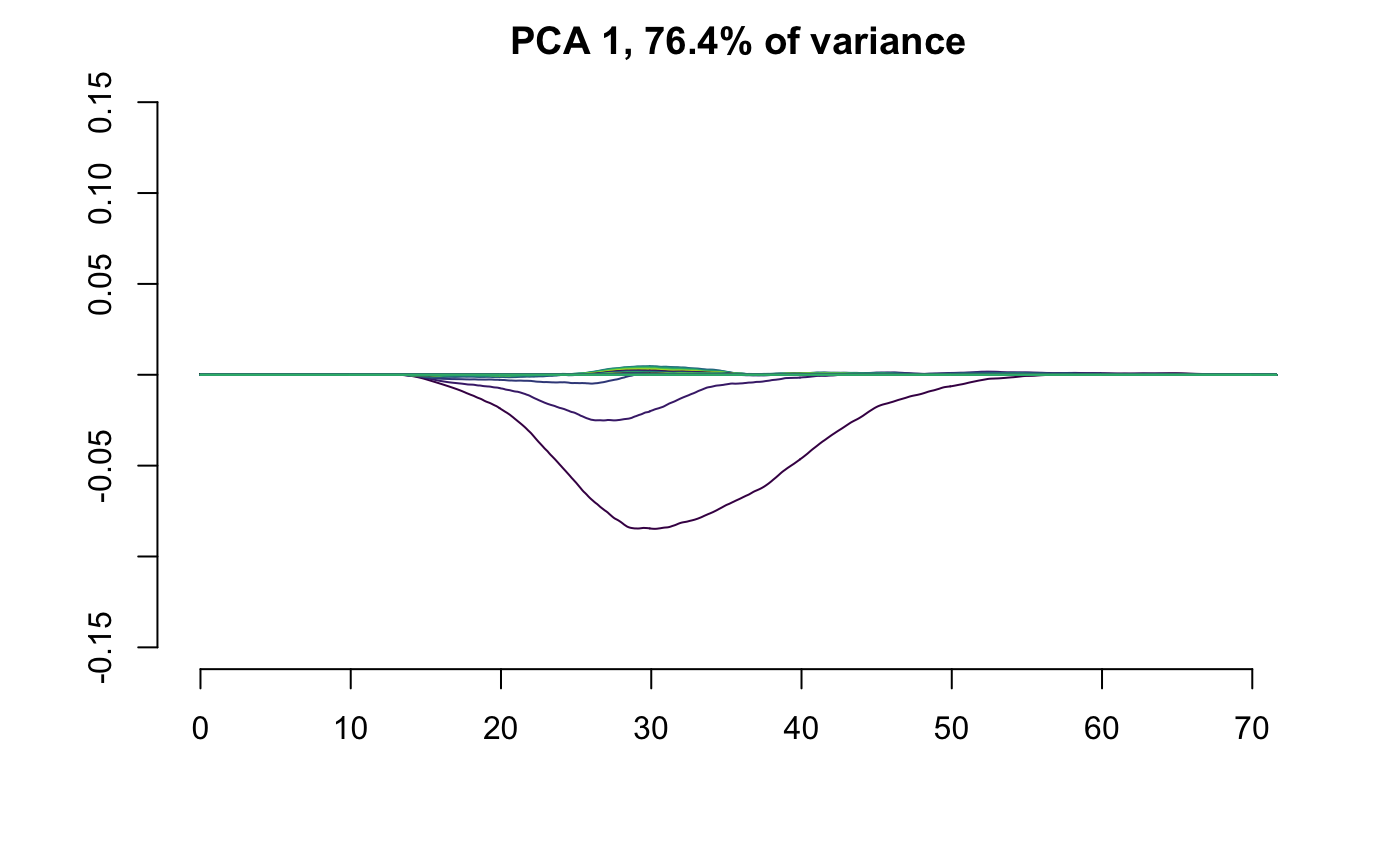}}\\[-2.5em]
        \raisebox{3em}{\includegraphics[trim = 165 205 50 250,clip,width=0.3\columnwidth]{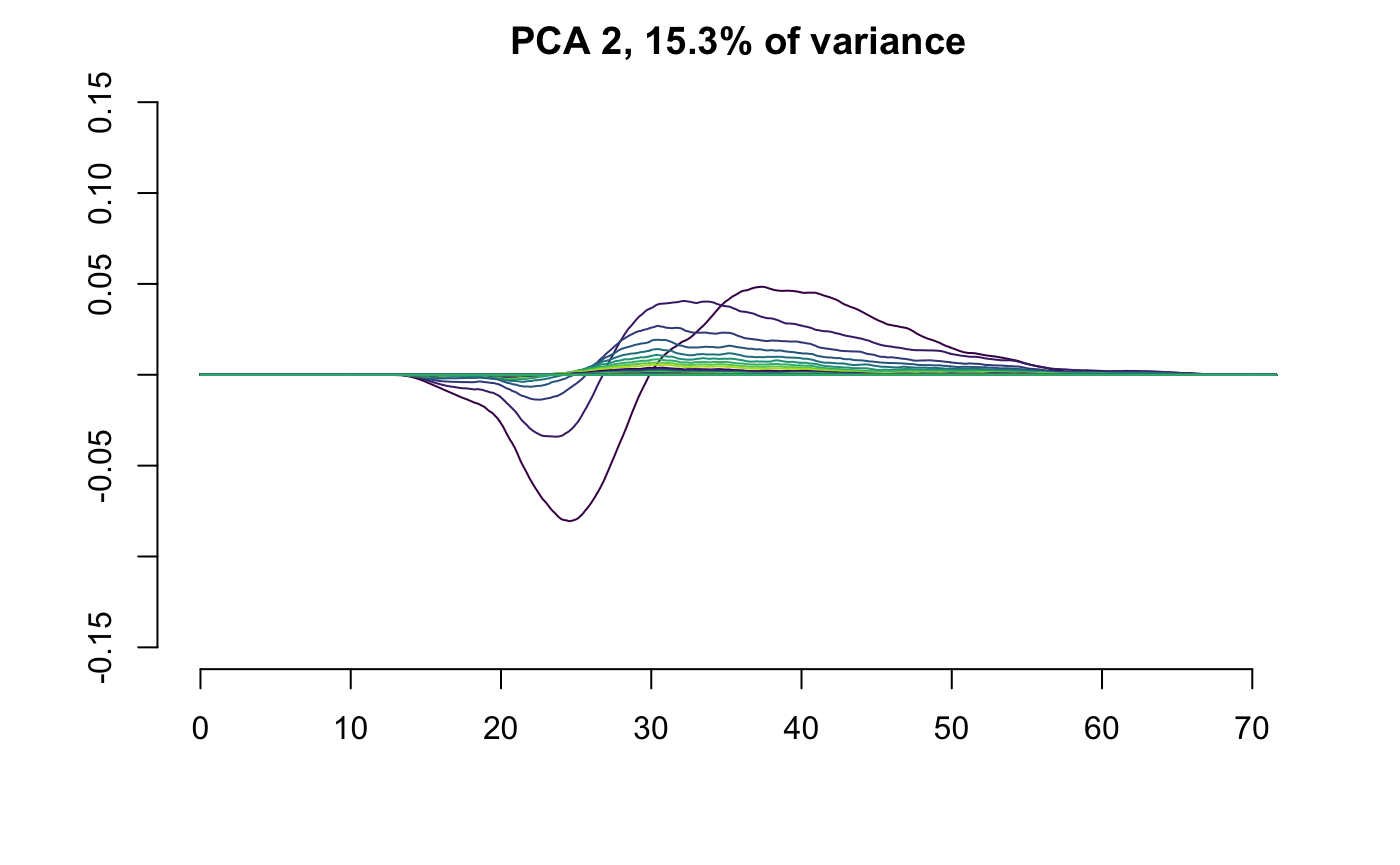}}
      \end{tabular}
    \end{tabular}
  \end{center}
  \vspace{-2em}
  \caption[PCA on all samples]{(A) Projection of average persistence landscapes onto $2$ principal components, with average persistence landscapes of videos given by outlined symbols and average persistence landscapes of classes given by solid symbols. (B) The first two PCA eigenvectors.}
  \label{fig:full-pca}
\end{figure}

Principal component analysis on the average persistence landscapes for each sample gives the graphs in \Cref{fig:full-pca}. In (A), the projections of the average persistence landscapes of the samples onto the first two principal components are given by hollow symbols and projections of the average persistence landscapes of the classes are given by solid symbols. Here we see results similar to those from the multidimensional scaling in \Cref{fig:multid-scaling}: the low-viscosity class landscapes are far from each other and the high-viscosity classes, while the high-viscosity class landscapes are quite close. We can also see some of the variance within classes. The low-viscosity classes have much more variability than the high-viscosity classes, with the highest-viscosity class, $3\%$ methylcellulose, having very little variation in these first two principal components. 

These conclusions about variation in each of the classes are supported by the standard deviations of each coordinate in the average (discrete) persistence landscapes of  each class. In \Cref{fig:class-avg-PL-variation}(B) the standard deviations of each coordinate are graphed as sequences of functions so that the standard deviations can be easily matched up with their corresponding locations on the average persistence landscapes. We conclude that there is little variation in the $3\%$ class, slightly more in the $2\%$ class, and much more in the $0.5\%$ and $1\%$ classes. The $1\%$ class showed more variation in higher-depth landscapes than the $0.5\%$ class, suggesting that \emph{C. elegans} can produce slightly more complex behaviors in a slightly higher viscosity environments. The variances of the $0.5\%$ and $1\%$ classes also exhibit a distinct pattern; the $0.5\%$ samples varied more towards the lower radius parameters (the left side of the graph), whereas the $1\%$ samples varied more towards higher (more to the right) radius parameters. 

\begin{figure}[H] 
  \def\po{\hspace{4.5em}}
  \begin{center}
    \begin{tabular}{cccc}
      $0.5\%$  & $1\%$  
      & $2\%$  & $3\%$ \\[0.5em]
      \multicolumn{1}{l}{\hspace{-1em}(A)}\\[0.5em]
      \includegraphics[width=\quarterwidth]{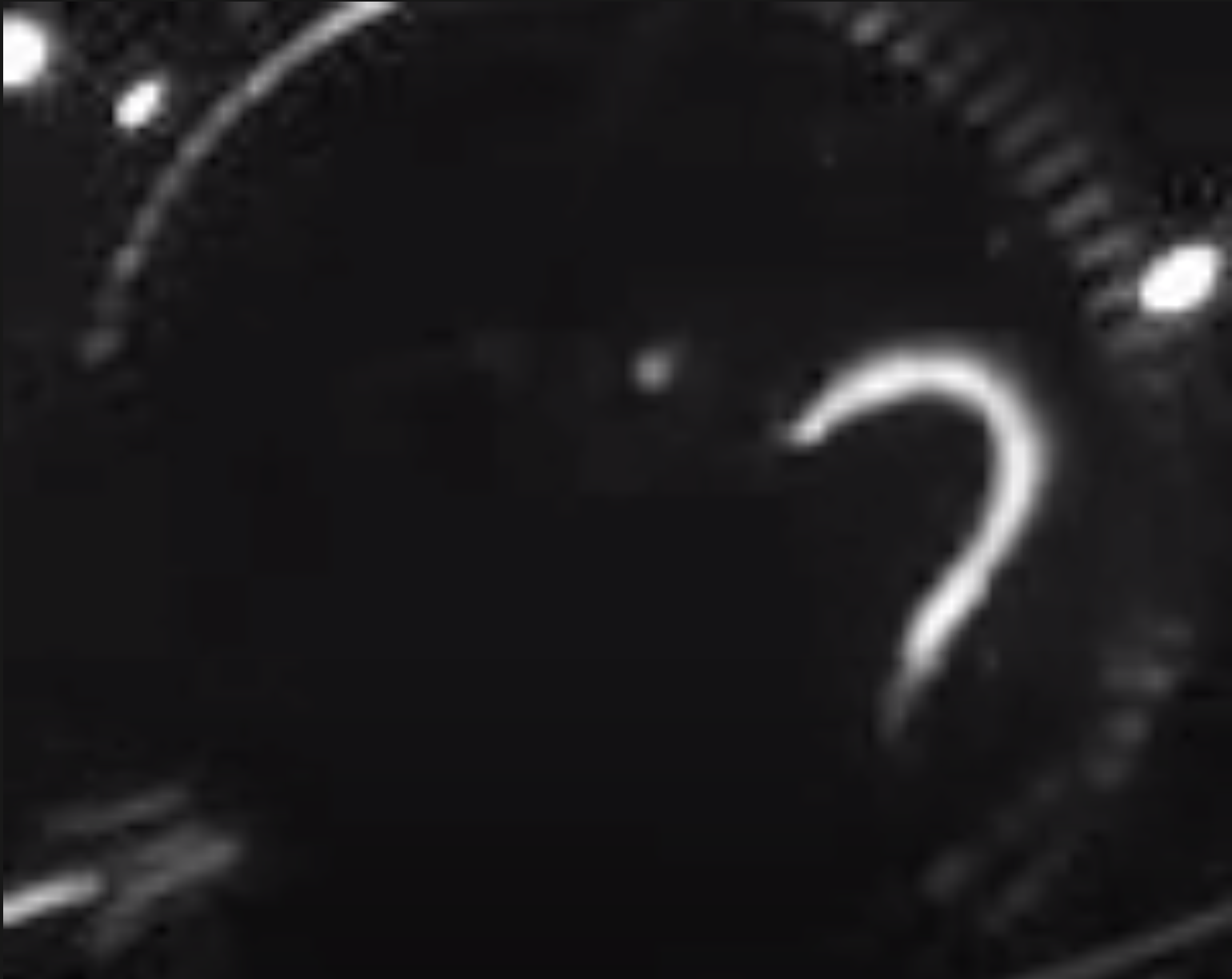}
      &\includegraphics[width=\quarterwidth]{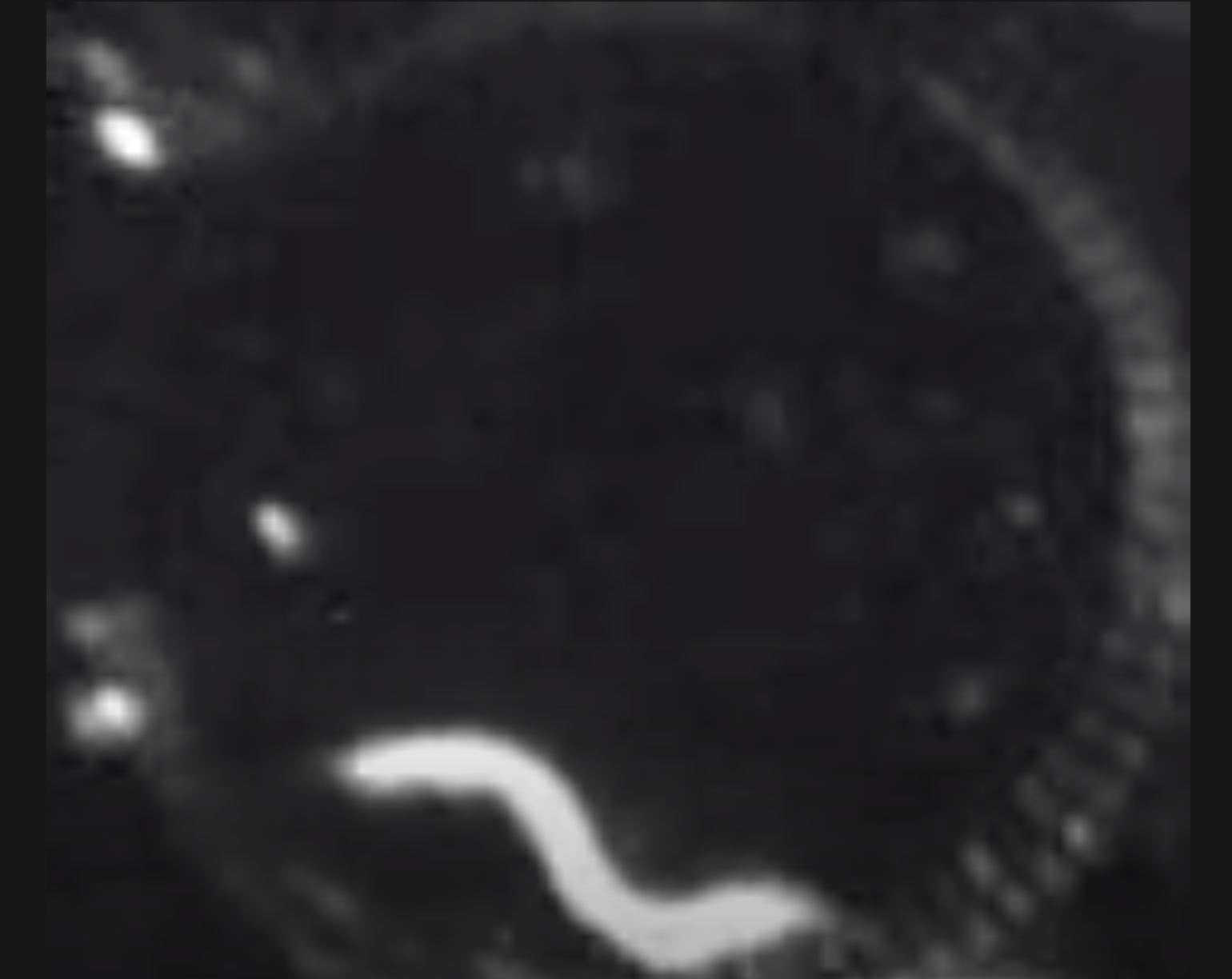}
      &\includegraphics[width=\quarterwidth]{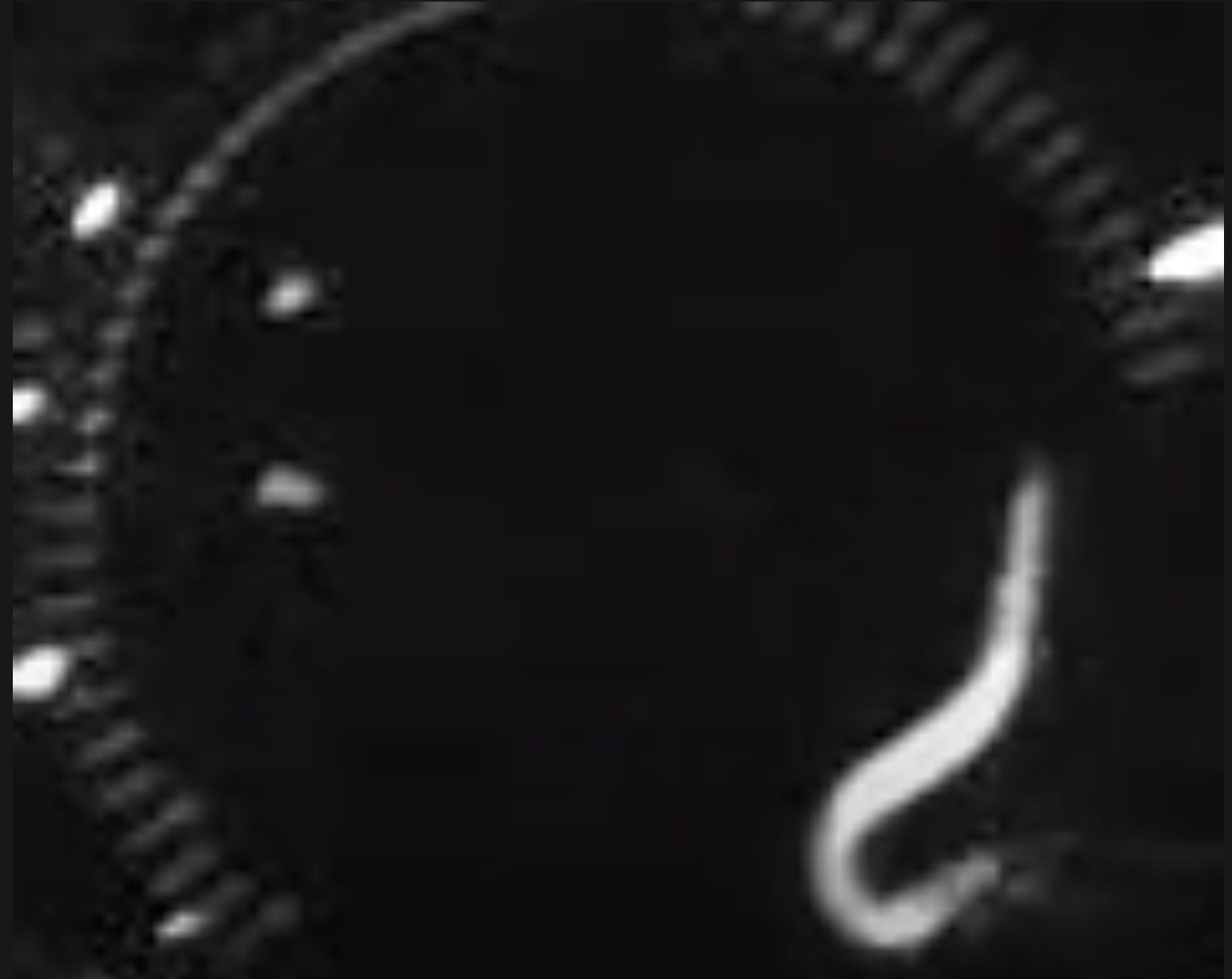}
      &\includegraphics[width=\quarterwidth]{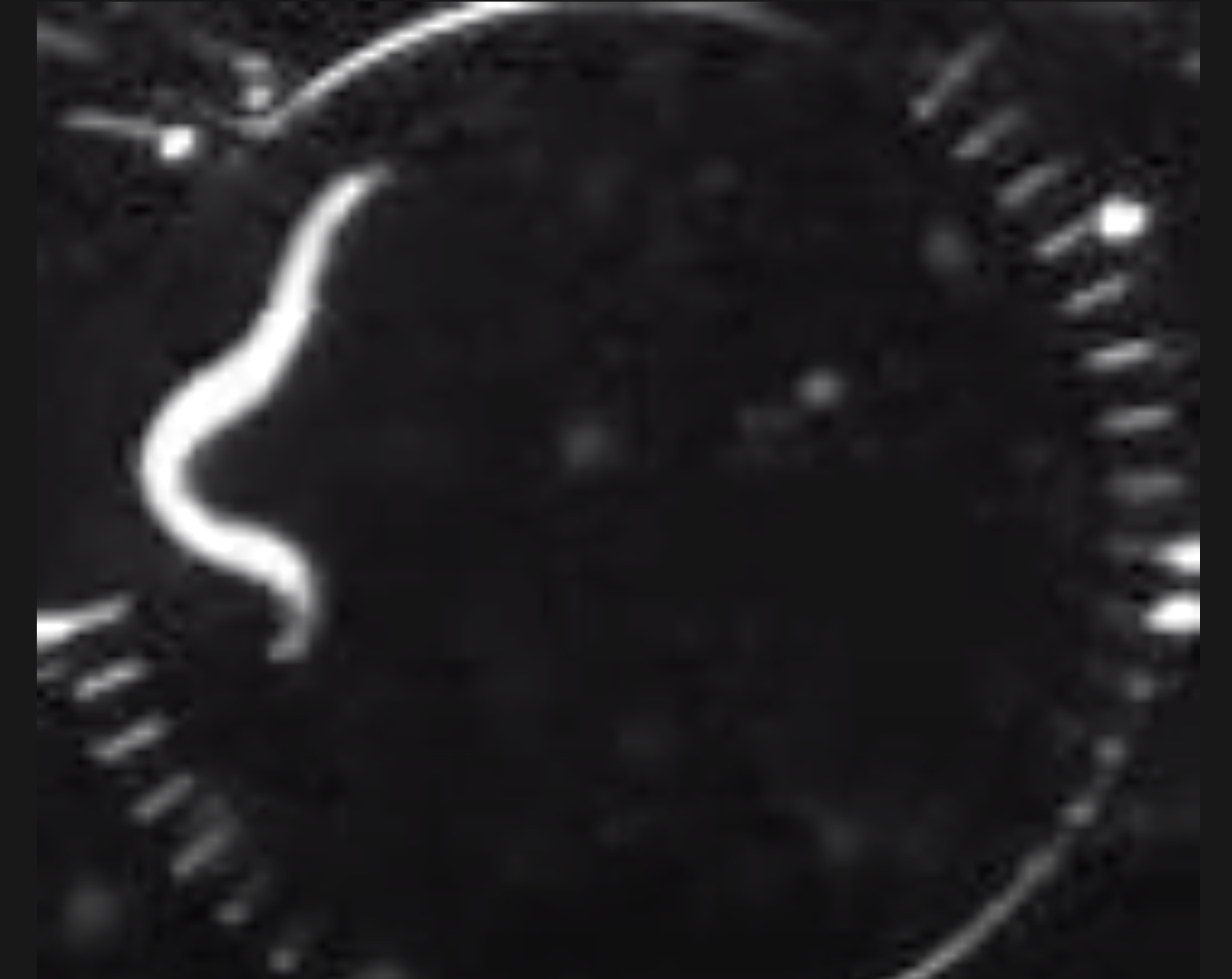}
      \\[1.5em]
      \multicolumn{1}{l}{\hspace{-1em}(B)}\\[0.5em]
      \trimmedwithaxesgraphic[width=\quarterwidth]{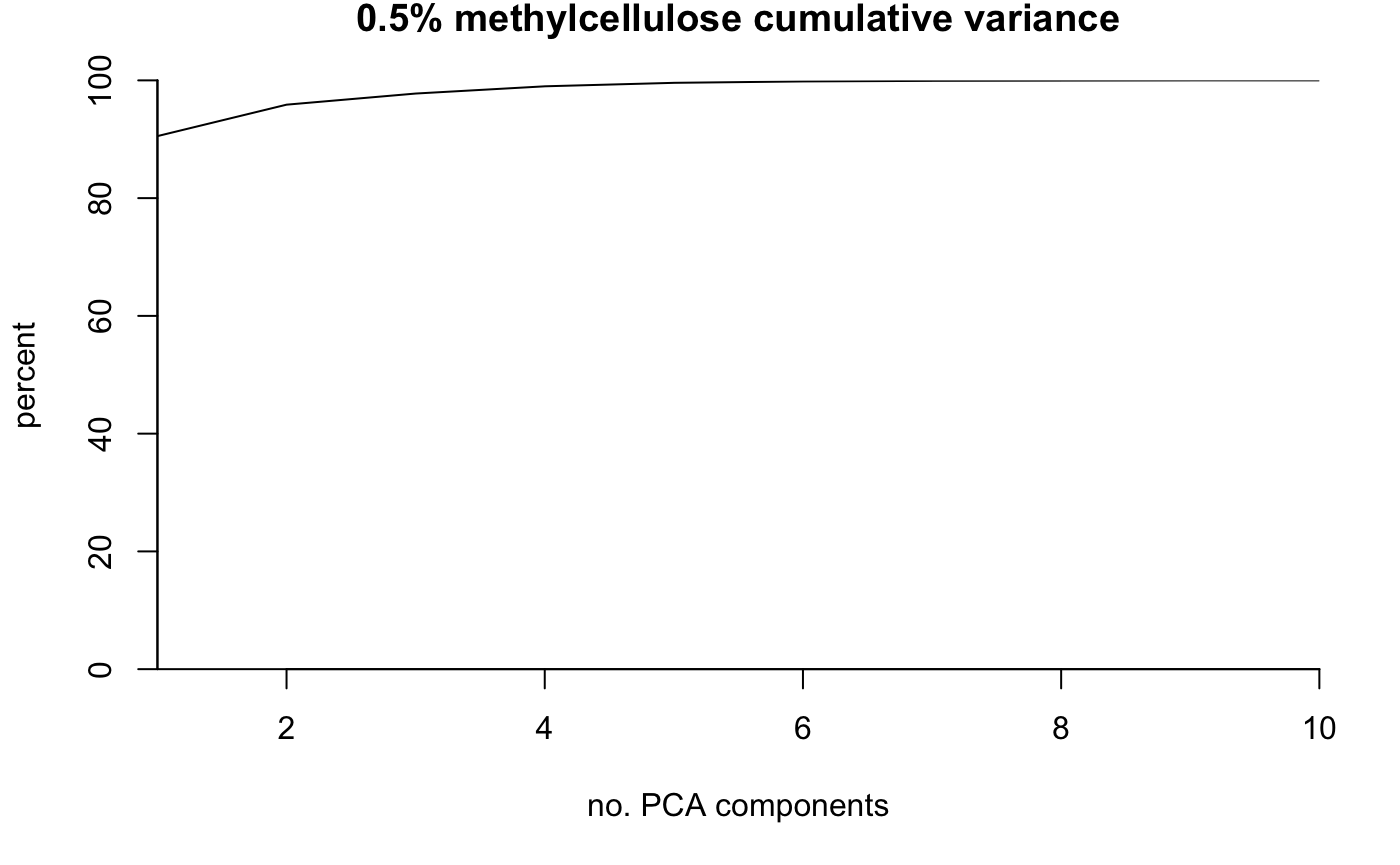}
      &\trimmedwithaxesgraphic[width=\quarterwidth]{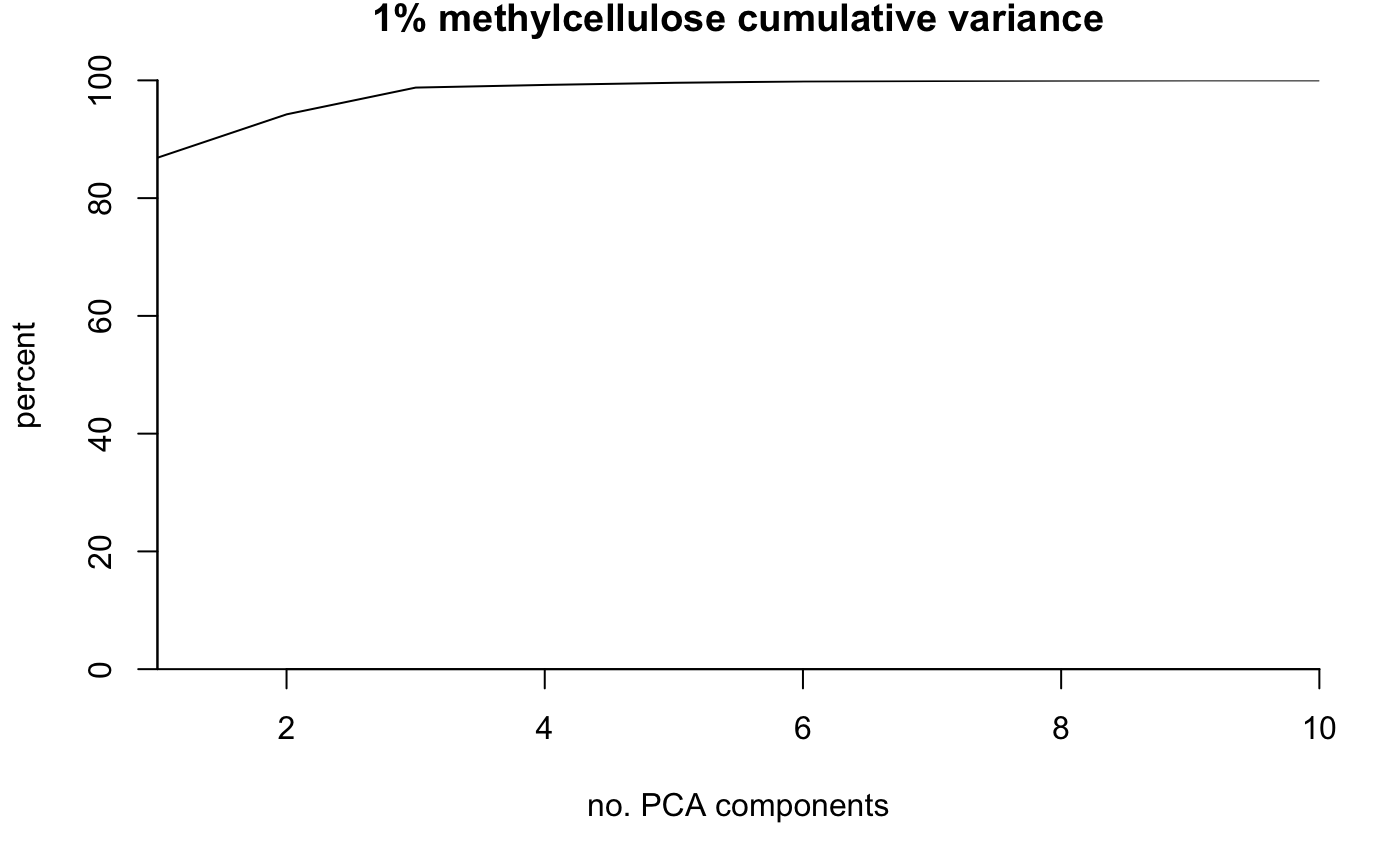}
      &\trimmedwithaxesgraphic[width=\quarterwidth]{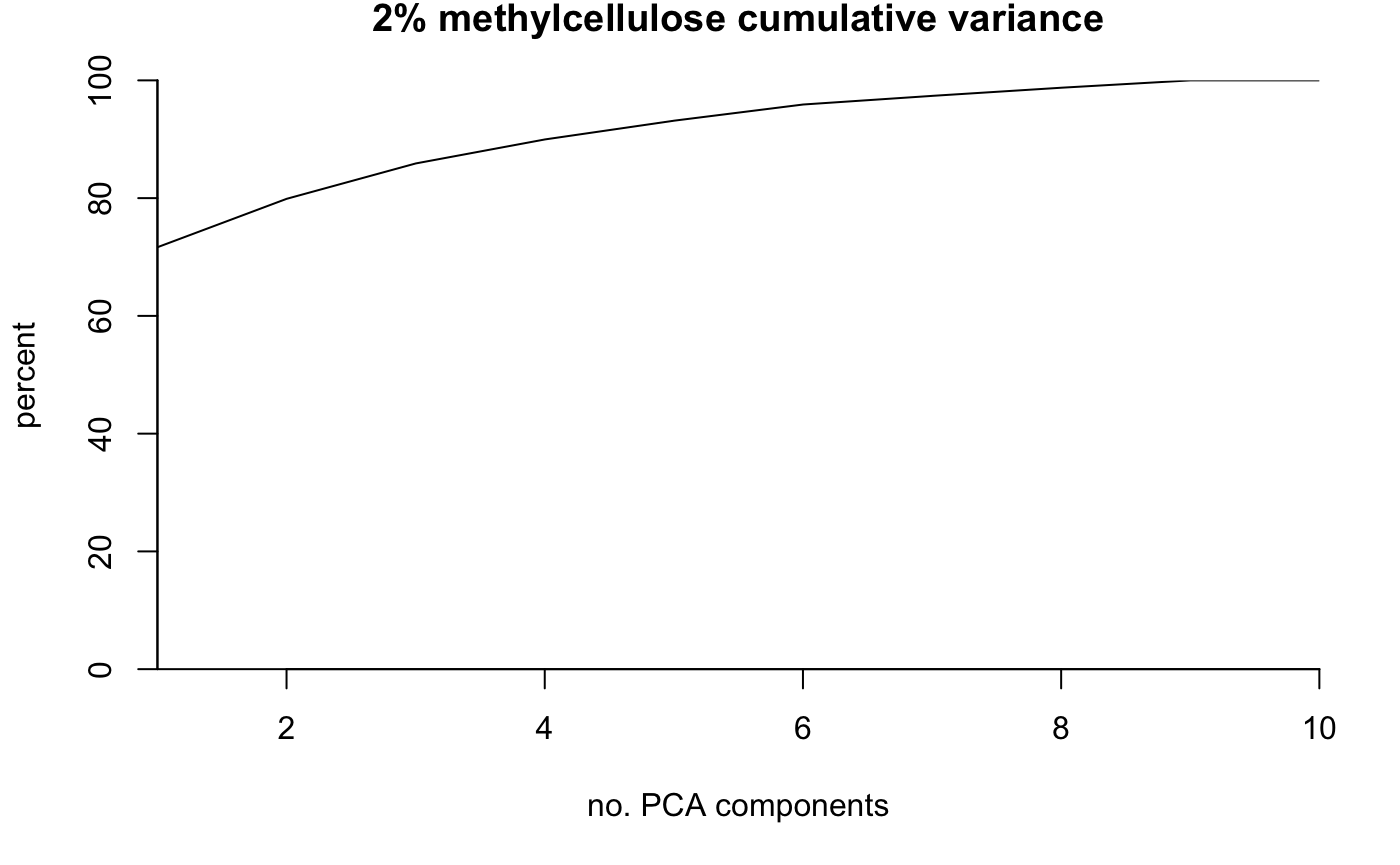}
      &\trimmedwithaxesgraphic[width=\quarterwidth]{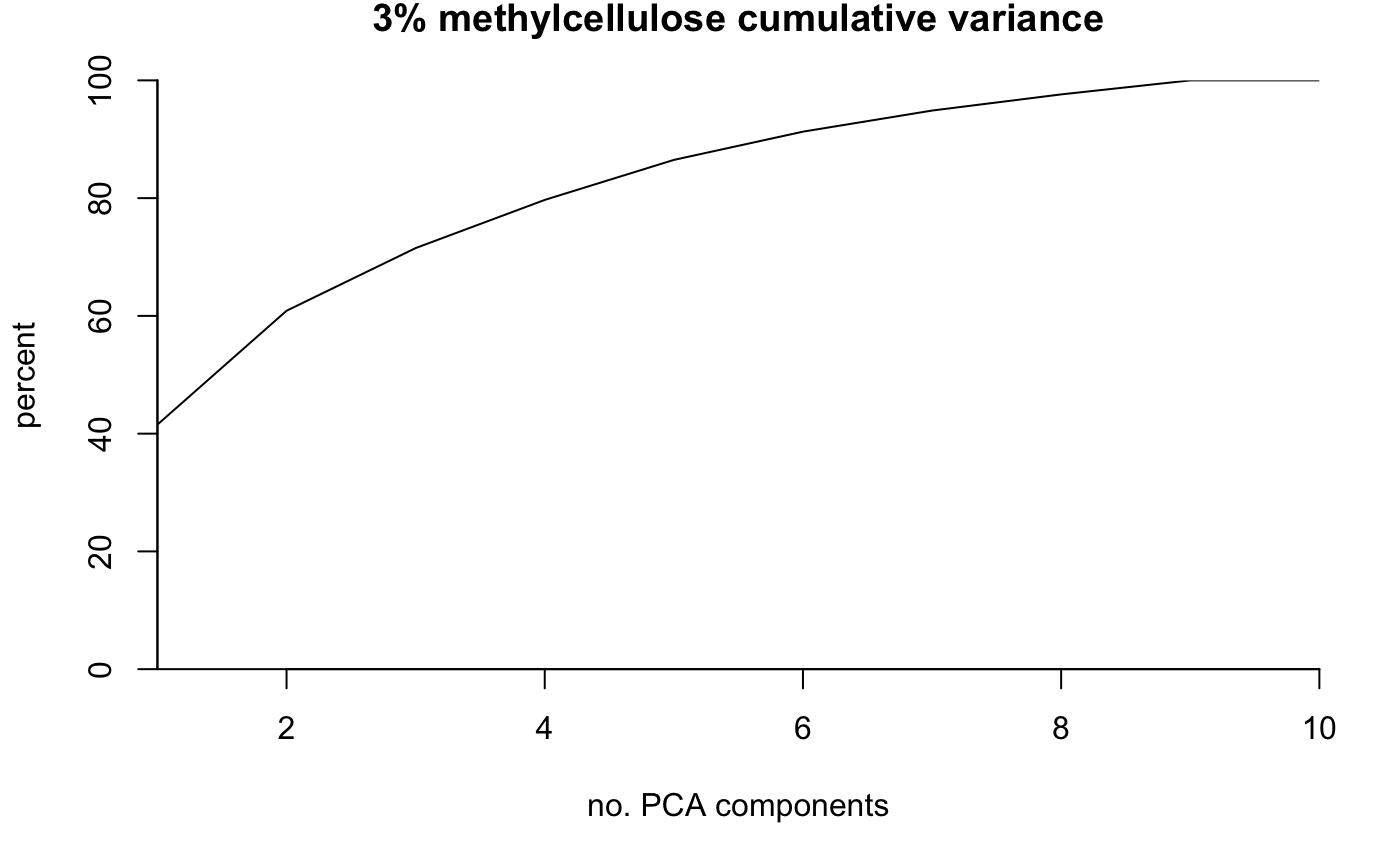}
      \\[1.5em]
      \multicolumn{1}{l}{\hspace{-1em}(C)}\\
      \po $90.5\%$ & \po $86.9\%$ & \po $71.7\%$ &\po $41.5\%$\\
      \trimmedgraphic[width=\quarterwidth]{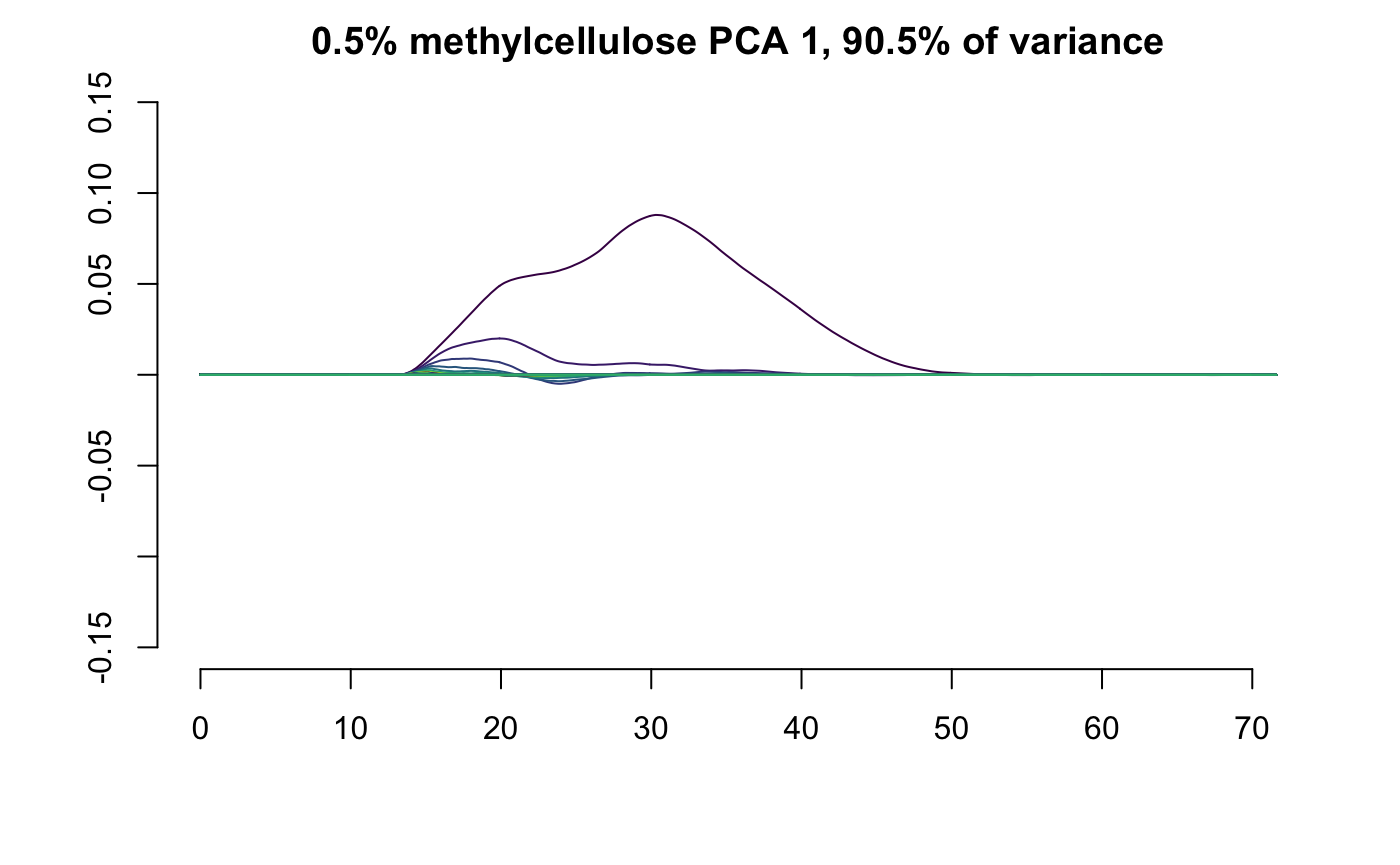}
      &\trimmedgraphic[width=\quarterwidth]{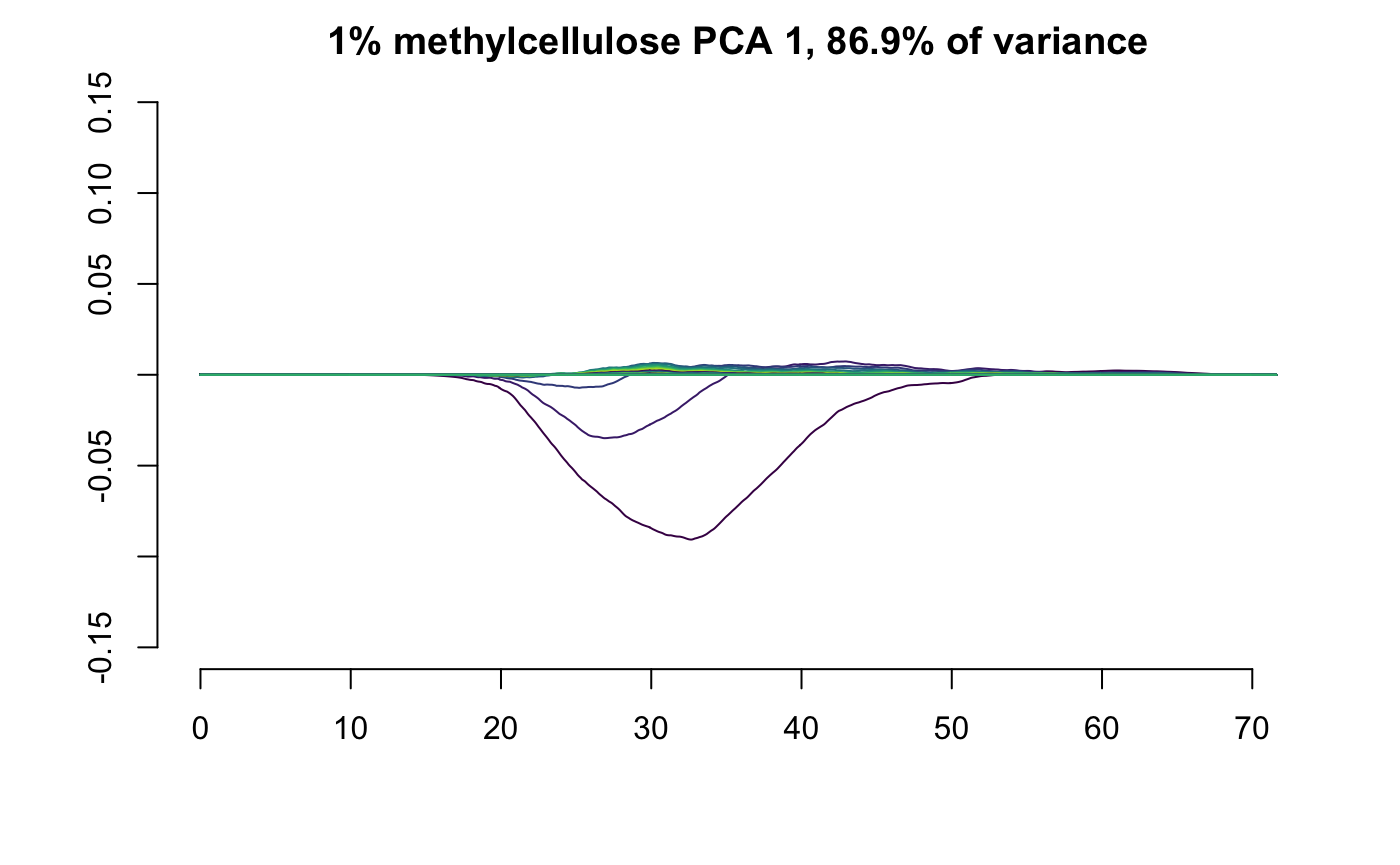}
      &\trimmedgraphic[width=\quarterwidth]{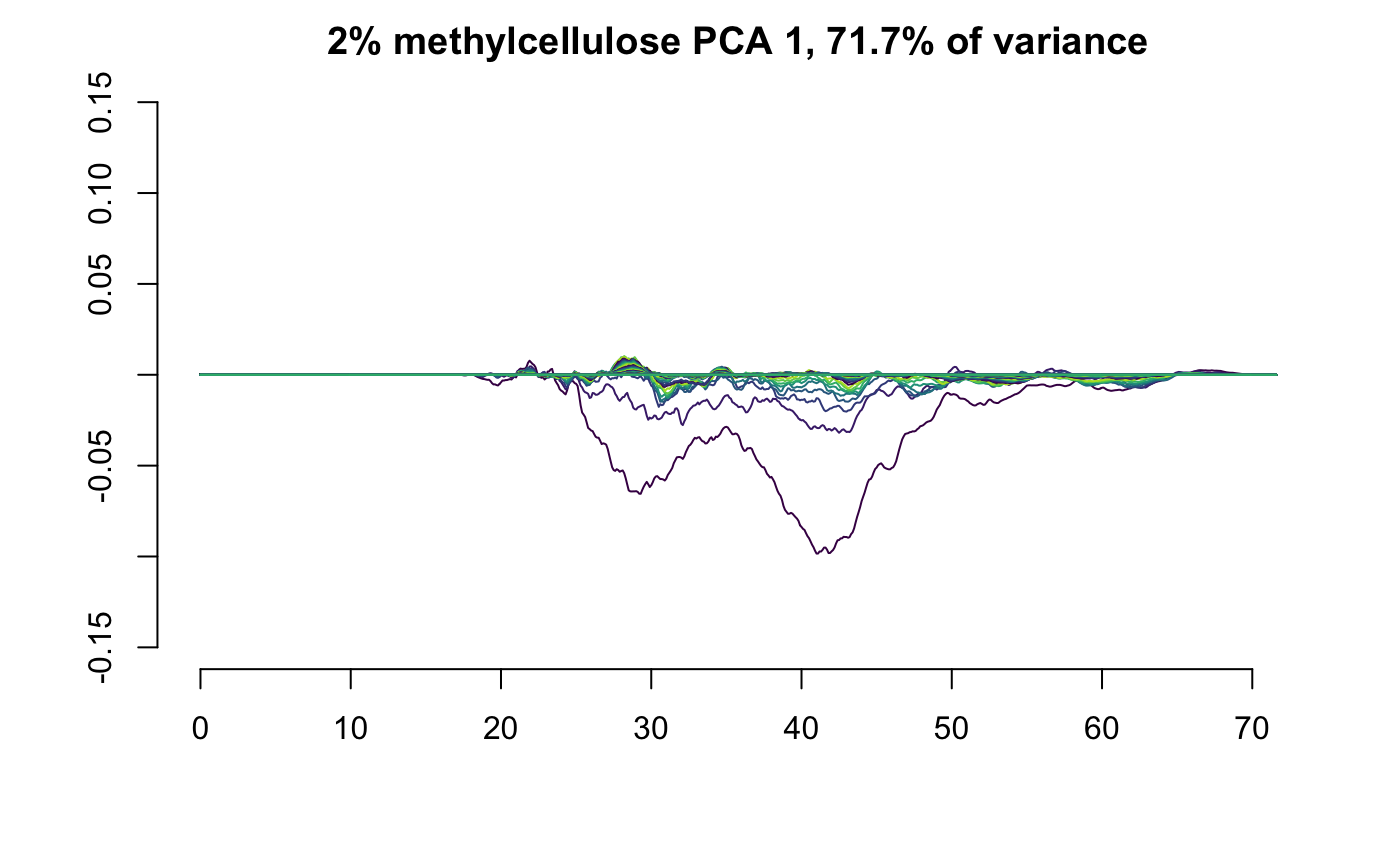}
      &\trimmedgraphic[width=\quarterwidth]{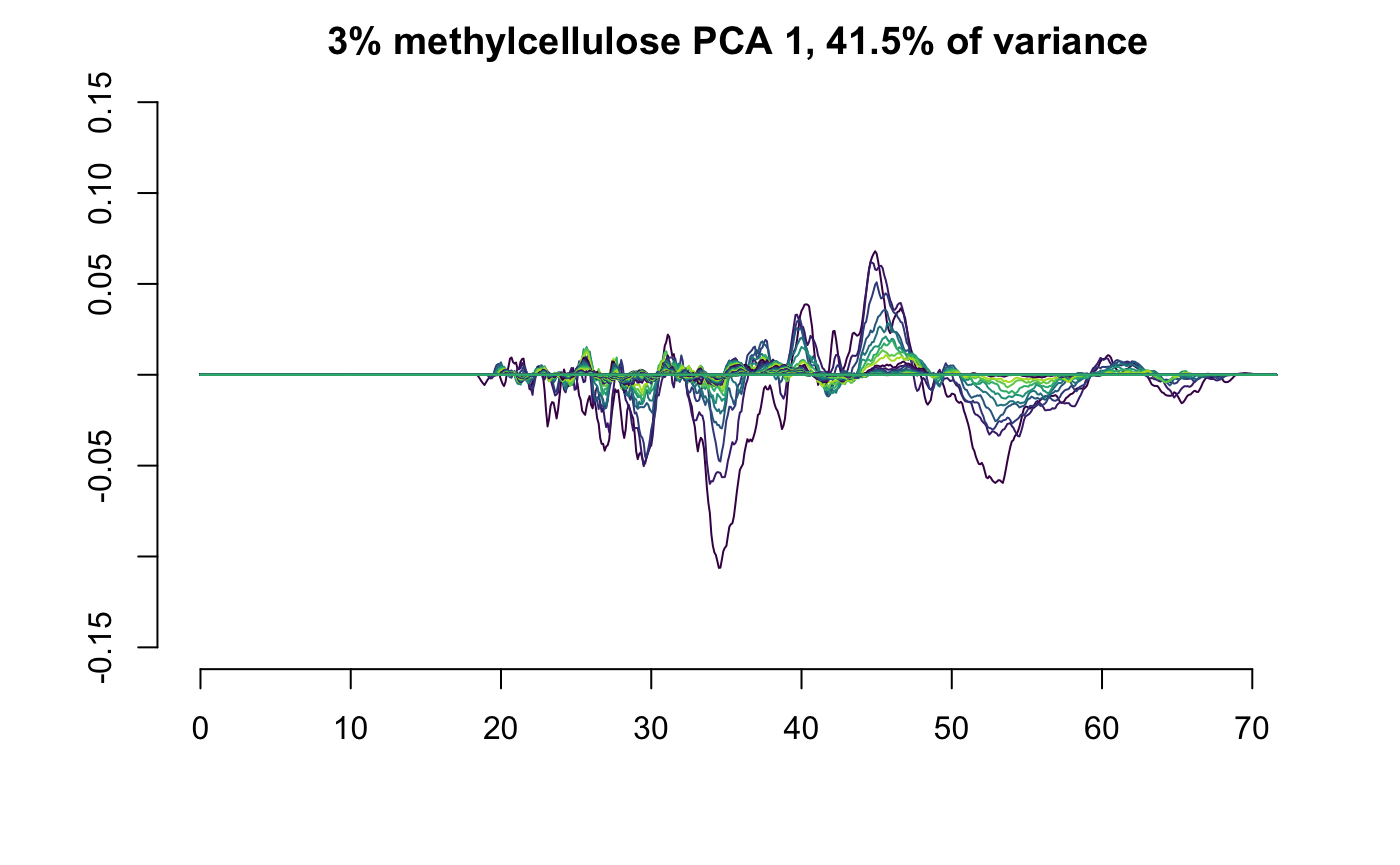}
      \\
      \multicolumn{1}{l}{\hspace{-1em}(D)}\\
      \po $5.3\%$ & \po $7.4\%$ & \po $8.2\%$ & \po $19.3\%$\\
      \trimmedgraphic[width=\quarterwidth]{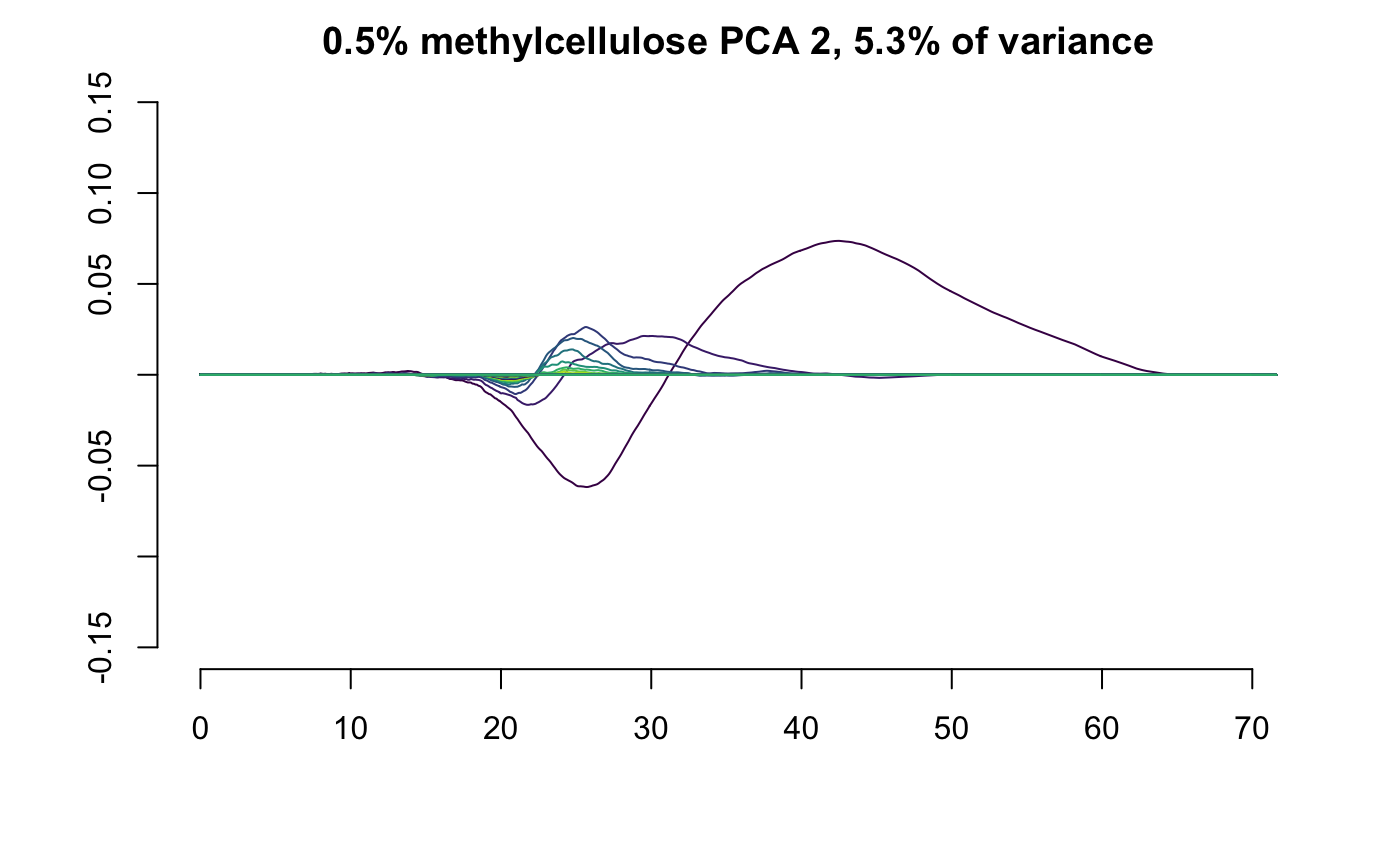}
      &\trimmedgraphic[width=\quarterwidth]{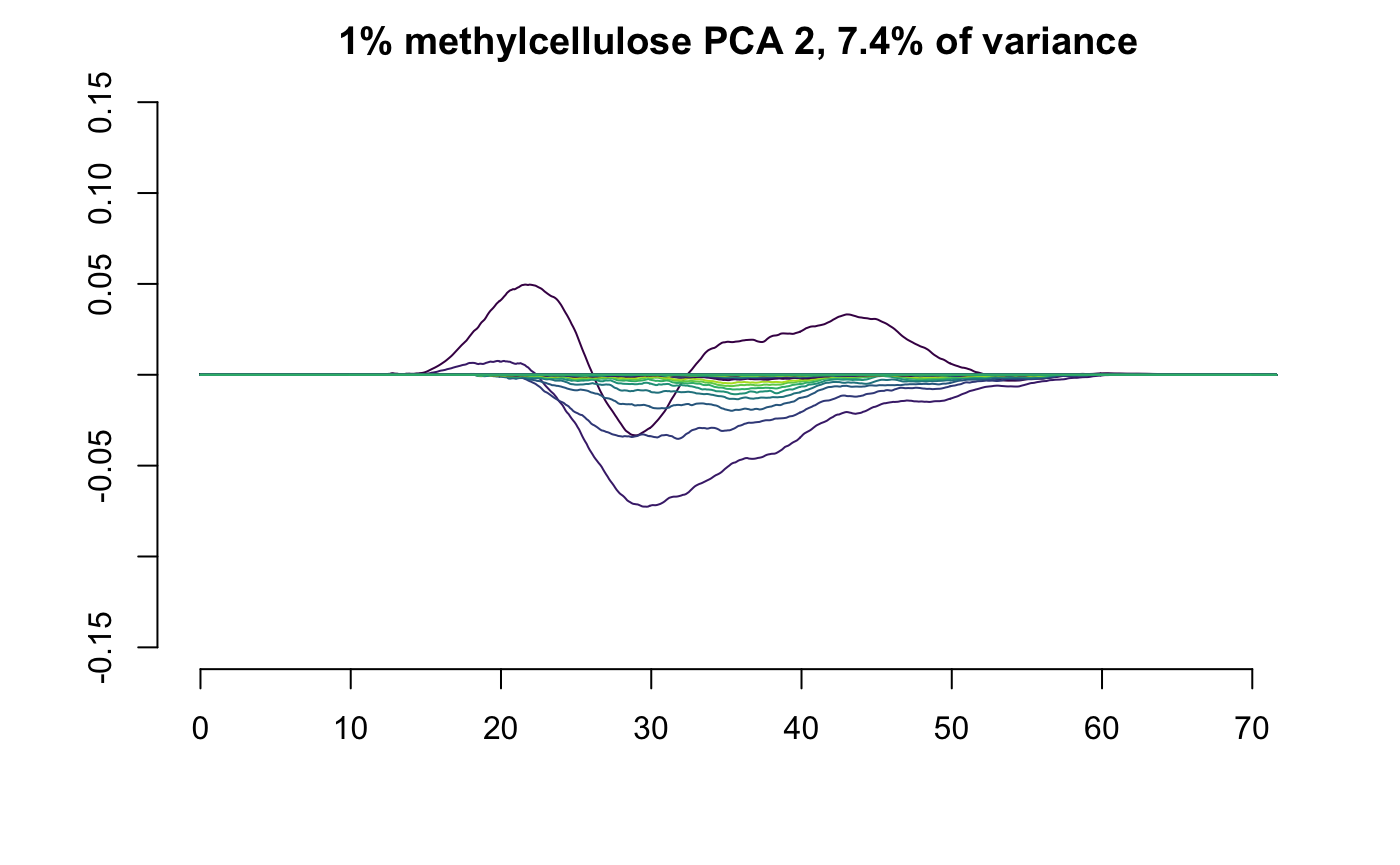}
      &\trimmedgraphic[width=\quarterwidth]{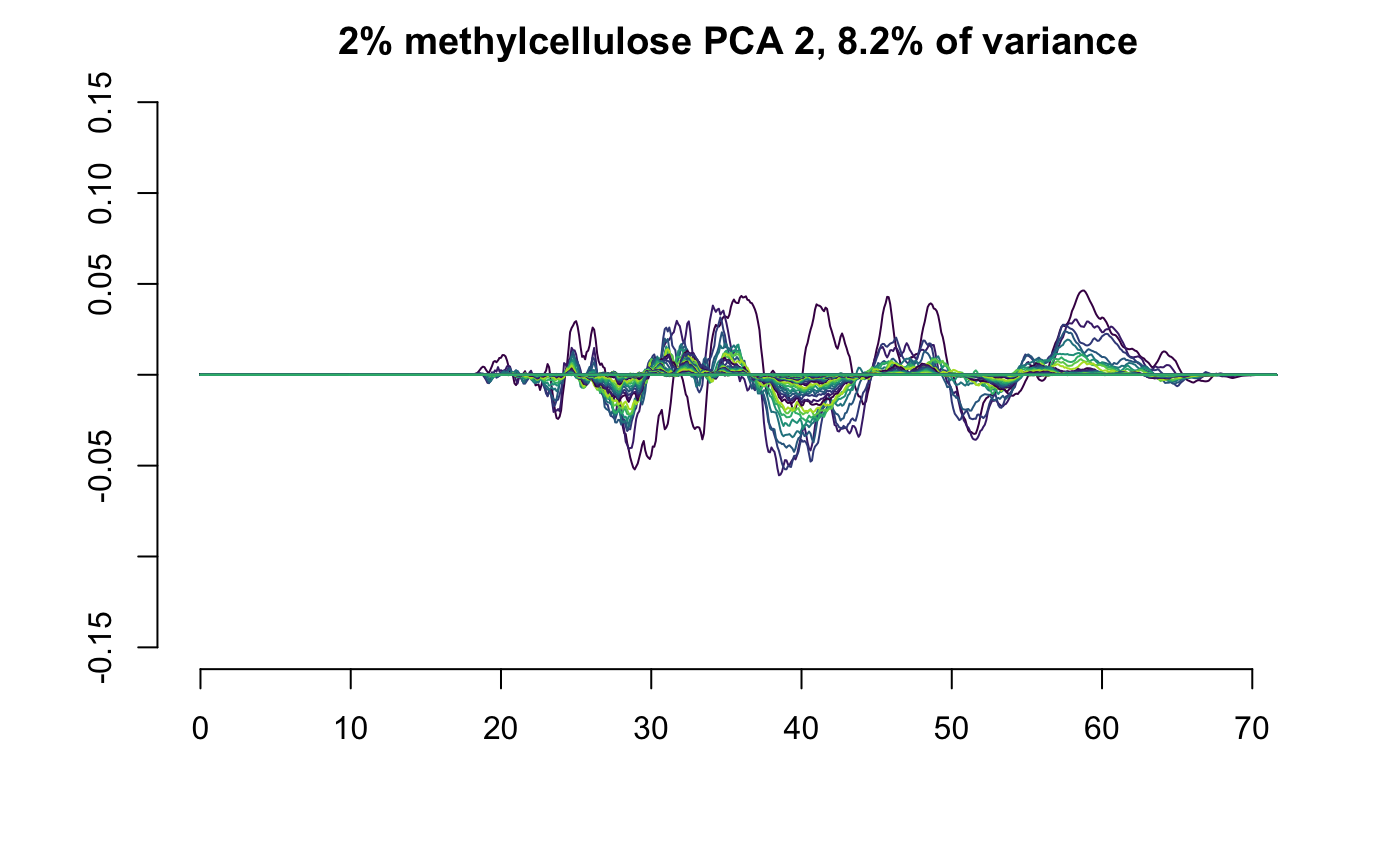}
      &\trimmedgraphic[width=\quarterwidth]{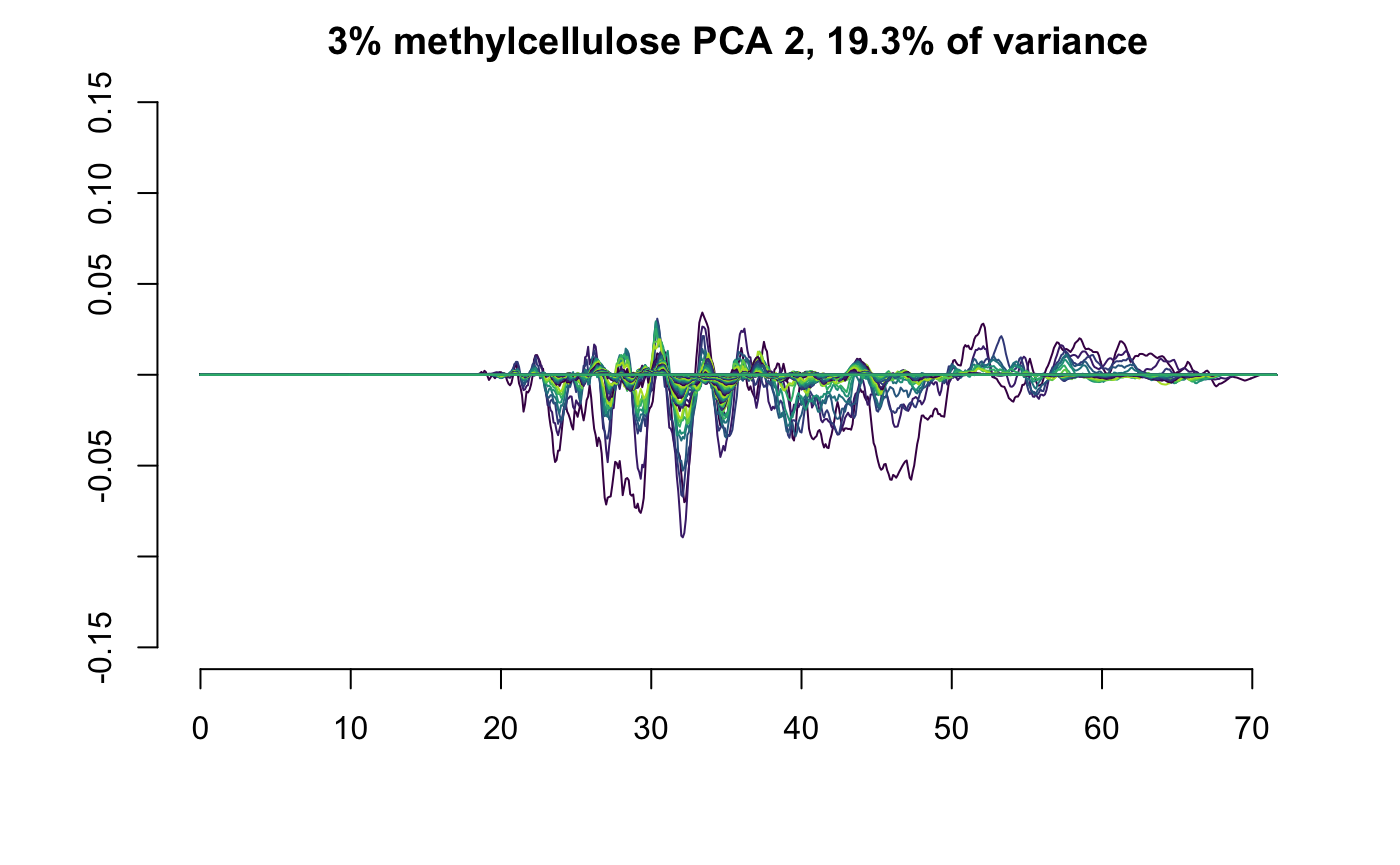}
      \\
      \multicolumn{1}{l}{\hspace{-1em}(E)}\\
      \po $1.9\%$ & \po $4.5\%$ & \po $6.0\%$ & \po $10.6\%$\\
      \trimmedgraphic[width=\quarterwidth]{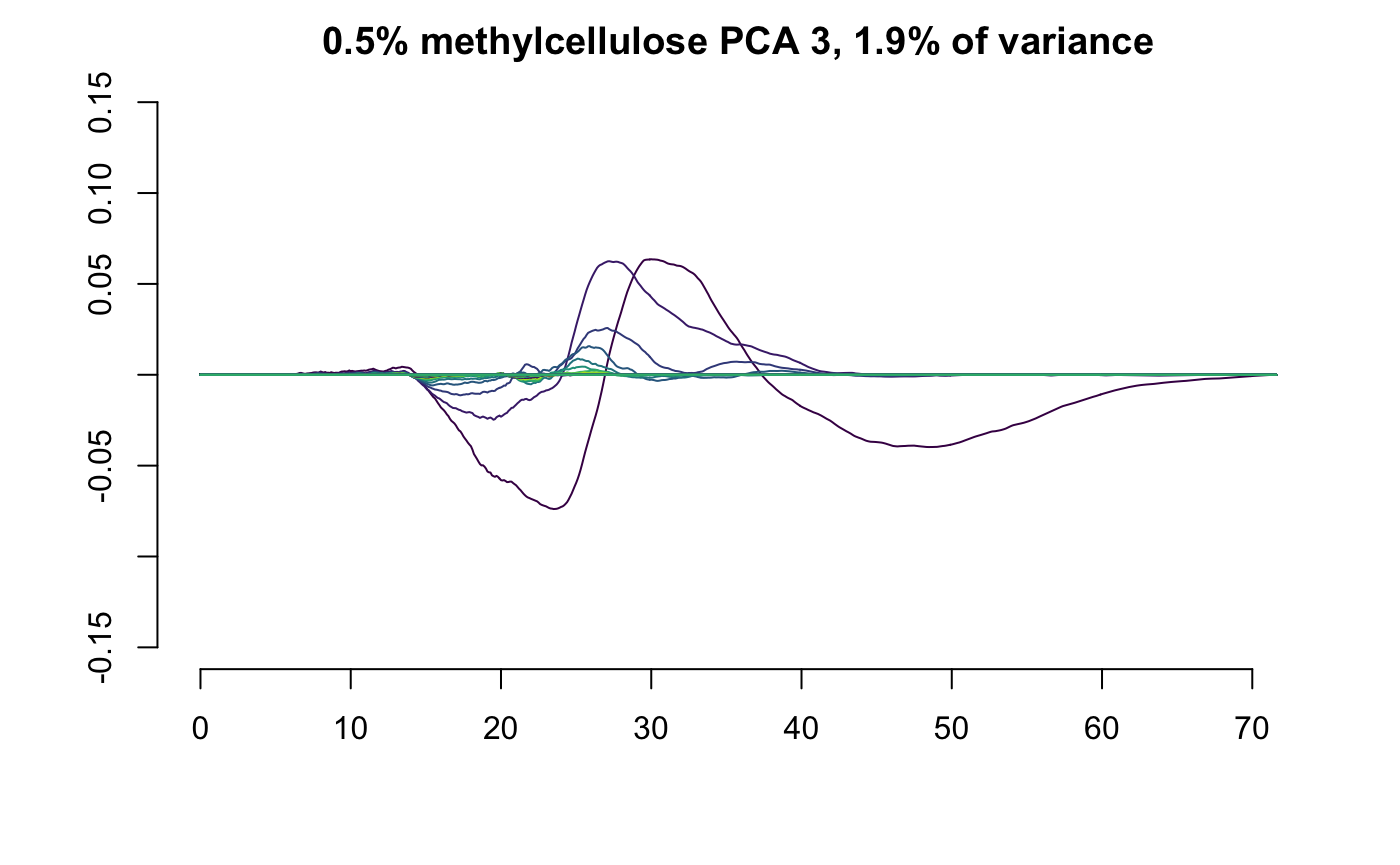}
      &\trimmedgraphic[width=\quarterwidth]{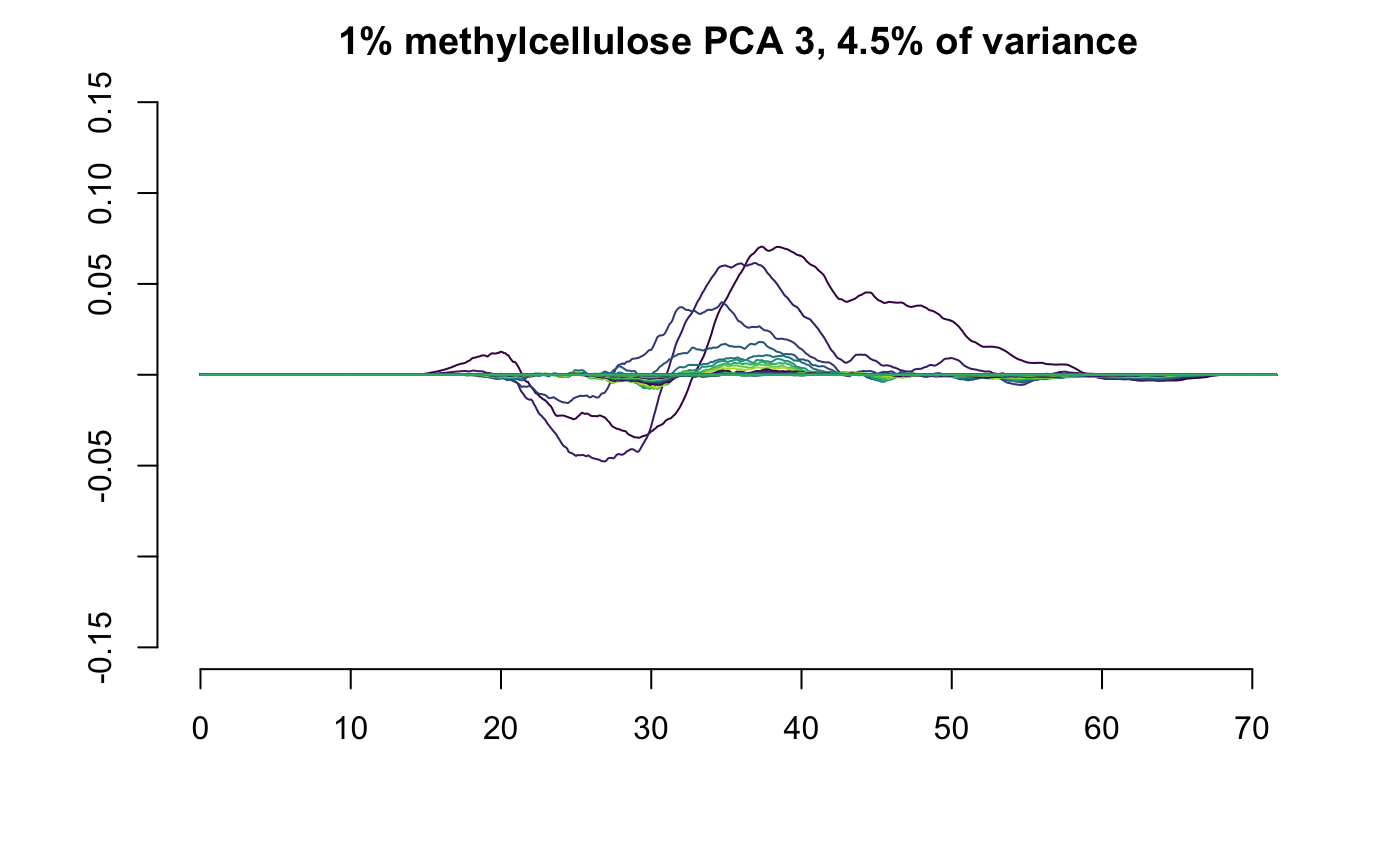}
      &\trimmedgraphic[width=\quarterwidth]{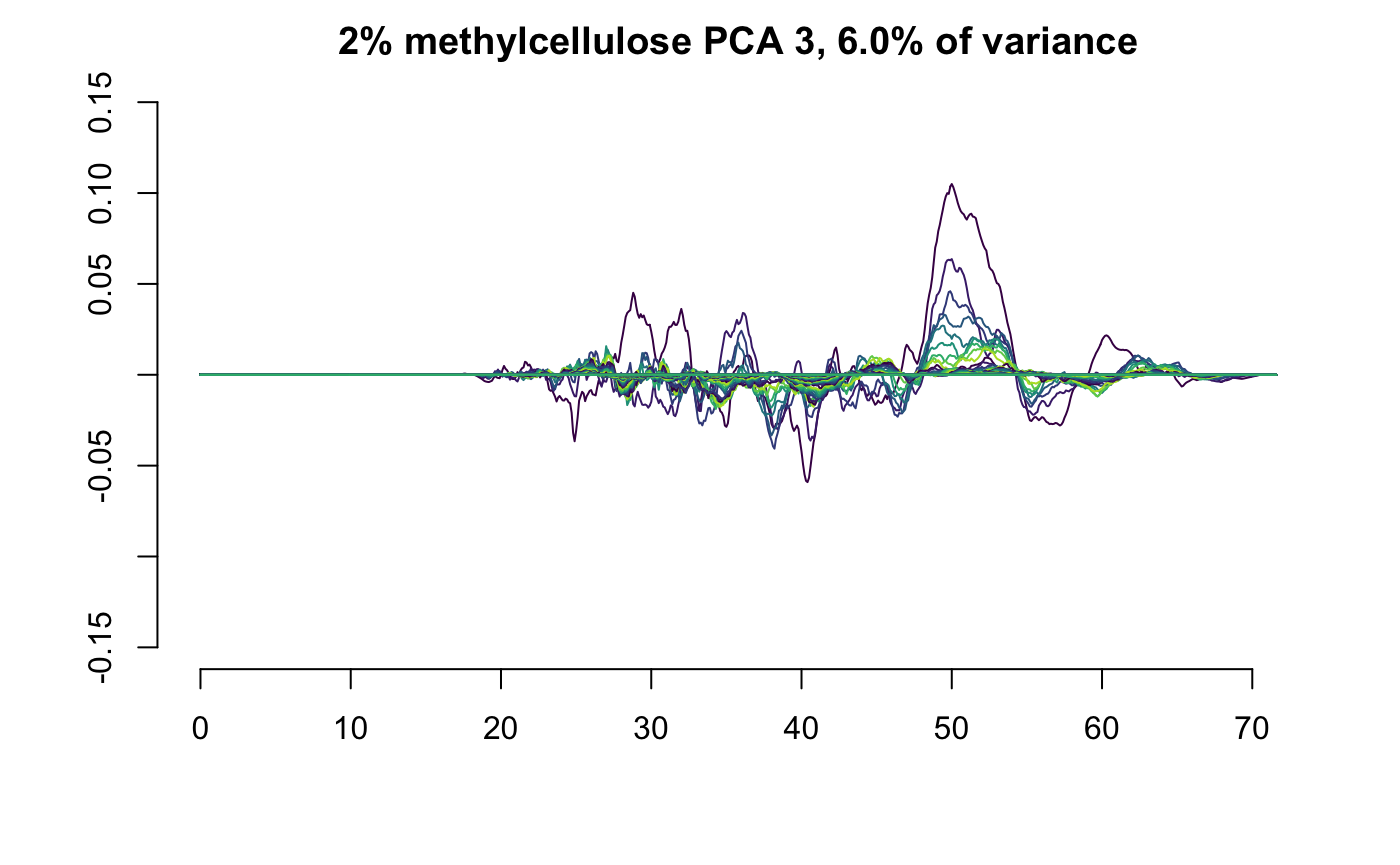}
      &\trimmedgraphic[width=\quarterwidth]{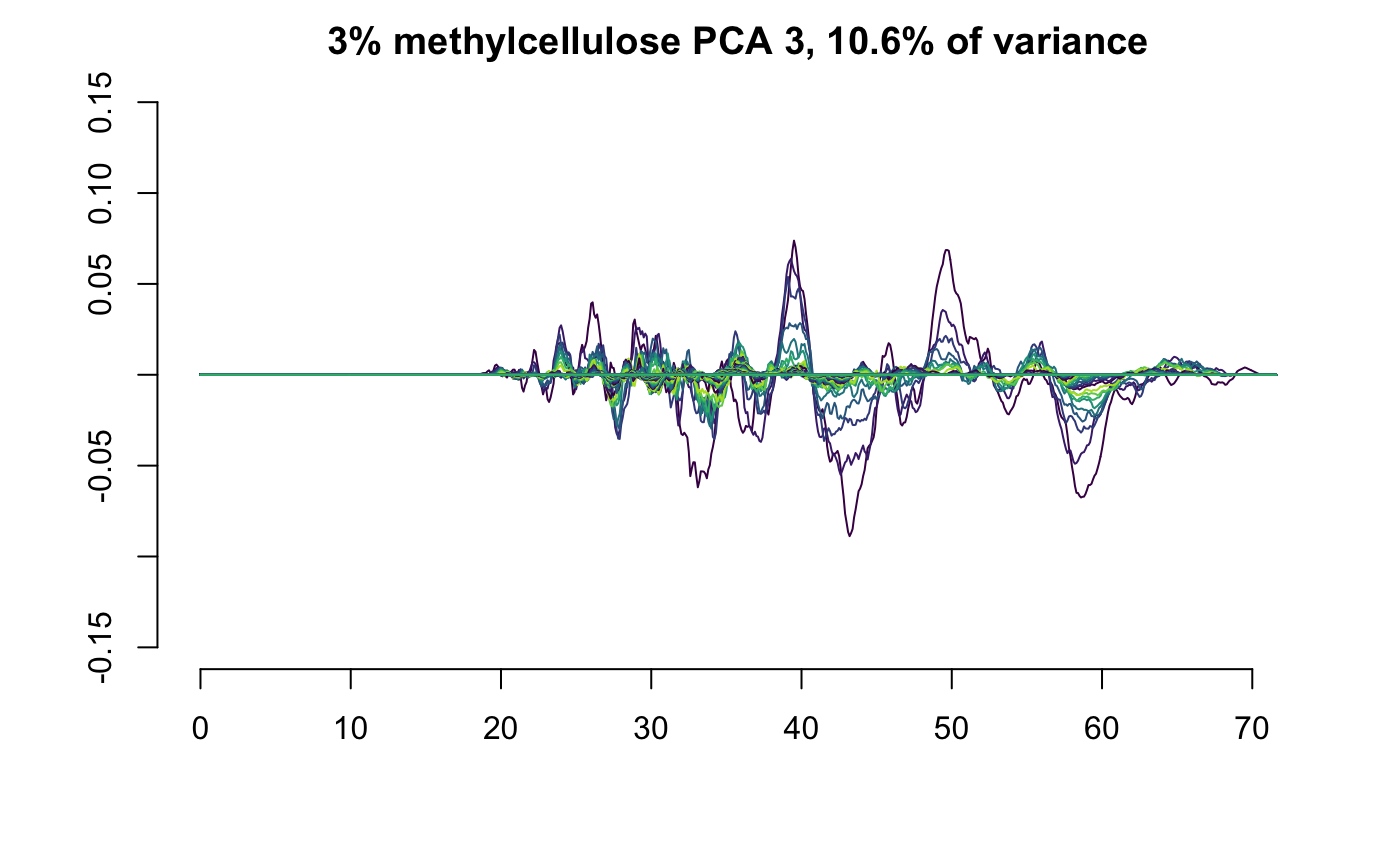}
    \end{tabular}
  \end{center}
  \caption[Average landscapes and PCA for each environmental condition]{PCA on the videos in each class. For each class: (A) video frames showing representative postures. (B) The cumulative variance of the first $n$ principal components. (C-E) The first $3$ PCA eigenvectors, labeled with the percent of the variance described by that component. These results show increasing complexity of shape and behavior as the environment becomes more viscous. }
  \label{fig:pca-by-class}
\end{figure}

In the following analysis we study the complexity of behavior expressed in each class. \Cref{fig:pca-by-class} shows results from using PCA on the videos in each class. The viscosity of the environment is negatively correlated with the percent of variance explained by the first principal component, which suggests that behaviors in low-viscosity environments are simpler than those in high-viscosity environments. Viscosity appears to be correlated with the number of nonzero landscapes, which also suggests that high-viscosity environments allow for more varied behaviors.

We conducted permutation tests between pairs of classes to determine how well average persistence landscapes can distinguish between samples from different classes. The $p$-values for these computations are shown in  \Cref{table:permutationtest-SVMconfusion}(A). The permutation test gives strong evidence of statistical significance, \ie that the topological summaries of samples from each class are significantly different. 

\begin{figure}[H]
  \begin{center}
    \includegraphics[trim = 0 80 50 64,clip, width=10cm]{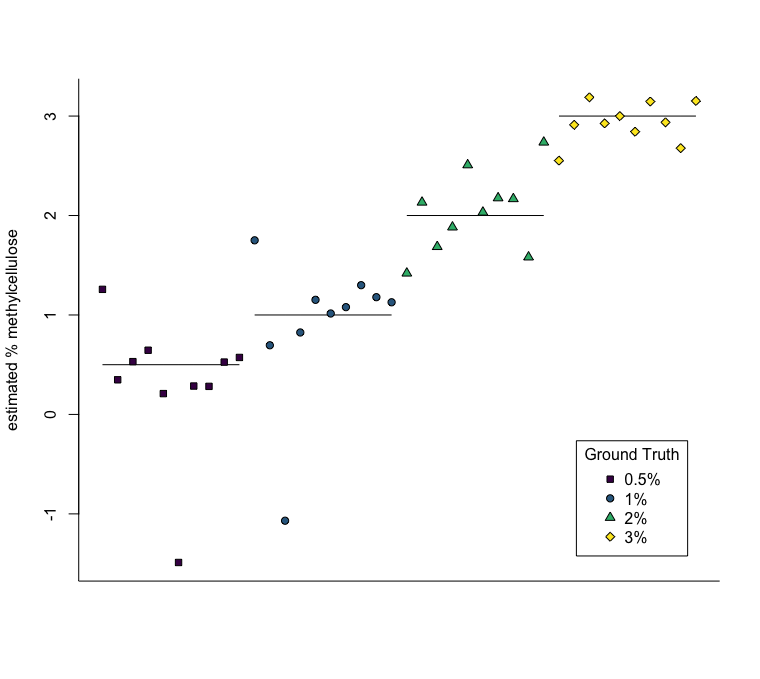}\\
  \end{center}
  \caption[SVR]{SVR estimates of methylcellulose content for each sample. Horizontal lines are at $0.5\%$, $1\%$, $2\%$, and $3\%$ methylcellulose. }
  \label{fig:SVR}
\end{figure}

We then used multiclass support vector machines (SVM) to build a classifier for the samples. The estimated accuracy of the classifier, computed by averaging accuracies across 20 instantiations of the multiclass SVM classifier, was 95.125\%. This indicates that persistent homology is able to produce meaningful, distinguishing features from the \emph{C. elegans} videos. A sample confusion matrix for one instance of SVM is shown in \Cref{table:permutationtest-SVMconfusion}(B). Finally, we used support vector regression to estimate the methylcellulose content in the environment for each sample. The results are plotted in \Cref{fig:SVR}. There are two outliers on this graph which are estimated as having negative methylcellulose content. The two animals in these samples moved much more quickly than their peers, so we believe that the SVR is picking up on the strong negative correlation between viscosity of the environment and speed, and based on these animals' fast speed, assigning a methylcellulose content that is so low that it is negative. 

\begin{table}[H]
  \begin{center}
    \hspace{-11em} (A) \hspace{18em} (B) \\
    \begin{tabular}{c | c c c}
      & $1\%$ & $2\%$ & $3\%$\\
      \hline
      $0.5\%$  & $0.0077$ & $0.0001$ & $0.0000$\\
      $1\%$ &  & $0.0000$ & $0.0000$\\
      $2\%$ &  &  & $0.0000$
    \end{tabular}
    \hspace{5em}
    \begin{tabular}{c | c c c c}
      & $0.5\%$ & $1\%$ & $2\%$ & $3\%$\\
      \hline
      $0.5\%$  & $10$ & $0$ & $0$ & $0$\\
      $1\%$ & $0$ & $9$ & $0$ & $0$\\
      $2\%$ &  $0$ & $1$ & $9$ & $0$ \\
      $3\%$ & $0$ & $0$ & $1$ & $10$
    \end{tabular}
  \end{center}
  \caption[Permutation test/example SVM confusion matrix]{Classification statistics. (A) Permutation test results. (B) Confusion matrix for one instance of SVM with 10-fold cross validation. }
  \label{table:permutationtest-SVMconfusion}
\end{table}

\subsubsection{Cross validation of window length} \label{sec:wl-cross-validation}
As was shown in \Cref{eg:figure-eight}, changes in the window length of the sliding window embedding can affect the conclusions drawn from our computational pipeline. 
In fact, sliding window embeddings have been critiqued for their need of this seemingly arbitrary parameter that can cause significant artifacts when dealing with volatile data \cite{Lindquist2014}. 
Our data is not particularly volatile since it is constrained by physical limitations of \emph{C. elegans}, but we have still justified our choice of window length by conducting cross validation. 

We assembled the results from running our computational pipeline with window lengths of $1$, $10$, $20$, and $30$. We found that different window lengths could be more useful for different tasks. 

The permutation test and multiclass SVM results show that using a window length of $1$ (\ie, not using a sliding window embedding) is the most predictive and that the smaller the window length, the better the accuracy. Window length $1$ had the largest multiclass SVM accuracy at $95.500\%$ compared to accuracies $91.750\%$, $83.375\%$, and $79.125\%$ for window lengths $10$, $20$, and $30$, respectively. The permutation test results were less definitive, with window length $1$ showing slightly more separation between the $0.5\%$ and $1\%$ methylcellulose classes but all window lengths giving strong results. The runtimes of each computation differed significantly, with window length $1$ data running for just over $5$ minutes while window length $30$ data took $90$ minutes. This indicates that for predictive purposes and doing statistics on large data sets, using persistence on as close to the raw data as possible is best, but sliding window embeddings can still be useful if they are desirable for other reasons. 

These statistical results contrast with the ease of visualization and interpretability of the analyses using the different window lengths. 
PCA projections of the cross validation data show much more differentiation between classes when the window length is $10$, $20$, or $30$ compared to when it is $1$. This is presumably because the first two principal components do not explain as much of the variation in the raw data, so separation between classes takes more principal components to describe. Meanwhile, sliding window embeddings consolidate variation in the samples into fewer principal components, so there are fewer principal components to visualize and interpret in terms of the original application. Essentially, sliding window embeddings give a slightly simplified but still predictive representation of the raw data. 

We can also see from \Cref{fig:case-study-persistence} that PCA projections of behavior data are easier to interpret when we use a sliding window embedding than when we look at just the raw time series. Recall that one of the behavioral features from the data --- the pause --- was not detectable over noise when looking at the PCA projection of the original time series but could be identified in both the PCA projection and the persistence landscape of the sliding window embedding. Because we were able to identify the topological feature we were also able to compute a corresponding representative cycle, which in turn allowed an estimate of the locomotion corresponding to the pause behavior to be constructed. Identification and construction of the corresponding locomotion of the pause behavior would have been more difficult using the original time series.

\subsubsection{Validation results} \label{sec:null-model-results}

\begin{table}[H]
  \begin{center}
    \hspace{-11em} (A) \hspace{18em} (B) \\
    \begin{tabular}{c | c c c}
      & $1\%$ & $2\%$ & $3\%$\\
      \hline
      $0.5\%$  & $0.3497$ & $0.0000$ & $0.0000$\\
      $1\%$ &  & $0.0492$ & $0.0480$\\
      $2\%$ &  &  & $0.1026$
    \end{tabular}
    \hspace{5em}
    \begin{tabular}{c | c c c}
      & $1\%$ & $2\%$ & $3\%$\\
      \hline
      $0.5\%$  & $0.0002$ & $0.0000$ & $0.0000$\\
      $1\%$ &  & $0.9961$ & $0.0625$\\
      $2\%$ &  &  & $0.0411$
    \end{tabular}
  \end{center}
  \caption[Permutation tests]{Permutation test results on simpler techniques which use rough measures of each worm's average movement speed and variation in posture. (A) Speed: averages of $2$-norms of the difference between consecutive frames of vector angles. (B) Posture change: standard deviations of vectors angles. }
  \label{table:permutationtest_for_two_techniques}
\end{table}

We conducted permutation tests using $10000$ permutations and applied multiclass SVM using $10$-fold cross validation with cost set $10$ to two techniques based on vector angles which are described in \Cref{sec:null-model}. Results of the permutations test are described in Table \Cref{table:permutationtest_for_two_techniques}. For the method based on averages of $2$-norms of consecutive frames of vector angles the cross validation error was $12.5\%$   and training error was $7.5\%$. For the method based on standard deviations of vector angles the cross validation error was $70\%$  and training error was $35\%$. Though more computationally taxing, persistence and sliding window embeddings produced much more accurate results than either of these simpler techniques. 

Computations for the null model outlined in \Cref{sec:null-model} involve permuting the original time series of each sample and running our sliding window embedding and persistence techniques on that permuted data. Analysis of the null model showed that 
temporal information is necessary for our techniques to give good accuracy differentiating between samples taken in differing viscosity environments. 

For the null model, the average error across 20 iterations of 10-fold cross validation on SVMs was $51.625\%$. The permutation tests showed that we could distinguish the $0.5\%$ methylcellulose class without temporal information, but could not do as well differentiating between the three higher viscosity classes. It seems that the distribution of poses in the $0.5\%$ methylcellulose class was different enough from the other classes to produce noticeably larger topological features that resulted in larger landscapes. This is probably because the lowest viscosity class had more extreme poses and thus the diameter of the space of poses for that class was significantly larger.

\section{Discussion}\label{sec:discussion}

We have demonstrated that persistent homology is a viable technique for studying \emph{C. elegans} behavior and provides useful interpretations and visualizations. Our method consists of constructing sliding window embeddings of time series of piecewise linear \emph{C. elegans} skeletons and using degree $1$ persistent homology to create topological summaries for each patch of each video. These topological summaries, called persistence landscapes, are averaged over patches to produce a single average persistence for each video. These average persistence landscapes are our topological summary statistics and they are the statistics to which we apply further statistical analysis and machine learning, such as principal component analysis, multidimensional scaling, permutation tests, and multiclass support vector machines. As far as we are aware, this is the first application of persistent homology to \emph{C. elegans} behavior data. 

Our analysis showed that persistence is able to detect variability in \emph{C. elegans} behavior data, but also that it can provide interpretable conclusions and useful visualizations. The potential of persistence for interpretability and visualization results is demonstrated in the case study of \Cref{sec:case-study}, where topological features were connected directly to behavioral features and persistence was used to create synthetic behavior data corresponding to stereotyped behaviors such as forward crawling. Our analysis of experimental data shows that persistent homology can detect the variation of behavior induced by changes in the viscosity of the environment. It also suggests that persistence can measure complexity of behavior and that sliding window embeddings with low window lengths can be more predictive while sliding window embeddings with higher window lengths can be more useful for producing clear and interpretable visualizations, including video of synthetic data. 

Persistent homology produces powerful summaries of the ``shape'' of data. However, using persistent homology in a way that is interpretable by experimentalists is a challenge and a topic of current research. We take a step in this direction by using representative cycles  of the most persistent features of a sliding window embedding to produce synthetic videos of characteristic cyclic behaviors. At this time, there does not exist a straightforward way to similarly interpret our composite summaries, the average persistence landscapes. However, there is work in progress towards this goal~\cite{bubenikWagner:heatmap}, and our pipeline would be able to incorporate such advances.

Our analysis has implications for future experimental design. We observed that low-viscosity environments allow for the detection of variation between samples, while high-viscosity environments may allow animals to perform more complex and varying behaviors. Tuning the viscosity of the environment for an experiment or performing experiments in multiple fluid environments with varying viscosities could allow for more easily assessing results regarding variations within populations or variations in behavior. 

An extension to this experiment that could provide more validation for our techniques would be to include samples from two new environmental conditions: buffer, which would correspond to $0\%$ methylcellulose and a lower viscosity than appears in our current data; and agar, which provides a solid surface for the worms to crawl on and surrounding air as opposed to an aqueous environment to be submerged and swim in. We would expect the new buffer class to allow for only fast, simple behaviors in line with the experiments already done, and the agar environment to allow more complex  behaviors in the subjects. 

The method that we have developed for applying topological data analysis to C.\ elegans locomotion data will facilitate the future study of biological phenomena such as aging. In particular, our rich quantitative summary of locomotion suggests that we may be able to measure not just lifespan, but ``healthspan,'' the length of time an individual is healthy and physically capable. Many therapies and medicines for humans and other organisms have a goal of expanding healthspan, and therefore require a detailed measure of the ability to locomote, such as those provided by our methods.

\section{Supplementary Material}\label{sec:supp-material}
 
 Code, data, and sample video output associated with this analysis can be found at \\ \href{https://github.com/althomas/tda-for-worm-behavior}{github.com/althomas/tda-for-worm-behavior}.

\AtNextBibliography{\small}
\printbibliography

\end{document}